\documentclass[11pt, reqno]{amsart}

\usepackage[margin=1.2in,marginparwidth=1.5cm, marginparsep=0.5cm]{geometry}
\pdfoutput=1

\usepackage{amsmath,amssymb,amscd,amsthm,amsxtra, esint}
\usepackage{bm}
\usepackage{microtype}

\usepackage[implicit=true]{hyperref}

\makeatletter
\AtBeginDocument{%
  \@ifundefined{SK@@ref}{}{%
    \ifx\SK@ref\@empty\else
      \renewcommand*\NR@ref@showkeys[1]{\SK@\SK@@ref{#1}}%
    \fi
  }%
}
\makeatother

\usepackage{mathtools} 
\mathtoolsset{showonlyrefs}

\usepackage{upgreek}

\usepackage{booktabs} % nice tables

\usepackage{mathrsfs}

\usepackage{color}

\usepackage{xsavebox}

\usepackage{cases}%%%

\allowdisplaybreaks[2]

\definecolor{gr}{rgb}   {0.,   0.69,   0.23 }
\definecolor{bl}{rgb}   {0.,   0.5,   1. }
\definecolor{mg}{rgb}   {0.85,  0.,    0.85}
\definecolor{yl}{rgb}   {0.8,  0.7,   0.}
\definecolor{or}{rgb}  {0.7,0.2,0.2}

\usepackage{tikz} 
\usepackage{marginnote}
\usepackage{scalerel} % to rescale the paraproduct symbols

\usetikzlibrary{shapes.misc}
\usetikzlibrary{shapes.symbols}
\usetikzlibrary{decorations}
\usetikzlibrary{decorations.markings}

\tikzset{
	dot/.style={circle,fill=black,draw=black,inner sep=0pt,minimum size=0.5mm},
	>=stealth,
	}

\tikzset{
	dot2/.style={circle,fill=black,draw=black,inner sep=0pt,minimum size=0.2mm},
	>=stealth,
	}

\tikzset{
	ddot/.style={circle,fill=black,draw=black,inner sep=0pt,minimum size=0.8mm},
	>=stealth,
	}

\tikzset{decision/.style={ % requires library shapes.geometric
        draw,
        diamond,
        aspect=1.5
    }}

\tikzset{dia2/.style
={diamond,fill=white,draw=black,inner sep=0pt,minimum size=1mm},
	>=stealth,
	}

\tikzset{dia/.style
={star,fill=black,draw=black,inner sep=0pt,minimum size=1mm},
	>=stealth,
	}

\tikzset{dia/.style
={diamond,fill=black,draw=black,inner sep=0pt,minimum size=1.3mm},
	>=stealth,
	}

\makeatletter
\def\DeclareSymbol#1#2#3{\xsavebox{#1}{\tikz[baseline=#2,scale=0.15]{#3}}}
\def\<#1>{\xusebox{#1}}
\makeatother

\newsavebox{\peA}
\newsavebox{\pneA}
\newsavebox{\plA}
\newsavebox{\pgA}
\newsavebox{\pleA}
\newsavebox{\pgeA}
\newsavebox{\pezA}

\savebox{\peA}{\tikz \draw (0,0) node[shape=circle,draw,inner sep=0pt,minimum size=8.5pt] {\scriptsize  $=$};}
\savebox{\pneA}{\tikz \draw (0,0) node[shape=circle,draw,inner sep=0pt,minimum size=8.5pt] {\footnotesize $\neq$};}
\savebox{\plA}{\tikz \draw (0,0) node[shape=circle,draw,inner sep=0pt,minimum size=8.5pt] {\scriptsize $<$};}
\savebox{\pgA}{\tikz \draw (0,0) node[shape=circle,draw,inner sep=0pt,minimum size=8.5pt] {\scriptsize $>$};}
\savebox{\pleA}{\tikz \draw (0,0) node[shape=circle,draw,inner sep=0pt,minimum size=8.5pt] {\scriptsize $\leqslant$};}
\savebox{\pgeA}{\tikz \draw (0,0) node[shape=circle,draw,inner sep=0pt,minimum size=8.5pt] {\scriptsize $\geqslant$};}
\savebox{\pezA}{\tikz \draw (0,0) node[shape=circle,draw,
fill=white, % color = white, 
inner sep=0pt,minimum size=8.5pt]{} ;}

\def \peB{\mathchoice
{\scalebox{.7}{{\usebox{\peA}}}}
{\scalebox{.7}{{\usebox{\peA}}}}
{\scalebox{.7}{{\usebox{\peA}}}}
{}
}

\def \plB{\mathchoice
{\scalebox{.7}{{\usebox{\plA}}}}
{\scalebox{.7}{{\usebox{\plA}}}}
{\scalebox{.7}{{\usebox{\plA}}}}
{}
}
\def \pgB{\mathchoice
{\scalebox{.7}{{\usebox{\pgA}}}}
{\scalebox{.7}{{\usebox{\pgA}}}}
{\scalebox{.7}{{\usebox{\pgA}}}}
{}
}

\def \pgeB{\mathchoice
{\scalebox{.7}{{\usebox{\pgeA}}}}
{\scalebox{.7}{{\usebox{\pgeA}}}}
{\scalebox{.7}{{\usebox{\pgeA}}}}
{}
}

\newcommand{\pe}{\mathbin{{\peB}}}

\newcommand{\pl}{\mathbin{{\plB}}}
\newcommand{\pg}{\mathbin{{\pgB}}}

\newcommand{\pge}{\mathbin{{\pgeB}}}

\usetikzlibrary{decorations.pathmorphing, patterns,shapes}

\DeclareSymbol{dot}{-2.7}{\draw (0,0) node[dot]{};}

\newtheorem{theorem}{Theorem} [section]

\newtheorem{lemma}[theorem]{Lemma}
\newtheorem{proposition}[theorem]{Proposition}
\newtheorem{remark}[theorem]{Remark}

\newtheorem{corollary}[theorem]{Corollary}

\DeclareMathOperator{\HS}{HS}
\DeclareMathOperator{\op}{op}

\newcommand{\1}{\hspace{0.2mm}\text{I}\hspace{0.2mm}}
\newcommand{\II}{\text{I \hspace{-2.8mm} I} }
\newcommand{\noi}{\noindent}
\newcommand{\Z}{\mathbb{Z}}
\newcommand{\R}{\mathbb{R}}

\newcommand{\T}{\mathbb{T}}
\newcommand{\bul}{\bullet}

\newcommand{\wick}[1]{:\!{#1}\!:}

\let\Re=\undefined\DeclareMathOperator*{\Re}{Re}

\newcommand{\TT}{\mathcal{T}}

\let\P= \undefined
\newcommand{\P}{\mathbf{P}}

\newcommand{\E}{\mathbb{E}}

\renewcommand{\L}{\mathcal{L}}

\newcommand{\Kc}{\mathcal{K}}
\newcommand{\F}{\mathcal{F}}

\newcommand{\al}{\alpha}
\newcommand{\be}{\beta}
\newcommand{\dl}{\delta}

\newcommand{\nb}{\nabla}

\newcommand{\Dl}{\Delta}
\newcommand{\eps}{\varepsilon}
\newcommand{\kk}{\kappa}
\newcommand{\g}{\gamma}
\newcommand{\G}{\Gamma}
\newcommand{\ld}{\lambda}

\newcommand{\s}{\sigma}
\newcommand{\Si}{\Sigma}
\newcommand{\ft}{\widehat}

\newcommand{\wt}{\widetilde}
\newcommand{\cj}{\overline}

\newcommand{\dt}{\partial_t}

\newcommand{\ta}{\theta}
\newcommand{\Ta}{\Theta}

\renewcommand{\l}{\ell}
\renewcommand{\o}{\omega}
\renewcommand{\O}{\Omega}

\newcommand{\les}{\lesssim}
\newcommand{\ges}{\gtrsim}

\newcommand{\jb}[1]
{\langle #1 \rangle}

\renewcommand{\b}{\beta}
\newcommand{\ind}{\mathbf 1}

\renewcommand{\S}{\mathcal{S}}

\newcommand{\M}{\mathcal{M}}

\newcommand{\N}{\mathbb{N}}

\renewcommand{\H}{\mathcal{H}}

\newtheorem*{ackno}{Acknowledgements}

\newcommand{\ze}{\zeta}

\newcommand{\Zf}{\mathfrak{Z}}
\newcommand{\zf}{\mathfrak{z}}

\newcommand{\Xs}{\mathsf{X}} 

\newcommand{\ws}{\mathsf{w}}

\newcommand{\Ys}{\mathsf{Y}}

\newcommand{\too}{\longrightarrow}

\newcommand{\I}{\mathcal{I}}

\newcommand{\If}{\mathfrak{I}}

\newcommand{\RR}{\mathcal{R}}
\newcommand{\B}{\mathcal{B}}

\newcommand{\A}{\mathcal{A}}
\newcommand{\D}{\mathcal{D}}

\newcommand{\C}{\mathcal{C}}

\newcommand{\PP}{\mathbb{P}}

\newcommand{\fn}{f_{N, \psi}}

\newcommand{\Law}{\textup{Law}}

\newcommand{\Xc}{\mathcal{X}}

\newcommand{\Ec}{\mathcal{E}} 

\newcommand{\Arg}{\textup{Arg}}

\numberwithin{equation}{section}
\numberwithin{theorem}{section}

\makeatletter
\@namedef{subjclassname@2020}{%
  \textup{2020} Mathematics Subject Classification}
\makeatother

\begin{document}
\baselineskip = 14pt

\title[Beyond variance blowup for NLW]
{Probabilistic well-posedness of dispersive PDEs \\ beyond variance
blowup II:\\ quadratic nonlinear wave equation}

\author[G.~Li, J.~Li, S.~Liu, T.~Oh, and N.~Tzvetkov]
{Guopeng Li, Jiawei Li, Shao Liu, Tadahiro Oh, and Nikolay Tzvetkov}

\address{
Guopeng Li, 
School of Mathematics and Statistics, Beijing Institute of Technology, Beijing 100081, China}

\email{guopeng.li@bit.edu.cn}

\address{
Jiawei Li, 
School of Mathematics\\
The University of Edinburgh\\
and The Maxwell Institute for the Mathematical Sciences\\
James Clerk Maxwell Building\\
The King's Buildings\\
Peter Guthrie Tait Road\\
Edinburgh\\ 
EH9 3FD\\
 United Kingdom}

\email{jiawei.li@ed.ac.uk}

\address{
Shao Liu, 
School of Mathematics and Statistics, Beijing Institute of Technology, Beijing 100081, China\\
and
Mathematical Institute\\
University of Bonn\\
Endenicher Allee 60\\
53115\\
Bonn\\
Germany}

\email{shaoliu@bit.edu.cn}

%
%\address{
%Tadahiro Oh, 
%School of Mathematics and Statistics, Beijing Institute of Technology, Beijing 100081, China\\
%and
%School of Mathematics\\
%The University of Edinburgh\\
%and The Maxwell Institute for the Mathematical Sciences\\
%James Clerk Maxwell Building\\
%The King's Buildings\\
%Peter Guthrie Tait Road\\
%Edinburgh\\
%EH9 3FD\\
% United Kingdom}

%
%	
%
\address{
Tadahiro Oh, School of Mathematics\\
The University of Edinburgh\\
and The Maxwell Institute for the Mathematical Sciences\\
James Clerk Maxwell Building\\
The King's Buildings\\
Peter Guthrie Tait Road\\
Edinburgh\\ 
EH9 3FD\\
 United Kingdom, 
 and 
 School of Mathematics and Statistics, Beijing Institute of Technology,
Beijing 100081, China
}

\email{hiro.oh@ed.ac.uk}

\address{
Nikolay Tzvetkov,
Ecole Normale Sup{\'e}rieure de Lyon\\
UMPA\\
UMR CNRS-ENSL 5669\\
46, all{\'e}e d'Italie\\ 
69364-Lyon Cedex 07\\ 
France}

\email{nikolay.tzvetkov@ens-lyon.fr}

% 35Q35, 


\begin{abstract} 

We study the quadratic nonlinear wave equation (NLW),
posed on the two-dimensional torus,
with rough Gaussian random initial data
beyond the variance blowup threshold
observed by Okamoto and the fourth author (2021).
By introducing a suitable vanishing
multiplicative renormalization constant
on the frequency-truncated Gaussian random initial data,
we prove that the corresponding solutions converge in law
to a solution to the stochastic NLW
with zero initial data,
driven by a Gaussian space-time noise
whose temporal covariance at each spatial Fourier mode
is given by the Bessel function of order~$0$.
We also show that the solution to the limiting equation almost surely blows up
in a finite time, and 
that the convergence above holds 
on every compact time interval before the maximal existence time.
A main new ingredient is the analysis
of the oscillatory covariance
of the second Picard iterates.
We identify the limiting covariance
by averaging the leading oscillatory phase
over the angular variable.
In contrast to our previous work
on the Benjamin-Bona-Mahony equation (2025),
multilinear dispersive smoothing plays an essential role
both in proving tightness of the second Picard iterates
and in controlling the rough nonlinear interactions
via paracontrolled analysis.
The resulting gain of a quarter of a derivative
persists in the limit,
where it is reflected in the decay
of the Bessel covariance.

\end{abstract}

%\date{\today}

\maketitle

\tableofcontents

\section{Introduction}
\label{SEC:1}

\subsection{Overview}
Over the last fifteen years, 
we have seen significant developments
in 
the study of singular stochastic PDEs, broadly 
interpreted
with rough stochastic forcing and\,/\,or random initial data.
The basic overarching 
idea
for constructing solutions to 
singular stochastic PDEs
(see, for example, 
\cite{Hairer0, GIP, Kup, MW2}
in the parabolic case; see 
Subsection \ref{SUBSEC:1.2}
for a discussion in the dispersive case)
can be summarized as follows:

\smallskip

\begin{itemize}
\item
[\bf Step 1.]
By imposing a suitable structure on the unknown, 
identify an associated  enhanced data set, 
consisting of stochastic terms (such as random distributions
and random operators).
Then, use stochastic analysis to 
give a precise meaning to each element in the enhanced data set.

\smallskip

\item
[\bf Step 2.]
Use deterministic analysis to construct 
a deterministic continuous
solution map, sending the enhanced data set to a solution.

\end{itemize}

\smallskip

\noi
This idea dates back to \cite{McK, BO96}
(see also \cite{Lyons98}).
As observed in \cite{Hoshino, Forlano, OO, Liu}, however, 
there are situations
where the
variance of a certain stochastic term diverges  at a threshold 
when tested against an arbitrary non-zero test function\footnote{Or
the variance of each Fourier coefficient of the stochastic term diverges.}
(namely, Step 1 above breaks down), 
while 
 the analytical framework in Step 2 still continues to hold beyond the threshold. 
 We refer to such a phenomenon as {\it variance blowup}.
 See the introduction
in \cite{Hairer} for a further discussion.
See also 
Remark \ref{REM:VB}.

In a recent work~\cite{Hairer}, 
Hairer 
introduced 
a new   renormalization procedure beyond  variance blowup.
By endowing a regularized noise
with 
a vanishing multiplicative renormalization constant, 
he established a ``{\it beyond variance blowup}'' result;
 solutions to the fractional KPZ equation with the renormalized
regularized noise converge {\it in law}
to a solution of the KPZ equation as we remove the regularization.
See \cite{Hairer} for the details.
See also \cite[Subsection~1.3]{LLOT}.
Following Hairer's work \cite{Hairer}, 
the first, second, fourth, and fifth authors~\cite{LLOT}
initiated the study of  
``beyond variance blowup'' 
for dispersive PDEs
in pursuit of  a possible extension of probabilistic well-posedness
theory for  singular random initial data beyond a variance blowup threshold;
see \cite[Section 1]{LLOT} for a further discussion.
In \cite{LLOT}, the authors considered
 the  Benjamin-Bona-Mahony equation (BBM)
 with rough random initial data
  and established 
 a ``beyond variance blowup'' result
 where solutions to BBM with regularized
 random initial data
 (endowed with vanishing multiplicative renormalization constants)
 converge in law to a 
 solution of the {\it stochastic} BBM forced by the derivative of a spatial white noise
 (with zero initial data),
 exhibiting a new phenomenon where random initial data
 turns into a stochastic forcing.

 In this paper, we continue our study on 
 dispersive PDEs with random initial data beyond variance blowup.
 In particular, we consider the quadratic nonlinear wave equation (NLW)
 on the two-dimensional torus $\T^2 = (\R/2 \pi \Z)^2$ (see \eqref{NLW1} below), 
 for which variance blowup was observed in \cite{OO};
 see Subsection \ref{SUBSEC:2.1}.
 By following the methodology developed in \cite{Hairer, LLOT}, 
we establish a ``beyond variance blowup'' result
for the quadratic NLW \eqref{NLW1};  
solutions to the quadratic 
NLW with the renormalized frequency-truncated random initial data
(see \eqref{NLW1a})
converge in law to a solution to a stochastic NLW
forced by a space-time Gaussian noise.
See Theorem~\ref{THM:1} for a precise statement.
We note that, 
in \cite{LLOT}, 
we in fact considered BBM   
  as a toy model for the quadratic NLW~\eqref{NLW1}
in the sense
that, for both equations,  the nonlinearity is quadratic, enjoying
an explicit smoothing of one degree in the Duhamel formulation.
There are, however, several crucial differences.

\smallskip

\begin{itemize}
\item
In contrast to BBM studied in \cite{LLOT}, 
the quadratic NLW enjoys
multilinear dispersive smoothing;
 the second Picard iterate\footnote{Strictly speaking, 
it is the second Picard iterate minus the linear solution.
For simplicity, however, we refer to it as the second Picard  iterate.
The same convention will be used throughout the paper.}  
$Z_N$ in   \eqref{ZN}
for the (Wick renormalized)
NLW with the renormalized (frequency-truncated) random initial data
(see \eqref{NLW1a})
gains $\frac 14$ derivative thanks to multilinear dispersion;
see Proposition~\ref{PROP:rand1}.
This second Picard iterate turns into 
the stochastic convolution $Z$ (see \eqref{Z1} and \eqref{NLW3a})
with 
a Gaussian space-time  noise $\ze$ in the limit,
where 
this multilinear smoothing  of $\frac14$ degree is encoded 
(at the linear level)
 in the decay
of the Bessel covariance;
see Theorem \ref{THM:ZNZ}.  See also \eqref{XI1} and \eqref{nxib}.

\smallskip
\item 
In \cite{LLOT}, the limiting equation was shown to be almost surely globally well-posed.
For our NLW problem, the situation is quite different;
we prove that there is almost surely 
a singularity formation (= finite time blowup)
for the limiting stochastic NLW  \eqref{NLW1b};
see Proposition \ref{PROP:B1}.
Our convergence result (Theorem \ref{THM:1})
is optimal and ``global'' in the sense that 
convergence  holds 
on every compact time interval before the (almost surely finite) maximal existence time.

\end{itemize}

\smallskip

\noi
We also point out that, 
in controlling rough nonlinear interactions, 
we employ  paracontrolled analysis, 
developed in \cite{GKO2, OOT1, OOT2}
in the context of wave equations.
See Subsection \ref{SUBSEC:IPa}.

\subsection{Main result}
%{Quadratic nonlinear wave equation}
\label{SUBSEC:1.1}

We consider the following
quadratic NLW
on $\T^2 $:
\begin{align}
\begin{cases}
\dt^2 u + (1- \Dl) u = u^2 \\
(u, \dt u)|_{t = 0} = (u_0, u_1),
\end{cases}
\label{NLW1}
\end{align}

\noi
where $u$ is a real-valued unknown.
Here, 
$(u_0, u_1)$ denotes
 the Gaussian random initial data
of the form:\footnote{By convention, we endow $\T^2$ with the normalized Lebesgue 
measure $dx_{\T^2} = (2\pi)^{-2} dx$ such that 
we do not need to carry factors involving $2\pi$.}
\begin{align}
(u_0, u_1)   
 = (u_0^\o(\al), u_1^\o(\al)) = 
\bigg( \sum_{n\in \Z^2} \frac{g_n(\o) }{\jb{n}^{1-\al}} e_n,
\sum_{n\in \Z^2} \jb{n}^\al h_n(\o)  e_n \bigg)
\label{data1}
\end{align}

\noi
for  $\al \in \R$, 
where  $\jb{n} = \sqrt{1 + |n|^2}$,
$e_n (x) = e^{i n\cdot x}$, and
$\{g_n, h_n \}_{n\in \Z^2}$ is a family of independent
complex-valued standard Gaussian random variables on 
some probability space $(\O, \F, \PP)$ 
conditioned that $g_{-n} = \cj{g_n}$
and $h_{-n} = \cj{h_n}$, $n \in \Z^2$.  
Given $s \in \R$, let $\mu_s$ be  the (massive) fractional Gaussian field on $\T^2$
with the formal density:
\begin{align}
d \mu_s = Z_s^{-1} e^{-\frac 12 \|u \|_{H^s}^2} du. 
\label{mu1}
\end{align}

\noi
Then, we have $\Law ((u_0, u_1)) = \mu_{1-\al}\otimes \mu_{-\al}$
for the Gaussian random initial data $(u_0, u_1)$ in~\eqref{data1}.
Here, $\Law(X)$ denotes the law of a random variable $X$.

Given $s \in \R$ and $1 \le p \le \infty$, 
set 
\[\vec W^{s, p} (\T^2)
=  W^{s, p} (\T^2) \times  W^{s-1, p} (\T^2) 
\qquad \text{and}\qquad 
\vec H^{s}(\T^2) = \vec W^{s, 2} (\T^2).
 \]

\noi
Then, 
it is easy to see that the random initial data $(u_0, u_1)$
in \eqref{data1} belongs
almost surely to $\vec W^{s, p} (\T^2) \setminus  \vec W^{-\al, p} (\T^2)$
for any $s < - \al$ and $1 \le p \le \infty$.
Recall from 
 \cite{LS, GKO}
that \eqref{NLW1} is locally well-posed
in $\vec H^s(\T^2)$ for $s > 0$, 
which trivially
yields
almost sure local well-posedness
of~\eqref{NLW1}
with the random 
initial data $(u_0, u_1)$
in \eqref{data1} for $\al < 0$.
On the other hand, it was shown in~\cite{FO}
that 
\eqref{NLW1} is ill-posed
in $\vec H^s(\T^2)$
for $s < 0$ 
and thus it is of great interest to 
study~\eqref{NLW1} with random initial data
(in particular, of the form \eqref{data1})
to seek  a possible extension of well-posedness theory in $\vec H^s(\T^2)$
for $s < 0$.

%When $\al = 0$, 
The case  $\al = 0$ is of particular importance
since, in this case, the Gaussian random functions
$u_0$ and $u_1$ in \eqref{data1} correspond
to the massive Gaussian free field
and the white noise on $\T^2$, respectively, 
and 
$\Law(u_0, u_1) = \mu_1 \otimes \mu_0$
serves as the base Gaussian measure
for the Gibbs measure associated with \eqref{NLW1};
see \cite{BO99, OTh2, OST}.
Note that for $\al \ge 0$, 
a solution $u$ is
expected to be merely a distribution in the spatial variable
and thus a renormalization is required
in giving a proper meaning to the nonlinearity $u^2$.
For simplicity of presentation, however, we ignore the issue
of renormalization for NLW \eqref{NLW1} in this introduction;
see Subsection~\ref{SUBSEC:2.1},
in particular~\eqref{NLW1x}, for the Wick renormalization in this context.
When $\al = 0$,
it was shown in 
\cite{OTh2, GKO, GKOT}\footnote{In \cite{OTh2},
the authors only treated the odd power nonlinearity $u^{2k+1}$, $k \in \N$, 
 but 
similar analysis applies to the quadratic nonlinearity, 
giving almost sure local well-posedness
of \eqref{NLW1} when $\al = 0$.
In 
\cite{GKO, GKOT}, the authors considered the stochastic (damped) NLW
but their analysis equally applies to NLW \eqref{NLW1} with the Gaussian random data in~\eqref{data1}
for $\al = 0$.
} 
that~\eqref{NLW1} 
is almost surely locally well-posed
with the Gaussian random initial data $(u_0, u_1)$
in \eqref{data1}.

In \cite{OO}, Okamoto and the fourth author investigated
the case $\al > 0$
and proved the following results:\footnote{In \cite{OO},
the authors in fact studied the  stochastic NLW (SNLW) on $\T^2$
(see \eqref{SNLW2} below)
but their analysis equally applies to 
the deterministic NLW~\eqref{NLW1} with the Gaussian random data $(u_0, u_1)$ in~\eqref{data1}.}

\smallskip
\begin{itemize}
\item
for $0 < \al < \frac 12$:
almost sure local well-posedness of 
%(the renormalized version of)~
\eqref{NLW1}
with the Gaussian random initial data in \eqref{data1}, 

\smallskip

\item 
for  $\al \ge \frac 12$:
 variance blowup
 for the second Picard iterate
   associated with~\eqref{NLW1} with the random initial data in \eqref{data1};
see 
\eqref{div1} and 
\eqref{div2}.
See also \cite{Deya, Deya2}.

\end{itemize}

\smallskip

\noi
See Subsection \ref{SUBSEC:2.1}
for further details on these results.
We note that the variance blowup
result for $\al \ge \frac 12$ in \cite{OO} implies
probabilistic ill-posedness 
in the mild sense
of \eqref{NLW1} 
with the Gaussian random initial data \eqref{data1}
for $\al \ge \frac 12$;  see \cite[Section~3]{OTz}
for a further discussion.

Our main goal in this paper is 
to establish a ``beyond variance blowup'' result
for~\eqref{NLW1} when $\al \ge \frac 12$, 
following the work \cite{Hairer, LLOT}.
More precisely, 
given $N \in \N$, we consider 
the following  NLW
with the renormalized (frequency-truncated) random initial data:
\begin{align}
\begin{cases}
\dt^2 u_N  + (1 -  \Dl)  u_N = \, :\!u_N^2\! :\\
(u_N, \dt u_N) |_{t = 0}
  = (C_{\al, N} \pi_N u_0, C_{\al, N} \pi_N u_1),
\end{cases}
\label{NLW1a}
\end{align}

\noi
where 
$(u_0, u_1) = (u_0(\al), u_1(\al))$ is the Gaussian random initial data in \eqref{data1}
and $\pi_N$ denotes the Dirichlet frequency projector
onto (spatial) frequencies $\{|n|\le N\}$; see \eqref{Diri1}.
Here, 
\[:\!u_N^2\! : \, = u_N^2 - \s_N\]

\noi
 denotes the standard Wick power, 
where $\s_N$ is as in  \eqref{sigN}, 
and 
$C_{\al, N} > 0$ denotes a %(suitably chosen) 
vanishing  
multiplicative renormalization constant, defined by 
\begin{align}
C_{\al, N}
= \bigg(\sum_{\substack{n \in \Z^2\\|n| \le N}} \frac {1}{\jb{n}^{4-4\al}}\bigg)^{-\frac 14}
\sim_\al
\begin{cases}
(\log \jb{N})^{-\frac 14} & \text{for } \al = \frac 12,  \\
N^{\frac 12 -\al}
  &  \text{for }   \al>\frac 12. 
\end{cases} 
\label{CN}
\end{align}

\noi
Note that 
\begin{align}
C_{\al, N} \too 0\quad \text{as}\quad N \too \infty.
\label{CN1}
\end{align}

 \noi
On top of the standard Wick renormalization on the quadratic nonlinearity, 
we applied a renormalization on the (frequency-truncated) rough random initial data
in \eqref{NLW1a}.

Fix $\al \ge \frac 12$.
Given $N \in \N$, 
let $z_N$ be the random linear solution for~\eqref{NLW1a}:
\begin{align}
\begin{split}
z_N (t)
& = 
W(t)\big(C_{\al, N} \pi_N u_0, C_{\al, N} \pi_N u_1\big), %\\
%& = \cos (t \jb{\nb}) C_{\al, N} \pi_N u_0
%+ \frac {\sin (t \jb{\nb})}{\jb{\nb}}  C_{\al, N} \pi_N u_1, 
\end{split}
\label{zN}
\end{align}

\noi
where
$W(t)$ is  the linear wave propagator given by 
\begin{align}
W(t)(f, g)
 = \cos (t \jb{\nb}) f
+ \frac {\sin (t \jb{\nb})}{\jb{\nb}}  g.
\label{Wx}
\end{align}

\noi
In the following, 
we restrict our attention to spatial regularities.
Moreover, in referring to regularity
of an $N$-dependent term, 
regularity means uniform (in $N$) regularity.
A standard computation shows 
that  
 $z_N$ has regularity\footnote{We
use the notation $a-$ 
%(and $a+$) 
to denote $a- \eps$ %(and $a+ \eps$, respectively)
for arbitrarily small $\eps > 0$.
When this notation appears in an estimate,  
the implicit constant is allowed to 
depend on $\eps> 0$ (and it usually diverges as $\eps \to 0$).} 
 $-\frac12-$.
Moreover,  
$z_N$
converges to $0$ as $N \to \infty$; 
see 
 Lemma \ref{LEM:randln}.

Define $\s_N$ by setting
\begin{align}
\s_N = \E \big[ z_{N}^2 (t, x)\big]
= C_{\al,N}^2
\sum_{|n| \le N} \frac {1}{\jb{n}^{2 - 2\al}}
\sim_\al C_{\al,N}^2  N^{2 \al}.
\label{sigN}
\end{align}

\noi
Note that $\s_N$ is independent of $(t, x) \in \R\times \T^2.$
Moreover, from \eqref{CN}, 
we see that $\s_N \to \infty$ as $N \to \infty$.
We then define  the second Picard iterate $Z_N$ by 
\begin{align}
 Z_N = \I(\, \wick{z_N^2} \, )
 = \int_{0}^{t} \frac {\sin {((t - t') \jb{\nabla})}}{\jb{\nabla}}
 \wick{z_N^2}  \!(t')\,  d t', 
\label{ZN}
\end{align}

\noi
where $\I$ denotes the wave Duhamel integral operator defined by 
\begin{align}
\I (F) (t) =  \int^t_0
\frac {\sin {((t - t') \jb{\nabla})}}{\jb{\nabla}}
F (t') dt'
\label{Duhamel}
\end{align} 

\noi
and 
$\wick{z_N^2}$ denotes the Wick power given by
\begin{align*}
\wick{z_N^2}  \, =  z_N^2  -  \s_N.
\end{align*}

\noi
It is easy to see that the Wick power $\wick{z_N^2}$ has regularity
$-1-=2\big(-\frac12-\big)$. 
If one proceeds with ``parabolic thinking'', one might  expect that $Z_N$ has
regularity $0-=2\big(-\frac12-\big)+1$, gaining
 one derivative  from the Duhamel integral operator $\I$  in \eqref{Duhamel}.
As observed in \cite{GKO2, OO}, 
 however, 
there is multilinear dispersive smoothing 
of $\frac 14$ degree
and $Z_N$  has regularity $\frac 14 -$; 
see Proposition~\ref{PROP:rand1}.

As in \cite{Hairer, LLOT}, 
the first key step for proving a ``beyond variance blowup'' result for NLW~\eqref{NLW1}
is to investigate convergence properties  of $Z_N$.
For the current NLW problem, 
the multilinear dispersive smoothing
for $Z_N$ mentioned above makes
it more difficult to identify a limiting process
(as
compared to  \cite{Hairer, LLOT}).
For simplicity of presentation,  we restrict our attention to positive times
in the remaining part of this paper.

Define  a Gaussian space-time noise 
 $\ze$ 
on $\R_+ \times \T^2 =[0, \infty)\times \T^2$ by 
\begin{align}
\ze(t) = \sum_{n\in \Z^2} \ze_n(t) e_n
\label{XI1}
\end{align}
\noi

\noi
whose Fourier coefficients $\{\ze_n\}_{n \in \Z^2}$
are independent complex-valued
centered stationary Gaussian processes
conditioned that
$\ze_{-n} = \cj{\ze_n}$, $n \in \Z^2$,
with 
the space-time covariance  given by
\begin{align}
\E \big[\ze_n(t_1) \cj{\ze_m(t_2)}\big]
= \ind_{n = m}\cdot J_0\big((t_1-t_2)|n|\big)
\label{nxib}
\end{align} 

\noi
for $t_1, t_2 \in \R_+$ and 
$n, m \in \Z^2$, 
where $J_0$ is the Bessel function of order $0$
as in~\eqref{Bessel1}.
Using  \cite[Lemma~2.6]{LLOT} with~\eqref{nxib}, 
it is easy to see that 
$\ze \in C(\R_+; H^{-1-\eps}(\T^2))$, almost surely.
Then, define 
a process $Z$ by setting
\begin{align}
Z(t)
= 
\int^t_0 \frac{\sin((t-t')\jb{\nb})}{\jb{\nb}}
\ze(t') dt', 
\label{Z1}
\end{align}

\noi
satisfying  the following  stochastic linear wave equation:
\begin{align}
\begin{split}
\begin{cases}
\dt^2 Z + (1 - \Dl) Z =  \ze \\
(Z, \dt Z)|_{t = 0} = (0, 0).
\end{cases}
\end{split}
\label{NLW3a}
\end{align}

\noi
In 
Theorem \ref{THM:ZNZ}, 
we show that,
as $N \to \infty$, 
the second Picard iterate  $Z_N$ in \eqref{ZN}
converges in law to the Gaussian process 
 $Z$  in \eqref{Z1} (with respect to a suitable topology).
Here,  the 
Bessel covariance structure~\eqref{nxib}
of the driving noise $\ze$
for the limiting process $Z$ in~\eqref{Z1}
encodes 
the multilinear dispersive smoothing of $\frac 14$-degree
for the second Picard iterate $Z_N$.
Indeed, recalling from \cite[Appendix B.8]{Gra} that 
$|J_0(r)| \les r^{-\frac 12}$
for $r \gg 1$, 
 the 
Bessel covariance structure~\eqref{nxib}
and time integration in~\eqref{Z1}
yield
extra $\frac 14$ regularity
for $Z$
(in addition to the one degree of smoothing
from the wave Duhamel integral operator);
see \eqref{Z3} below.
See also Lemma \ref{LEM:randZ}.
Immediately after Theorem \ref{THM:ZNZ}, 
we provide a brief discussion
on ideas of the proof of 
the convergence in law of $Z_N$ to $Z$,
where we exploit 
angular averaging in carrying out summations
(see Lemma \ref{LEM:angavg}).

Once we establish convergence in law of $Z_N$ 
to $Z$  (and also of
the random linear solution $z_N$ in~\eqref{zN} to $0$; see 
Lemma \ref{LEM:randln}), 
we can proceed as in \cite{LLOT}
via the  second order expansion
 for the solution $u_N$ to \eqref{NLW1a}:
\begin{align}
u_N = z_N + Z_N + v_N
\label{exp1}
\end{align}

\noi
and identify a limiting equation.
As $N \to \infty$, the second order expansion \eqref{exp1}
formally converges to the following 
 first order expansion:
\begin{align}
u = Z + v.
\label{exp1a}
\end{align}

\noi
By formally taking $N \to \infty$ 
in the equation satisfied by $v_N$ (see \eqref{NLW3})
and using the relation~\eqref{exp1a}, 
we obtain 
the following stochastic NLW\textup{:}
\begin{align} 
\begin{cases}
\dt^2 u + (1 - \Dl) u = u^2 +  \ze \\ 
(u, \dt u)|_{t = 0} = (0, 0)
\end{cases} 
\label{NLW1b}
\end{align}

\noi
as the limiting equation.  Here, $\ze$ denotes
the Gaussian space-time noise defined in \eqref{XI1}
and~\eqref{nxib}.
We note that in view  of the extra $\frac 14$
gain of regularity 
thanks to the Bessel covariance~\eqref{nxib}, 
  a solution $u$
to \eqref{NLW1b}
is expected to have spatial regularity $\frac 14 -$
and, thus,
 there is no renormalization on the nonlinearity $u^2$ in
 the limiting equation~\eqref{NLW1b}.

In the case of KPZ and BBM studied
in \cite{Hairer} and \cite{LLOT}, 
the limiting equations were globally well-posed.
For the current NLW problem, 
however, 
the limiting equation \eqref{NLW1b}
is {\it not} globally well-posed.
This is due to the (non-defocusing) quadratic nonlinearity, 
which yields ODE-type finite time blowup solutions
for~\eqref{NLW1} in the classical setting.

%while  \eqref{NLW1b} is locally well-posed, 
%its solution almost surely blows up in finite time.

\begin{proposition}
\label{PROP:B1}

Fix 
  $0 < s < \frac{1}{4}$.

\smallskip

\noi
\textup{(i) (local well-posedness).}
The stochastic NLW \eqref{NLW1b}
is locally well-posed
such that the local-in-time solution
 $(u, \dt u)$ belongs to the class
 $C([0, \tau ); \vec H^s(\T^2))$, 
 where $\tau$ denotes the almost surely positive local existence time.

\smallskip

\noi
\textup{(ii) (almost sure finite time blowup).}
Let $T_{\max}$ denote
the almost surely positive maximal existence time
for the solution 
 $(u, \dt u) $ to~\eqref{NLW1b}
 such that 
  $(u, \dt u)  \in C([0, T_{\max}); \vec H^s(\T^2))$.
Then, we have
\begin{align}
 \PP (  T_{\max} < \infty  ) = 1.
\label{F1}
\end{align}

\noi
Moreover,  we have
\begin{align}
\lim_{t \to T_{\max}-}
\|(u(t), \dt u(t))\|_{\vec H^s} = \infty, \quad \text{almost surely.}
\label{F2}
\end{align}
\end{proposition}

Given the regularity properties
of the limiting process $Z$ in \eqref{Z1}
(see Lemma \ref{LEM:randZ}), 
local well-posedness of \eqref{NLW1b} follows
from a standard contraction argument.
See also  Remarks~\ref{REM:uniq1} and \ref{REM:uniq2}.
We prove Proposition~\ref{PROP:B1}\,(ii)
on almost sure finite time blowup
by a contradiction argument, using  the test function method
of Zhang \cite{Zha1, Zha2, Oh4, OOP} and exploiting
the non-negativity of the nonlinearity $u^2$
and the stationarity (in time) of the noise $\ze$.
See Section~\ref{SEC:B} for the proof of 
Proposition~\ref{PROP:B1}.

\begin{remark}\label{REM:LWP1}
\rm 

Let $N \in \N$. Then, \eqref{NLW1a}
is almost surely locally well-posed
since the random initial data is smooth.
An adaptation of small data global well-posedness
in $\vec H^1(\T^2)$
for NLW~\eqref{NLW1}
shows that \eqref{NLW1a}
is globally well-posed with positive probability.
On the other hand, by noting that 
the zeroth Fourier mode $\ft u_N(0)$
is increasing and strictly convex
when $\ft u_N(0)|_{t = 0}$ is sufficiently large
and $\dt \ft u_N(0)|_{t = 0} > 0$, 
an ODE argument
 shows that
the corresponding solution to \eqref{NLW1a}
blows up in finite time.
Namely, the solution to \eqref{NLW1a}
blows up in finite time with positive probability.
\end{remark}

We are now ready to state
our main theorem 
on a ``beyond variance blowup'' 
convergence result
for the quadratic NLW \eqref{NLW1}.

\begin{theorem}\label{THM:1}

Let $\al \ge \frac 12$.
Given $N \in \N$, 
let $\vec u_N = (u_N, \dt u_N)$ 
be 
the \textup{(}smooth\textup{)} solution to~\eqref{NLW1a}
with the renormalized initial data
$\big(C_{\al, N} \pi_N u_0, C_{\al, N} \pi_N u_1\big)$.
Let $T$ be a  random time such that  $0 < T < T_{\max}$
almost surely, 
where
$T_{\max}$ is 
the almost surely positive
 and finite
maximal existence time for \eqref{NLW1b}
from Proposition \ref{PROP:B1}\,(ii).
Then,  
there exists a sequence $\{T_N\}_{N \in \N}$
of almost surely positive random times
such that

\smallskip

\begin{itemize}
\item
for each $N \in \N$, the solution 
$\vec u_N$ to \eqref{NLW1a} almost surely exists on $[0, T_N]$, and

\smallskip
\item as $N \to \infty$, 
the pair $\big(T_N, \vec u_N(\,\cdot\, \wedge T_N)\big)$
converges in law to the pair
$\big(T, \vec u(\,\cdot\, \wedge T)\big)$
in\footnote{In the remaining part of this paper, we use $\eps > 0$
to denote an (arbitrarily) small constant.} $\R\times C(\R_+; \vec H^{-\frac 12 - \eps}(\T^2))$, 
where 
 $\vec u = (u, \dt u)$ is the solution to 
the stochastic NLW~\eqref{NLW1b}
and
$C(\R_+; \vec H^{-\frac 12 - \eps}(\T^2))$
is endowed with the compact open topology in time.

\end{itemize}

\smallskip

\noi
Here, $a \wedge b = \min (a, b)$.

\end{theorem}

As in the BBM case studied in \cite{LLOT}, 
 the effect of the random initial data
in \eqref{NLW1a}
turns into the  stochastic forcing $\ze$ in the limiting equation~\eqref{NLW1b}.
While the noise in the limiting equation
was independent of time in the BBM case, 
our noise $\ze$ 
in the limiting equation~\eqref{NLW1b}
is 
time-dependent, reflecting
the multilinear dispersive nature of  \eqref{NLW1a} in the limit as $N \to \infty$. % (with finite $N \in \N$).
Such embedding of multilinear smoothing
in the limiting {\it linear} process seems
to be a new phenomenon in the study of dispersive equations.

For the BBM problem studied in \cite{LLOT}, 
both  equations in \cite[(1.18) and (1.29)]{LLOT}
are (almost surely) globally well-posed. 
On the other hand, 
for the current NLW problem,
the solution $u$
to 
the limiting equation~\eqref{NLW1b}
almost surely blows up in finite time (Proposition \ref{PROP:B1});
see also Remark \ref{REM:LWP1}.
For this reason, 
in Theorem \ref{THM:1}, 
we considered the stopped processes
$\vec u_N(\,\cdot\, \wedge T_N)$
and $\vec u(\,\cdot\, \wedge T)$, 
which allows
us to work 
in one fixed space
$C(\R_+; \vec H^{-\frac 12 - \eps}(\T^2))$
(rather than a random space 
$C([0, T]; \vec H^{-\frac 12 - \eps}(\T^2))$
with a random time
  $0 < T < T_{\max}$).

Let us briefly describe our strategy 
for proving Theorem \ref{THM:1}, 
pointing out 
 similarities and differences with the BBM case.
As in the BBM case, 
the first crucial step  
is  to establish convergence in law of
the second Picard iterate $Z_N$
for \eqref{NLW1a} (see~\eqref{ZN})
to the limiting process $Z$ in \eqref{Z1}; see Theorem \ref{THM:ZNZ}.
Contrary to the BBM case, 
the bilinear structure of $Z_N$ in~\eqref{ZN} induces
bilinear dispersive  smoothing 
of $\frac 14$-degree 
(as in \cite[Proposition~1.6]{GKO2}
and 
\cite[Proposition~1.4]{OO});
see 
Proposition~\ref{PROP:rand1}.
Once we establish convergence in law of $Z_N$ in~\eqref{ZN}
to $Z$ in \eqref{NLW3a} (and also of
the random linear solution $z_N$ in~\eqref{zN} to $0$), 
we invoke the Skorokhod representation theorem (Lemma~\ref{LEM:Sk})\footnote{Strictly speaking, 
upon the application of the Skorokhod representation theorem, 
we need to introduce new random variables such as $\wt Z_N$
with $\Law(\wt Z_N) = \Law(Z_N)$, etc.
For simplicity of presentation, however, we  ignore
this point here.}
to upgrade  the convergence in law to almost sure convergence
on a new probability space $\wt \O$. 
As in \cite{LLOT}, 
we then aim to 
prove that 
% almost sure convergence claimed in Theorem~\ref{THM:1}
% by establishing 
 the remainder term $v_N$ in  
 the second order expansion
\eqref{exp1} for the solution $u_N$ to \eqref{NLW1a}
converges almost surely 
to the remainder term $v$
in the first order expansion~\eqref{exp1a}
for  the solution $u$ to~\eqref{NLW1b}.
In the current NLW setting, however, there is 
a deficiency of spatial regularities
in making sense of 
the product $v_N z_N$ (in the limiting sense), 
which was not present in the BBM case studied in \cite{LLOT}.
In order to overcome this difficulty, 
we employ
 paracontrolled calculus
for the wave equation developed in \cite{GKO2, OOT1}
by 
further exploiting multilinear dispersion.
We point out that 
the aforementioned almost sure local well-posedness result in \cite{OO} 
%mentioned above
for 
(the Wick renormalized version of)~\eqref{NLW1}
with the Gaussian random initial data in \eqref{data1}
for $0 < \al < \frac 12$
(namely,  below the variance blowup threshold) 
follows from 
a straightforward modification of the argument in \cite{OPTz}
based on 
the second order expansion (analogous to \eqref{exp1})
{\it without} multilinear smoothing %of the second Picard iterate
or paracontrolled calculus.
In particular, the need for 
paracontrolled calculus arises only for $\al \ge \frac 12$.

Finally, given a random time $0 < T < T_{\max}$,
 the almost sure convergence 
 of $u_N$ to $u$
 via paracontrolled analysis
 on the new probability space $\wt\O$
 (Theorem \ref{THM:uNu}) 
allows us to construct 
random times $T_N$
such that the  pair
$\big(  T_N,  \vec u_N(\,\cdot\, \wedge T_N)  \big)$
converges almost surely to
$\big(  T,  \vec u(\,\cdot\, \wedge T)  \big)$
on the new probability space $\wt \O$, 
where we use the identifications
 $u_N = z_N + Z_N + v_N$
and $u = Z+v$.
This in turn implies
the convergence in law 
on the original probability space $\O$
claimed in Theorem \ref{THM:1}.

In Section \ref{SEC:sketch}, 
we present a sketch of the  proof of 
Theorem \ref{THM:1}, 
except for the last step.
See Subsection \ref{SUBSEC:main2} for the last step
where we transfer the almost sure convergence on the new probability space $\wt \O$
to the convergence in law 
on the original probability space $\O$, 
as claimed in Theorem \ref{THM:1}.

\subsection{``Beyond variance blowup'' as probabilistic ill-posedness}
\label{SUBSEC:1.2}

In the foundational works 
\cite{BO94, BO96}
 on invariant Gibbs measures for nonlinear
Schr\"odinger equations (NLS), 
Bourgain initiated the study of dispersive PDEs
with random initial data.
In
\cite{BT1}, Burq and the fifth author 
studied  the cubic NLW
 on a three-dimensional compact Riemannian manifold with randomized
initial data that are given as randomization of (the Fourier coefficients of) a given pair of
functions of the scaling supercritical regularity.
 Proceeding with the first
order expansion (as in \eqref{exp1a})
 and exploiting the almost sure averaging property of the random
linear solution, they established almost sure local well-posedness of the cubic NLW below the
scaling critical regularity. 
These works
set a foundational basis of a (somewhat informal) notion of probabilistic
local well-posedness, namely, 
almost sure existence and uniqueness of a solution with random initial data
with stability under a suitable class of perturbations.
Since then, there has been an extensive literature on probabilistic well-posedness of 
dispersive equations with  random initial data and\,/\,or (additive) stochastic forcing; 
see, for example, 
\cite{CO, BT3, BOP2, GKO, OTh2,
STz1, Bring0, STz2, OPTz, OOT1, Bring1,
OOT2, STzX, GKO2, ORTz, DNY2, DNY3, BDNY, BLL, BCST, LLO}.
See also surveys
\cite{BOP4, Tzv1, DNY5}.
We also mention a recent progress on pathwise well-posedness
of stochastic dispersive PDEs with multiplicative noises;
see \cite{CLO, CGLLO}
and the references therein.

In  \cite{DNY2}, 
Deng, Nahmod, and Yue introduced the notion of probabilistic scaling criticality
which provides useful heuristics 
for  threshold regularities
for probabilistic well-posedness
of dispersive PDEs with random initial data.
Nonetheless, 
in recent works
(see, for example,  \cite{Forlano, OO, Liu}), 
 it has become clear that even within the probabilistic framework, 
  one may encounter intrinsic
obstructions to constructing solutions, 
in particular before reaching probabilistic
scaling critical thresholds.
As mentioned above,  
in the case of the quadratic NLW~\eqref{NLW1} on~$\T^2$, 
variance blowup occurs for $\al \ge \frac 12$, 
whereas the 
probabilistic
scaling critical threshold is given by $\al_{\text{crit}} = \frac 54$;
see \cite[Subsubsection~2.2.2]{DNY5}.

%On the one hand, 
As explained in \cite{OO}, a variance blowup result implies that 
the standard probabilistic local well-posedness theory 
based on the first-order
expansion or its higher-order variants breaks down.
In a recent work \cite{OTz}, 
this viewpoint of probabilistic ill-posedness
was further explored by the fourth and fifth authors
by incorporating continuity\,/\,smoothness
of the solution map at the origin
in the enhanced notion of  probabilistic well-posedness.
In particular, the variance blowup result 
in \cite{OO}
for NLW \eqref{NLW1}
with the Gaussian random initial data $(u_0, u_1)$
in \eqref{data1}
when $\al \ge \frac 12$
implies
probabilistic ill-posedness in the mild sense
for $\al \ge \frac 12$;
 see \cite[Section~3]{OTz}.

The enhanced notion of probabilistic local well-posedness
in \cite[Definition 1.8]{OTz} includes the following continuity at the origin
(besides existence, uniqueness, and stability under the frequency truncation); 
given 
 any sequence 
 $\{\dl_N\}_{N \in \N}$
  of positive numbers  tending to $0$
 as $N \to \infty$, 
let 
$ (u_N, \dt u_N)$ be 
 the  (smooth) random solution  to 
 (the Wick renormalized version of) NLW~\eqref{NLW1} with the scaled frequency-truncated
 random initial data
$ (u_N, \dt u_N) |_{t = 0} = (\dl_N  \pi_N u_0, \dl_N  \pi_N u_1) $, 
where $(u_0, u_1)$ is as  in \eqref{data1}.
Then,  we have 
%
%we have, for each $\o \in \Si$, 
\begin{align}
( u_N, \dt u_N) \too (0, 0) 
\label{disp1}
\end{align}

\noi
 in  $C([0, T];\D'(\T^2)\times \D'(\T^2))$, 
 almost surely, as $N \to \infty$, 
 where  
 $T$ denotes the almost surely positive  local existence time.
We note that \eqref{disp1}
is satisfied by 
the aforementioned almost sure local well-posedness result  in~\cite{OO} 
on 
(the Wick renormalized version of)~\eqref{NLW1}
with the Gaussian random initial data in \eqref{data1}
for $0 < \al < \frac 12$.

As mentioned
at the end of Subsection \ref{SUBSEC:1.1}, 
the proof of 
our ``beyond variance blowup'' result 
(Theorem \ref{THM:1})
 yields the following almost sure convergence on the  new probability space
 $\wt \O$ (see Theorem \ref{THM:uNu}); 
given $\al \ge \frac12$, 
 by choosing $\dl_N = C_{\al, N}$ as in \eqref{CN}, 
$( u_N, \dt u_N) $ converges almost surely 
to the solution $(u, \dt u) \not \equiv(0, 0)$ of the stochastic NLW \eqref{NLW1b}, 
  violating \eqref{disp1}.
Hence, as a corollary to Theorem \ref{THM:1}, 
we obtain
 the following probabilistic ill-posedness result.
See \cite{OTz}
for a further discussion.

\begin{corollary}\label{COR:2}

Let $\al \ge \frac 12$.
Then, \textup{(}the Wick renormalized version of\textup{)} the quadratic NLW~\eqref{NLW1} on $\T^2$ with the
Gaussian random initial data $(u_0, u_1)$ in~\eqref{data1}
is probabilistically ill-posed.

\end{corollary}

In Appendix \ref{SEC:A}, 
we present a brief discussion on the proof of Corollary \ref{COR:2}.

\subsection{Remark \& comments}
We conclude this introduction by stating several remarks.

\begin{remark}
\label{REM:NLW2}
\rm 
(i)~Given $\al \ge 0$, we consider 
the following quadratic SNLW on $\T^2$: 
\begin{align} 
\begin{cases}
\dt^2 u + (1 - \Dl) u = u^2 + \jb{\nb}^\al \xi \\ 
(u, \dt u)|_{t =0} = (\phi_0, \phi_1), 
\end{cases} 
\label{SNLW2}
\end{align}

\noi
where $\xi$ denotes the space-time white noise on $\R_+ \times \T^2$
and  $(\phi_0, \phi_1)$ is sufficiently smooth deterministic initial data.
In \cite{GKO, OO}, it was shown that 
when $0 \le \al < \frac 12$, 
(the Wick renormalized version of) \eqref{SNLW2} is 
locally well-posed, 
whereas 
 variance blowup occurs 
when $\al \ge \frac 12$.

Fix $\al \ge \frac 12$. 
By studying 
SNLW \eqref{SNLW2} 
with the renormalized noise $C_{\al, N} \jb{\nb}^\al \pi_N \xi$, 
where $C_{\al, N}$ is as in \eqref{CN}, 
a straightforward modification of the proof of Theorem \ref{THM:1}
yields a ``beyond variance blowup'' convergence result
for SNLW.
See \cite[Appendix A]{LLOT} for a related discussion 
in the context of  the  stochastic BBM.

\medskip

\noi
(ii)
Let  $\al \ge \frac 12$.
Following \cite[Theorem 1.9]{LLOT} (see also \cite{GT25}), 
we can also consider the following weakly interacting NLW
on $\T^2$:
\begin{align}
\begin{cases}
\dt^2 u_N + (1 - \Dl) u_N = C_{\al, N}^2 u_N^2 - \s_N \\
(u_N, \dt u_N)|_{t = 0} = (\pi_N u_0, \pi_N u_1)
\end{cases} 
\label{WINLWnew}
\end{align}

\noi
for the Gaussian random initial data 
$(u_0, u_1) = (u_0^\o (\al), u_1^\o (\al))$ in \eqref{data1}, 
where $C_{\al, N}$ and $\s_N$ 
are as in \eqref{CN} and \eqref{sigN}, 
respectively.
Proceeding as in the proof of Theorem \ref{THM:1}, 
we can prove that the solution $(u_N, \dt u_N)$
to \eqref{WINLWnew}
converges (in an appropriate sense) 
to the solution $(u, \dt u)$ of the 
 following  stochastic linear wave equation: 
\begin{align}
\begin{cases}
\dt^2 u + (1 - \Dl) u =  \ze\\
(u, \dt u)|_{t = 0} = (u_0, u_1),  
\end{cases} 
\label{WINLWnew2}  
\end{align}

\noi
where 
  $\ze$ denotes the space-time noise 
 in \eqref{XI1} and \eqref{nxib},  independent of 
 the Gaussian random initial data $(u_0, u_1)$.
Unlike  the limiting equation \eqref{NLW1b}
in Theorem~\ref{THM:1}, 
the equation \eqref{WINLWnew2}
is linear and thus is 
 globally well-posed, 
 and hence the convergence holds
 (asymptotically) globally-in-time.

\end{remark}

\section{Sketch of the proof of Theorem \ref{THM:1}}
\label{SEC:sketch}

In this section, 
we briefly go over the results from \cite{OO}
and then describe the strategy for proving Theorem \ref{THM:1}.

\subsection{Second Picard iterate}
%Beyond variance blowup for the quadratic NLW}
\label{SUBSEC:2.1}

Let us briefly go over the local well-posedness and variance blowup results
from \cite{OO}.\footnote{As mentioned before, 
the paper \cite{OO} treated SNLW \eqref{SNLW2}
but their results equally apply to \eqref{NLW1} with the Gaussian initial data
in \eqref{data1}.}
Given $\al \in \R$, 
 let 
$(u_0, u_1) = (u_0^\o(\al), u_1^\o (\al))$
be the Gaussian random initial data in \eqref{data1}.
In the following, we fix  $\al \ge 0$, 
corresponding to the regime where \eqref{NLW1} is deterministically ill-posed.

Let $\zf$ be the random linear solution, given by 
\begin{align*}
\zf (t)= W(t) (u_0, u_1),
%\label{lin0}
\end{align*}

\noi
where $W(t)$ is as in \eqref{Wx}.
For $N \in \N$,  set
\begin{align*}
\zf_N  = \pi_N \zf .
%\label{zfN}
\end{align*}

\noi
Given any $(t, x) \in \R\times \T^2$, we have 
\begin{align*}
\ld_N =  \E \big[ \zf_N^2 (t, x ) \big]
= \sum_{|n| \le N}
\frac {1}{\jb{n}^{2 - 2\al}}
\sim
\begin{cases}
\log \jb{N},  &\text{if } \al = 0,  \\
N^{2 \al},   &\text{if } \al > 0, 
\end{cases}
\end{align*}

\noi 
which tends to $\infty$ as $N \to \infty$.
Then, we define the Wick power 
$\wick{\zf_N^2} $ by setting 
\begin{align*}
\wick{\zf_N^2} \, = \zf_N^2 - \ld_N.
%\label{Wickre}
\end{align*}

\noi
We also define the truncated second Picard iterate: 
\begin{align}
\begin{split}
\Zf_N (t)
= \I  (\, \wick {\zf_N^2}\,) (t), 
\end{split}
\label{exp2}
\end{align}

\noi
where $\I$ denotes the wave Duhamel integral operator defined in \eqref{Duhamel}.

We first consider the case $0 \le \al < \frac 12$. 
In this case, a standard computation shows that 
$\wick{\zf_N^2}$ and  $\Zf_N$
converge almost surely to non-trivial  limits, 
denoted by 
 $\wick{\zf^2}$ and $\Zf$, respectively,  
as $N \to \infty$, and we have 
\begin{align}
\Zf  = \I (\,\wick{\zf^2}\,) . 
\label{foe}
\end{align}

\noi
In view of one degree of smoothing from the wave Duhamel integral operator
$\I$ in \eqref{Duhamel}, 
one may naively expect that $\Zf$
is one degree smoother (in the spatial variable)
than $\wick{\zf^2}$\,.
However, by exploiting the explicit
product structure and multilinear dispersion, 
it was shown that $\Zf$ is in fact smoother; 
see \cite[Proposition 1.4]{OO}.
See also \cite[Proposition 1.6]{GKO2}, 
where such a multilinear smoothing was first established.
This multilinear smoothing plays an important role
in our ``beyond variance blowup'' result.

Once we construct the stochastic objects described above, 
almost sure local well-posedness of~\eqref{NLW1}
with the Gaussian initial data $(u_0, u_1)$ in~\eqref{data1}
follows from the second order expansion:
\begin{align}
u = \zf + \Zf + v
\label{exp2a}
\end{align}

\noi
and studying the equation satisfied by the remainder term $v$:
\begin{align}
\begin{cases}
\dt^2 v + (1- \Dl) v =  (v+\Zf)^2 + 2 \zf (v+\Zf)\\
(v, \dt v)|_{t = 0} = (0, 0).
\end{cases}
\label{NLW1y}
\end{align}

In particular, the smooth (random) solution $(u_N, \dt u_N)$
to 
\begin{align}
\begin{cases}
\dt^2 u_N + (1- \Dl) u_N = u_N^2 - \ld_N \\
(u_N, \dt u_N)|_{t = 0} = (\pi_N u_0, \pi_N u_1)
\end{cases}
\label{NLW1x}
\end{align}

\noi
converges almost surely to some non-trivial limit $(u, \dt u)$
in $C([0,T];\vec H^{-\al-\eps}(\T^2))$ as $N \to \infty$, 
where $T$ is an almost surely positive random time.
Here, the limit $u$ is given by \eqref{exp2a}, 
where the remainder term $v$ satisfies \eqref{NLW1y}.

\medskip

Next, we consider the case
 $\al \ge \frac 12$.
In this case,  
 it was shown 
 in \cite{OO}
that the
following variance blowup occurs:
\begin{align}
\lim_{N \to \infty}
\E\big[ \|  \ft \Zf_N (\,\cdot\,, n) \|_{L^\infty_T} ^2\big] = \infty
\label{div1}
\end{align}

\noi
for any $T > 0$ and any $n \in \Z^2$. %  (such that $t \jb{n} \ll 1$. 
Here, $\ft{f}$ denotes the Fourier transform of $f$ on $\T^2$.  
By adapting the argument in the proof of
Lemma \ref{LEM:COV}, 
we also have 
\begin{align}
\lim_{N \to \infty}
\E\big[\|\jb{  \Zf_N, \psi}\|_{L^\infty_T}^2\big] = \infty 
\label{div2}
\end{align}

\noi
for any $T > 0$
and any non-trivial test function $\psi \in \D(\T^2) \setminus \{0\}$, 
where 
$\jb{f, g}$
denotes the $\D'(\T^2)$-$\D(\T^2)$ duality pairing.
We note that 
 the divergence \eqref{div2},  arising from 
the `high $\times$ high into low' 
frequency interactions, 
takes place before 
reaching the probabilistic scaling critical threshold $\al_\text{crit} = \frac 54$.

\begin{remark}\label{REM:VB}\rm
As in \cite{Hairer, LLOT}, we use the term variance blowup
to describe the situation, 
where (i)~the variance of a certain stochastic term diverges 
(see \eqref{div2}
for the NLW case) at a threshold
{\it but} (ii)~the analytical framework still continues to hold beyond the threshold.
See the introduction in \cite{Hairer}
and \cite[Remark 1.2]{LLOT}
for a further discussion.

\end{remark}

\medskip

\noi
$\bul$ {\bf Convergence in law of the second Picard iterate.}
Fix $\al \ge \frac 12$.
In the remaining part of this section, 
 we go over our strategy for proving Theorem \ref{THM:1}
on  convergence 
of the solution $(u_N, \dt u_N)$
to \eqref{NLW1a}.
As mentioned in Section \ref{SEC:1}, 
the first key ingredient for proving Theorem \ref{THM:1}
is  the following  convergence in law of 
the second Picard iterate $Z_N$ in \eqref{ZN}
to the stochastic convolution  $Z$ defined in \eqref{Z1}.

\begin{theorem}
\label{THM:ZNZ}
Let $\al \ge \frac 12$
and  $s < \frac 14$.
Then,
as $N \to \infty$, 
 the second Picard iterate $Z_N$ defined in~\eqref{ZN}
converges in law to $Z$ defined in \eqref{Z1} in
$C(\R_+; W^{s, \infty} (\T^2)) \cap C^1 (\R_+; W^{s-1, \infty} (\T^2))$, 
endowed with the compact-open topology \textup{(}in time\textup{)}.
Moreover, 
$Z$ 
is independent of $(u_0,u_1)$.

\end{theorem}

The main novelty of the present paper
is the convergence in law of the second Picard iterate
$Z_N$~in Theorem~\ref{THM:ZNZ}.
While we follow the general strategy in \cite{LLOT},
based on the fourth moment theorem and tightness,
the dispersive nature of the wave equation
plays an essential role in our analysis.
In particular, we need to exploit the oscillatory nature
of the wave propagator
both to identify the limiting process $Z$
and to establish tightness in a space of {\it positive} regularity.
We briefly discuss these points below.

For the BBM equation studied in \cite{LLOT},
the dispersion relation $\varphi(n)=-in/(1+n^2)$ satisfies
\begin{align*}
\varphi(m)+\varphi(n-m)\longrightarrow0
\quad\text{as } |m|\to\infty
\end{align*}

\noi
for each fixed $n$.
This leads to the time-independent spatial noise
in the limiting equation.
For the wave equation, however, we have
\begin{align*}
\jb{m}-\jb{n-m}
=\frac{n\cdot m}{|m|}+O_n(|m|^{-1})
\quad\text{as } |m|\to\infty.
\end{align*}

\noi
In particular, in studying the limiting behavior
of  the 
space-time covariance:
\[\E\big[\jb{Z_N(r), \psi_1} \jb{ Z_N (t), \psi_2}\big]\]
in Lemma \ref{LEM:COV}, 
the oscillatory phase of 
the leading contribution
\begin{align}
C_{\al,N}^4
\sum_{\substack{n=n_1+n_2\\|n_1|,|n_2|\le N}}
\frac{\cos(t(\jb{n_1}-\jb{n_2}))}
{\jb{n_1}^{2-2\al}\jb{n_2}^{2-2\al}}
\label{Z2}
\end{align}

\noi
depends, to leading order,
on the direction of $m=n_2$.
By averaging over the angular variable
(see Lemma \ref{LEM:angavg}),
we show that
\begin{align*}
\eqref{Z2}
\longrightarrow
\frac{1}{2\pi}\int_0^{2\pi}
\cos(t|n|\cos\theta)d\theta
=J_0(t|n|)
\end{align*}

\noi
for each fixed $n\in\Z^2$ and $t\in\R$;
see Lemmas \ref{LEM:COVa} and \ref{LEM:COVb}.
When $\al>\frac12$,
we prove this convergence by a Riemann sum approximation,
after estimating the contribution near the origin separately.
At the endpoint $\al=\frac12$,
the radial integral diverges logarithmically,
and we instead use Lemma~\ref{LEM:angavg}.
In both cases, the normalization in~\eqref{CN}
cancels the radial contribution.
On the other hand, the contribution from the phase
$\jb{n_1}+\jb{n_2}$ vanishes 
(in the limit) after time integration.
This gives the Bessel covariance in \eqref{nxib};
see Subsections~\ref{SUBSEC:Riea} and \ref{SUBSEC:uniq}.

In proving tightness (Proposition~\ref{PROP:rand1}),
 we exploit
multilinear dispersive smoothing for $Z_N$
as in \cite{GKO2, OO}.
See Subsection~\ref{SUBSEC:tight}.
At the level of the limiting noise, 
this extra smoothing can also be seen directly
from its covariance.
Indeed, recalling  $|J_0(r)|\les\jb{r}^{-\frac12}$,
it follows from~\eqref{Z1} and \eqref{nxib}
that  for $n\ne0$, we have 
\begin{align}
\begin{split}
\E\big[|\ft Z(t,n)|^2\big]
&\les\jb{n}^{-2}
\int_0^t\int_0^t
\big|J_0((t_1-t_2)|n|)\big|dt_1dt_2\\
&\les\jb{n}^{-\frac52}
\int_0^t\int_0^t|t_1-t_2|^{-\frac12}dt_1dt_2
\les_t \jb{n}^{-\frac52}.
\end{split}
\label{Z3}
\end{align}

\noi
Thus, the time covariance of $\ze$
also gives the $\frac14$-smoothing for $Z$;
see Lemma \ref{LEM:randZ}.
This is a key difference from the BBM problem studied in~\cite{LLOT},
where the limiting noise is independent of time.

We present the proof of Theorem~\ref{THM:ZNZ}
in Section~\ref{SEC:Psecond}.
As for the last claim in Theorem~\ref{THM:ZNZ}
on the independence
of $Z$ and  $(u_0,u_1)$, 
see the proof of Lemma \ref{LEM:SEP}.
See also \cite[Remark 3.6]{LLOT}.

\begin{remark}\rm
\label{REM:ZNZ}

We point out that multilinear dispersive smoothing
induced by the Duhamel integral operator $\I$ in \eqref{Duhamel}
is essential for Theorem \ref{THM:ZNZ}.
While Theorem \ref{THM:ZNZ}
implies 
convergence in law of  $\wick{z_N^2}$  
as space-time distributions, 
%it seems that 
$\wick{z_N^2}$  does not converge
as continuous-in-time processes.
(For example, one can easily see that the Fourier zero mode is not tight.)

\end{remark}

\subsection
{Paracontrolled analysis}
%{Convergence of the full dynamics}
\label{SUBSEC:IPa}

Once we establish convergence in law of $(z_N, Z_N)$
to $(0, Z)$, 
our next goal is to study uniform (in $N$)
local well-posedness of \eqref{NLW1a}.
Under the second order expansion \eqref{exp1}, 
 the remainder term $v_N$ satisfies
\begin{align}
\begin{cases}
\dt^2 v_N + (1 - \Dl) v_N
= (v_N + Z_N)^2 + 2 v_N z_N + 2 Z_N z_N\\
(v_N, \dt v_N)|_{t=0} = (0, 0).
\end{cases}
\label{NLW3}
\end{align}

 \noi
 The following discussion is parallel to that in \cite[Subsection 1.3]{GKO2}, 
 which we include for readers' convenience.
 The worst term on the right-hand side of \eqref{NLW3}
 is $z_N$ with regularity $-\frac 12-$.
 From the paraproduct consideration (see \eqref{para1}
 and Lemma \ref{LEM:para}), 
we expect 
the product
  $Z_N z_N$  has regularity $-\frac12-$, 
inheriting the worse regularity of $z_N$. 
Thus,  $v_N$ is expected to have regularity at most $\frac12-=\big(-\frac12-\big)+1$. In
particular, the product $v_N z_N$ is not well defined, since
the sum of the regularities is negative.

In \cite{GKO2}, 
Gubinelli, Koch, and the fourth author encountered the same issue, 
which they overcame by adapting
paracontrolled calculus, originally introduced
in the parabolic context  by
Gubinelli, Imkeller, and Perkowski \cite{GIP}
(see also \cite{MW2}), 
to the wave setting.
We 
proceed as in~\cite{GKO2} and introduce our paracontrolled ansatz;
write  $v_N$ as 
\begin{align}
v_N = X_N + Y_N,
\label{exp4}
\end{align}

\noi
where $X_N$ and $Y_N$ satisfy
\begin{align}
 \dt^2 X_N + (1 - \Dl) X_N  & =  2 (X_N +Y_N + Z_N)\pl z_N, 
 \label{NLW4a}\\
 \dt^2 Y_N + (1 - \Dl) Y_N
&  = (X_N + Y_N + Z_N)^2 +  2(X_N + Y_N + Z_N) \pge z_N
\label{NLW4b}
\end{align}

\noi
with the initial data $(X_N, \dt X_N, Y_N, \dt Y_N) = (0, 0, 0, 0)$.
See Subsection \ref{SUBSEC:space}
for the definitions of 
 the paraproducts and resonant products.
We note that in \eqref{exp4}, 
 \eqref{NLW4a}, and  \eqref{NLW4b}, 
we are in fact postulating that 
$\dt^2 v_N + (1 - \Dl)v_N$ is paracontrolled by $z_N$; 
see  \cite[Remark~1.7]{GKO2} 
for a further discussion.

Assuming that 
 $X_N$ and $Y_N$ have positive regularities, 
 it follows 
from the first equation~\eqref{NLW4a} that 
$X_N$ has regularity $\frac 12- = \big(-\frac{1}{2}-\big) + 1$.
For now, let us ignore the problematic resonant product
$2(X_N + Y_N + Z_N) \pe z_N$ in \eqref{NLW4b}
and discuss the regularity of $Y_N$.
Recalling that $Z_N$ has regularity $\frac 14-$,
we expect that  $Y_N$ has regularity $\frac 34 - = \big(-\frac 14 -\big) + 1$, 
which comes from the 
paraproduct $Z_N \pg z_N$ (with regularity $- \frac 14 -$).
In this discussion, we ignored
 the resonant product
$(X_N + Y_N + Z_N) \pe z_N$, 
to which we give a meaning below:

\smallskip
\begin{itemize}
\item
$Y_N \pe z_N$
 makes sense 
deterministically
via the paraproduct estimate (Lemma \ref{LEM:para}), 
since the sum of regularities is positive.

\smallskip
\item
Proceeding with stochastic analysis 
and multilinear oscillatory sum analysis as in~\cite{GKO2}, 
we can make sense of 
the following  resonant product:
\begin{align}
R_N := Z_N \pe z_N
\label{Zrz}
\end{align}

\noi
with  regularity
$- \frac 14 - = \big(\frac{1}{4}-\big) + \big(-\frac 12 -\big)$;
see Lemma \ref{LEM:rand2}.

\end{itemize}

\smallskip
\noi
Furthermore, both 
$Y_N \pe z_N$ and $R_N$ tend to $0$ as $N \to \infty$.

Due to the deficiency of regularity, 
we can not apply the paraproduct estimate (Lemma~\ref{LEM:para})
to make sense of  the resonant product $X_N \pe z_N$.
We instead use the fact that $X_N$ is the solution to 
the equation~\eqref{NLW4a}.
Namely, $X_N$ satisfies 
 the following Duhamel formulation:
\begin{align}
 X_N  =
 2\I \big((X_N +Y_N + Z_N) \pl z_N\big)
\label{X1}
 \end{align}

\noi
whose structure we can exploit in giving a meaning to 
 the resonant product $X_N \pe z_N$.

Following \cite{GKO2}, we now introduce
the paracontrolled operator
$\If_{\pl}^N$ by setting
\begin{align}
\begin{split}
\If_{\pl}^N (w) 
 =  \I (w \pl z_N)
\end{split}
\label{X2}
\end{align}

\noi
for  a function $w \in C(\R_+; H^{s_1}(\T^2) )$ with some $s_1 > 0$.
Next, we 
 decompose 
$\If_{\pl}^N$ into two pieces.
Let  $n_1$ and $n_2$ be the spatial frequencies
of $w$ and $z_N$ as in \eqref{X2}.
Given
 small $\theta  > 0$, 
we write
\begin{align}
\If_{\pl}^N = 
 \If_{\pl}^{N, (1)}+ \If_{\pl}^{N, (2)},
\label{X2a}
\end{align}

\noi
 where
$ \If_{\pl}^{N, (1)}$ and $ \If_{\pl}^{N, (2)}$
 are
the restrictions of $\If_{\pl}^N$ onto $\{ |n_1| \ges |n_2|^{\theta}\}$
and  $\{ |n_1| \ll |n_2|^{\theta}\}$.
As we see in Lemma~\ref{LEM:rand3},
 the first paracontrolled operator $ \If_{\pl}^{N, (1)}$
enjoys extra regularity
thanks to the lower bound  $|n_1| \ges |n_2|^\theta$
and the positive regularity of $w$.
As a consequence, 
 the resonant product
\begin{align}
 \If_{\pl}^{N, (1)}(X_N + Y_N + Z_N) \pe z_N
\label{X2b} 
\end{align}

\noi
makes sense 
via the paraproduct estimate (Lemma \ref{LEM:para}), 
 and tends to $0$ as $N \to \infty$.
See also Remark \ref{REM:Z1}.

\begin{remark}\label{REM:para1}\rm
Given a space-time distribution $\Xi_1$, let us  define $\If^{\Xi_1}_{\pl}$ by 
\begin{align*}
\If^{\Xi_1}_{\pl} (w) 
 =  \I (w \pl \Xi_1).
%\label{XY2}
\end{align*}

\noi
Then, as in \eqref{X2a}, we define $ \If_{\pl}^{\Xi_1, (1)}$
as the restriction of 
$\If^{\Xi_1}_{\pl}$
onto $\{ |n_1| \ges |n_2|^{\theta}\}$.
In particular, we have
$\If^{N, (1)}_{\pl} = \If^{z_N, (1)}_{\pl}$.
\end{remark}

In order to handle the remaining contribution to $X_N \pe z_N$, 
we introduce
the second paracontrolled operator
$\If_{\pl, \pe}^N$  by setting
\begin{align}
\If_{\pl, \pe}^N (w) 
  =
  \If_{\pl}^{N, (2)}(w) \pe z_N
\label{X6}
\end{align}

\noi
for $w$ as in \eqref{X2}.
In Lemma \ref{LEM:rand4}, 
we directly study its almost sure mapping property via stochastic analysis
and exploiting multilinear dispersion
and the sine cancellation.
This allows us to give a meaning
to  the remaining contribution:
\begin{align}
\If_{\pl}^{N, (2)}(X_N + Y_N + Z_N) \pe z_N
= \If_{\pl, \pe}^N (X_N + Y_N + Z_N) , 
\label{X2c} 
\end{align}

\noi
where we interpret the right-hand side
as the definition of the left-hand side.
Moreover, the right-hand side tends to $0$ as $N \to \infty$.
Putting \eqref{X1}, \eqref{X2b}, and \eqref{X2c} together, 
we interpret the resonant product $X_N \pe z_N$
as
\begin{align}
X_N \pe z_N
= 2  \If_{\pl}^{N, (1)}(X_N + Y_N + Z_N) \pe z_N
+ 2 \If_{\pl, \pe}^N (X_N + Y_N + Z_N) .
\label{X2d}
\end{align}

Therefore, from \eqref{NLW4a}, \eqref{NLW4b}, and \eqref{X2d}, 
we obtain the following paracontrolled  system
for $X_N$ and $Y_N$:
\begin{align}
\begin{split}
\dt^2 X_N  + (1- \Dl) X_N
&  =  2 (X_N + Y_N + Z_N) \pl z_N, \\
 \dt^2 Y_N + (1- \Dl) Y_N
&   =  (X_N + Y_N + Z_N)^2 + 2(X_N + Y_N + Z_N) \pg z_N \\
& \hphantom{X}
  + 2 Y_N \pe z_N + 2 R_N 
%& \hphantom{X}
  + 4 \If_{\pl}^{N, (1)}(X_N + Y_N + Z_N) \pe z_N \\
& \hphantom{X}
  + 4 \If_{\pl, \pe}^N (X_N + Y_N + Z_N), \\
(X_N, \dt X_N, Y_N, \dt Y_N)|_{t = 0} & = (0, 0, 0, 0).
\end{split}
\label{NLW5}
\end{align}

\noi
The (random) enhanced data set for the  paracontrolled system \eqref{NLW5} is  given by 
\begin{align}
\begin{split}
\Xi_N = (z_N, Z_N, R_N,
\If_{\pl, \pe}^N), 
\end{split}
\label{data2}
\end{align}

\noi
where $z_N$, 
$Z_N$, $R_N$,
and $\If_{\pl, \pe}^N$
are as in \eqref{zN}, \eqref{ZN}, 
\eqref{Zrz}, and 
\eqref{X6}, respectively.

\medskip

In the following, we consider the following deterministic version of the
paracontrolled  system~\eqref{NLW5}:
\begin{align}
\begin{split}
\dt^2 \Xs  + (1- \Dl) \Xs
&  =  2 (\Xs + \Ys + \Xi_2) \pl \Xi_1, \\
 \dt^2 \Ys + (1- \Dl) \Ys
&   =  (\Xs + \Ys + \Xi_2)^2 + 2(\Xs + \Ys + \Xi_2) \pg \Xi_1 \\
& \hphantom{X}
  + 2 \Ys \pe \Xi_1 + 2 \Xi_3 
  + 4 \If_{\pl}^{\Xi_1, (1)}(\Xs + \Ys + \Xi_2) \pe \Xi_1 \\
& \hphantom{X}
  + 4 \Xi_4 (\Xs + \Ys + \Xi_2), \\
(\Xs, \dt \Xs, \Ys, \dt \Ys)|_{t = 0} & = (0, 0, 0, 0)
\end{split}
\label{NLW5x}
\end{align}

\noi
with 
 a deterministic enhanced data set: 
 \begin{align*}
 \Xi = (\Xi_1, \Xi_2, \Xi_3, \Xi_4), 
% \label{data2a}
 \end{align*}
 where 
$ \If_{\pl}^{\Xi_1, (1)}$ is as in Remark \ref{REM:para1} and 
 the elements $\Xi_1$, $\Xi_2$, $\Xi_3$, and  $\Xi_4$
 replace 
$z_N$, $Z_N$, $R_N$, and 
$\If_{\pl, \pe}^N$, 
appearing in \eqref{data2}, respectively.
Note that $\Xi_4$ 
(replacing $\If_{\pl, \pe}^N$) is an operator, acting on space-time functions.

Given $s, s_1, s_2 \in \R$ and $T > 0$, we define
$X^s(T)$ and $\vec X^{s_1, s_2}(T)$ by setting
\begin{align}
\begin{split}
X^{s}(T) &= C([0, T]; H^s(\T^2))
\cap C^1([0, T]; H^{s - 1} (\T^2)), \\
\vec X^{s_1, s_2}(T )
& =  X^{s_1}(T ) \times  X^{s_2}(T ).
\end{split}
\label{Xsp}
\end{align}

\noi
Given  $1 < q < \infty$, 
we also define the space $\L (q, T, s_2)$ of operators by setting
\begin{align}
 \L (q, T, s_2) = \L\big(L^q ([0, T]; L^2(\T^2) )
\, ;\,
C([0, T]; H^{s_2 -1}(\T^2))\big).
\label{L2}
\end{align}

\noi
Then, we define
 $\L_{a}(q,T, s_2)$ to be  the closed subspace
of ``adapted''\footnote{Here, the term ``adapted'' means that the output up to time $t$ depends only 
on the information of the input function $w$ up to time $t$.} operators in $\L(q,T, s_2)$.
Namely, $L\in\L_{a}(q,T, s_2)$ if
\begin{align}
(Lw)|_{[0,t]}
=\big(L(\ind_{[0,t]}\cdot w)\big)|_{[0,t]}
\label{adapt1}
\end{align}

\noi
for every $t\in[0,T]$ and $w$.
The norm on this subspace is the operator norm
in \eqref{L2}.
Given small $\eps > 0$, we   define the space $\Xc^{\eps, s_2}$
of enhanced data sets by setting
\begin{align}
\begin{split}
\Xc^{\eps, s_2}
& =
C([0,1]; W^{-\frac 12 - \eps, \infty}(\T^2))
\times
C([0,1]; W^{\frac 14 - \eps, \infty}(\T^2))
\\
& \hphantom{X}
\times
C([0,1]; W^{ - \frac 14 - \eps, \infty}(\T^2))
\times \L_{a}(q,1, s_2)
\end{split}
\label{L3}
\end{align}

\noi
for some (fixed) $1<q<\infty$.
Then, 
the following local well-posedness result 
for the deterministic 
 paracontrolled 
system \eqref{NLW5x} follows from 
 the standard energy estimate \eqref{A4}.

\begin{proposition}
\label{PROP:LWP1}
Let $\frac14<s_1<\frac12<s_2<\frac34$
with $s_2\le2s_1$.
Then, there exist small $\ta = \ta( s_2) > 0$
and $\eps = \eps (s_1, s_2) > 0$
such that, 
given $ \Xi = (\Xi_1, \Xi_2, \Xi_3, \Xi_4)
 \in \Xc^{\eps, s_2}$, 
 the paracontrolled  system~\eqref{NLW5x} is locally well-posed
 in the sense that 

\smallskip

\begin{itemize}
\item
there exists a unique solution $(\Xs, \Ys) \in   \vec X^{s_1, s_2}(T )$
to \eqref{NLW5x} on the time interval $[0, T]$
for some local existence time $T = T\big(\|\Xi\|_{\Xc^{\eps, s_2}}\big) \in (0, 1]$, 

\smallskip

\item 
 the solution map\,$:
\Xi \in \Xc^{\eps, s_2}
\mapsto 
(\Xs, \Ys) \in   \vec X^{s_1, s_2}(T )$
is locally Lipschitz continuous.

\end{itemize}

\smallskip

\end{proposition}

We present a proof of Proposition \ref{PROP:LWP1} 
in  Section \ref{SEC:conv1}.   
It follows from 
Theorem~\ref{THM:ZNZ}
and Lemmas~\ref{LEM:randln}, \ref{LEM:rand2}, 
 \ref{LEM:rand3}, and  \ref{LEM:rand4}
that the enhanced data set 
$\Xi_N$ in \eqref{data2} almost surely belongs to $\Xc^{\eps, s_2}$.
Hence, 
almost sure local well-posedness
of 
the paracontrolled  system \eqref{NLW5}
follows from Proposition~\ref{PROP:LWP1}.

\subsection{Convergence of the full dynamics}
\label{SUBSEC:full}

In view of the convergence
in law of the enhanced data set
$\Xi_N = (z_N, Z_N, R_N,
\If_{\pl, \pe}^N)$ in \eqref{data2}
to $(0, Z, 0, 0)$, 
where $Z$ is as in \eqref{Z1}, 
by formally taking $N \to \infty$, 
we see that the paracontrolled  system \eqref{NLW5} reduces to 
\begin{align}
\begin{split}
\dt^2 X  + (1- \Dl) X
&  = 0, \\
 \dt^2 Y + (1- \Dl) Y
&   =  (X + Y + Z)^2, \\
(X, \dt X, Y, \dt Y)|_{t = 0} & = (0, 0, 0, 0), 
\end{split}
\label{NLW5a}
\end{align}

\noi
which in particular implies $X \equiv 0$.
Namely, 
the solution $(X_N,  Y_N)$ to \eqref{NLW5}
formally converges to $(0, v)$, 
where $v$ is the solution to
\begin{align}
\begin{cases}
\dt^2 v + (1 - \Dl) v = (v +  Z)^2 \\
(v, \dt v)|_{t = 0} = (0, 0).  
\end{cases}
\label{NLW3b}
\end{align}

\noi
Given the almost sure regularity property of $Z$
from  Theorem \ref{THM:ZNZ}, 
we see that 
the enhanced data set $\Xi = (0, Z, 0, 0)$
almost surely belongs to $\Xc^{\eps, s_2}$ for any $\eps > 0$
and $s_2 < 1$.
Hence, 
it follows from 
Proposition \ref{PROP:LWP1}
that the system \eqref{NLW5a}
is almost surely locally well-posed
with $X\equiv 0$ and $Y \in X^{s_2}(T)$
for some $s_2 > \frac 12$ (sufficiently close to $\frac 12$), 
where $T$ denotes an almost surely positive local existence time.
Under the first order expansion $u = Z + v$ in \eqref{exp1a}, 
this provides an alternative argument
for local well-posedness of 
the limiting equation \eqref{NLW1b}
with  $v \in X^{s_2}(T)$,
$\frac12<s_2<\frac34$.

\begin{remark}\label{REM:uniq1}\rm

In Section \ref{SEC:B}, we present the proof of Proposition \ref{PROP:B1}\,(i)
on  local well-posedness
of the limiting equation \eqref{NLW1b}
by establishing local well-posedness
of  the $v$-equation \eqref{NLW3b} via  a contraction argument.
This argument yields
a local-in-time solution  $v \in X^s(\tau)$, 
$0 < s < \frac 14$, for an almost surely positive 
local existence time $\tau$, 
where $X^s(\tau)$ is as in \eqref{Xsp}.
Moreover, this contraction argument 
 yields the following blowup alternative:
\begin{align}
 T_{\max} = \infty \qquad \text{or}\qquad
\lim_{t \to T_{\max}-}
\|(v(t), \dt v(t))\|_{\vec H^s} = \infty, 
\label{F5b}
\end{align}

\noi
almost surely.
By  iterating the local well-posedness argument, 
we can extend $v$ to the maximal time interval $[0, T_{\max})$, 
where $T_{\max} < \infty$ almost surely 
and hence
\begin{align}
\lim_{t \to T_{\max}-}
\|(v(t), \dt v(t))\|_{\vec H^s} = \infty, 
\label{FF1}
\end{align}

\noi
almost surely.

On the other hand, 
as mentioned above, 
Proposition \ref{PROP:LWP1} also yields local well-posedness
of the $v$-equation \eqref{NLW3b}
by viewing it 
as \eqref{NLW5x}
with $(\Xs, \Ys) = (0, v)$
and the enhanced data set $\Xi = (0, Z, 0, 0)$.
The solution thus constructed, call it $\wt v$, 
belongs to $X^{s_2}(\wt \tau)$, 
$\frac12<s_2<\frac34$, 
 for an almost surely positive 
local existence time $\wt \tau$.
Let $\wt T_{\max}$ be a random maximal existence
time for $\wt v$, based on this construction at the spatial regularity $s_2$. 
This local well-posedness argument
via Proposition \ref{PROP:LWP1}
also yields a blowup alternative analogous to \eqref{F5b}:
\begin{align}
\wt  T_{\max} = \infty \qquad \text{or}\qquad
\lim_{t \to \wt T_{\max}-}
\|(\wt v(t), \dt \wt v(t))\|_{\vec H^{s_2}} = \infty, 
\label{FF2}
\end{align}

\noi
almost surely.

We first point out that the contraction argument
in  the proof of 
Proposition \ref{PROP:B1}\,(i) presented in Section \ref{SEC:B}
yields unconditional uniqueness of a solution to \eqref{NLW3b}.
In particular, we have $v = \wt v$ on $[0, \wt T_{\max})$, 
almost surely.
It is, however, possible that 
$\wt T_{\max} < T_{\max}$, namely, the $\vec H^{s_2}$-norm
of $(v(t), \dt v(t))$ diverges before reaching $T_{\max}$.
In Remark \ref{REM:uniq2},  we provide a brief argument
showing that this does not happen.
Namely, we have 
\begin{align}
\text{$T_{\max} = \wt T_{\max}$\quad 
and \quad $v = \wt v$ \quad on $[0, T_{\max})$, \ almost surely.}
\label{FF2a}
\end{align}

\end{remark}

\medskip

Let us continue with the discussion 
on convergence of the full dynamics.
From 
\eqref{exp1}
and~\eqref{exp4}, we have
\begin{align*}
u_N = z_N+Z_N + X_N + Y_N.
%\label{dec1}
\end{align*}

\noi
By formally taking $N \to \infty$, 
this converges to the following first order expansion:
\begin{align}
u = Z + v. 
\label{dec2}
\end{align}

\noi 
In view of 
\eqref{NLW3a} and \eqref{NLW3b}, 
we see that $u$ 
is the solution to the stochastic NLW \eqref{NLW1b}.

\medskip

In order to carry out pathwise analysis, 
we now
 apply
 the Skorokhod representation theorem 
(Lemma \ref{LEM:Sk})
and upgrade the convergence in law 
of the enhanced data set $\Xi_N$ in~\eqref{data2}
to almost sure convergence, but on a new probability space.
More precisely, 
we  apply
 the Skorokhod representation theorem 
(Lemma \ref{LEM:Sk})
to 
the extended enhanced data set:
\begin{align}
%(\pi_N u_0, \pi_N u_1, \Xi_N) = 
(u_0, u_1, z_N,  Z_N, R_N, 
\If_{\pl, \pe}^N) 
\label{data2x}
\end{align}

\noi
converging in law to 
\begin{align}
(u_0, u_1, 0,  Z, 0, 0), 
\label{data2y}
\end{align}

\noi
where $Z$ in \eqref{data2y}
is %taken to be 
independent of $(u_0,u_1)$.
See Subsection \ref{SUBSEC:Sk0}
for a discussion on the joint convergence in law of \eqref{data2x}
to \eqref{data2y}.
This creates
a new probability space $(\wt \O, \wt \F, \wt \PP)$, 
 new extended enhanced data sets:
\begin{align}
%(\pi_N u_0, \pi_N u_1, \Xi_N) = 
(\wt u_0^{(N)}, \wt u_1^{(N)}, \wt z_N,  \wt Z_N, \wt R_N, 
\wt \If_{\pl, \pe}^N)
\label{data2z}
\end{align}

\noi
and $(\wt u_0, \wt u_1, 0, \wt Z, 0 , 0 )$ such that

\smallskip
\begin{itemize}
\item[(i)] As $N \to \infty$, we have 
\begin{align}
(\wt u_0^{(N)}, \wt u_1^{(N)}, \wt z_N,  \wt Z_N, \wt R_N, 
\wt \If_{\pl, \pe}^N)
\too
(\wt u_0, \wt u_1, 0,  \wt Z, 0, 0), \quad \text{$\wt \PP$-almost surely}
\label{data2zz}
\end{align}

\noi
(with respect to the topology of $(\M, d_\M)$
defined in \eqref{A34} and \eqref{JW4}).

\smallskip
\item[(ii)]
For each $N \in \N$, we have 
\begin{align}
\begin{split}
\Law \big(\wt u_0^{(N)}, \wt u_1^{(N)}, 
\wt z_N,  \wt Z_N, \wt R_N, 
\wt \If_{\pl, \pe}^N\big)
& = 
\Law (u_0, u_1,  z_N,  Z_N, R_N, 
\If_{\pl, \pe}^N),\\ 
\Law \big((\wt u_0, \wt u_1), \wt Z\big)
& = 
\Law (u_0, u_1)\otimes 
\Law (Z).
\end{split}
\label{Law1}
\end{align}

\noi
In particular, for any $N \in \N$, we have 
\begin{align}
\begin{split}
\Law (\wt u_0^{(N)}, \wt u_1^{(N)})
= 
\Law (\wt u_0, \wt u_1)
= \Law (u_0, u_1)
=   \mu_{1-\al}\otimes \mu_{-\al}, 
\end{split}
\label{Law2}
\end{align}

\noi
where $\mu_s$ is as in \eqref{mu1}.

\smallskip

\end{itemize}

\noi
We point out   that the space $\L(q, T, s_2)$  of 
operators 
defined in  \eqref{L2}
is not separable
and thus we need to proceed with care.
See
Section \ref{SEC:Sk}
for a further discussion.

\begin{remark}\label{REM:IV2}\rm
In \eqref{data2z}, 
we introduced
the initial data 
$(\wt u_0^{(N)}, \wt u_1^{(N)})$, depending on $N \in \N$ (while satisfying \eqref{Law2}).
This is due to the fact that, 
using one fixed copied initial pair for every $N \in \N$
(namely, $(\wt u_0^{(N)}, \wt u_1^{(N)}) 
= (\wt u_0, \wt u_1)$, $N \in \N$), 
the $\wt \PP$-almost sure convergence~\eqref{data2zz}
with \eqref{Law1} 
fails.
See Remark \ref{REM:IV1} for a further discussion.

\end{remark}

Given $N \in \N$,  as in \eqref{X2}
and \eqref{X2a}, define
the paracontrolled operator $\wt \If_{\pl}^N$ by setting
\begin{align*}
\wt \If_{\pl}^N (\,\cdot \,) 
=  \If_{\pl}^{\wt z_N} (\,\cdot \,) 
 =  \I (\,\cdot \, \pl \wt z_N)
\end{align*}

\noi
and, given  small $\theta  > 0$, 
split $\wt \If_{\pl}^N $ as 
\begin{align*}
\wt \If_{\pl}^N = 
\wt  \If_{\pl}^{N, (1)}+ \wt \If_{\pl}^{N, (2)}.
%\label{Law4}
\end{align*}

\noi
Then, for each $N \in \N$, 
we have 
\begin{align}
\begin{split}
\wt z_N (t)
& = 
W(t)\big(C_{\al, N} \pi_N \wt u_0^{(N)}, C_{\al, N}\pi_N  \wt u_1^{(N)}\big),\\
\wt  Z_N & = \I(\, \wick{\wt z_N^2} \, ) =  \I(\wt z_N^2 - \s_N ), \\
\wt R_N & = \wt Z_N \pe \wt z_N, \\
\wt \If_{\pl, \pe}^N(\,\cdot\,)
& = \wt \If_{\pl}^{N, (2)}(\,\cdot \,) \pe \wt z_N, 
\end{split}
\label{Law5}
\end{align}

\noi
$\wt \PP$-almost surely, 
where $W(t)$, 
 $C_{\al, N}$, and 
$\s_N$ are as in 
\eqref{Wx}, 
 \eqref{CN}, and 
\eqref{sigN}.
Note that the left-hand sides of~\eqref{Law5}
are as in \eqref{data2z}, defined 
in an abstract manner via the application of 
 the Skorokhod representation theorem 
(Lemma~\ref{LEM:Sk}).
The identities in \eqref{Law5} show that 
each of $\wt z_N$, 
$\wt  Z_N$, 
$\wt R_N$, and 
$\wt \If_{\pl, \pe}^N$
has the right structure as in \eqref{zN}, \eqref{ZN}, \eqref{Zrz}, and \eqref{X6}
on the original probability space $\O$.
We prove the identities~\eqref{Law5}
in Subsection~\ref{SUBSEC:Sk2}.

Define $\wt \ze$
by setting 
\begin{align}
\wt \ze 
=(\dt^2+1-\Dl)\wt Z.
\label{Law6}
\end{align}

\noi
Since 
$\Law (\wt Z) = \Law (Z)$, we have 
 $(\wt Z, \dt \wt Z)|_{t = 0} = (0, 0)$, $\wt \PP$-almost surely, 
 and 
\begin{align}
\begin{split}
\Law (\wt \ze) 
& = \Law \big((\dt^2+1-\Dl)\wt Z\big)\\
& = \Law \big((\dt^2+1-\Dl) Z\big)
= \Law (\ze), 
\end{split}
\label{Law7}
\end{align}

\noi
where $\ze$ is as in \eqref{XI1} and~\eqref{nxib}. 
Hence, we obtain
\begin{align}
\wt Z(t)
= 
\int^t_0 \frac{\sin((t-t')\jb{\nb})}{\jb{\nb}}
\wt \ze(t') dt'.
\label{Law8}
\end{align}

\noi
See the proof of Proposition~\ref{PROP:Sk}\,(iii)
for further details.

%
% \wt Z_N, \wt R_N, 
%\wt \If_{\pl, \pe}^N), 

We now consider the following  NLW
with the renormalized (frequency-truncated) random initial data
on the new probability space $\wt \O$:
\begin{align} 
\begin{cases}
\dt^2 \wt u_N + (1 - \Dl) \wt u_N = \wt u_N^2 - \s_N \\  
(\wt u_N, \dt \wt u_N)|_{t = 0} = 
\big(C_{\al, N}\pi_N \wt u_{0}^{(N)}, C_{\al, N} \pi_N \wt u_{1}^{(N)}\big), 
\end{cases} 
\label{XNLW1}
\end{align}

\noi
where $\s_N$ is defined in \eqref{sigN}, 
and the limiting  SNLW: 
\begin{align} 
\begin{cases}
\dt^2 \wt u + (1 - \Dl) \wt u = \wt u^2 +  \wt \ze\\ 
(\wt u, \dt \wt u)|_{t = 0} = (0, 0).
\end{cases} 
\label{XNLW2}
\end{align}

\noi
With the second order expansion (as in \eqref{exp1})
and the paracontrolled ansatz (as in~\eqref{exp4}):
\begin{align}
\wt u_N = \wt z_N + \wt Z_N + \wt X_N + \wt Y_N, 
\label{exp5}
\end{align}

\noi
we see that the 
pair $(\wt X_N,  \wt Y_N)$
satisfies the  paracontrolled system
\eqref{NLW5x}
with 
\[   (\Xi_1, \Xi_2, \Xi_3, \Xi_4) = 
(\wt z_N, \wt Z_N, \wt R_N,
\wt \If_{\pl, \pe}^N) =:\wt \Xi_N
\]

\noi
Similarly, under the following expansion: 
\begin{align}
\wt u = 0 +  \wt Z  + \wt X + \wt Y  \quad \text{with }\,  (\wt X, \wt Y)  = (0, \wt v), 
\label{exp6}
\end{align}

\noi
 the 
pair $(\wt X, \wt Y)  = (0, \wt v)$
satisfies the  paracontrolled system
\eqref{NLW5x}
with 
\begin{align}
   (\Xi_1, \Xi_2, \Xi_3, \Xi_4) = 
(0, \wt Z, 0, 0)=:\wt \Xi.
\label{exp7}
\end{align}

%
%\noi
%the remainder term $\wt v$ satisfies
%\eqref{NLW3b}, where $Z$ is replaced by $\wt Z$ in \eqref{Law6}.

\noi
Hence, 
in view of 
\eqref{Law1}
with 
Theorem~\ref{THM:ZNZ}
and Lemmas \ref{LEM:randln}, \ref{LEM:rand2}, 
 \ref{LEM:rand3}, and  \ref{LEM:rand4}, 
it follows from Proposition \ref{PROP:LWP1}
that \eqref{XNLW1} for each $N \in \N$
and \eqref{XNLW2}
are almost surely locally well-posed.
We  note that, in view of \eqref{Law7},  
local well-posedness
of the limiting SNLW \eqref{XNLW2}
also follows from Proposition \ref{PROP:B1}\,(i) 
and that the solution $\wt u$ to \eqref{XNLW2}
almost surely blows up in finite time.
Let $\wt T_{\max}$ denote
the almost surely positive maximal existence time
for the solution 
 $(\wt u, \dt \wt u) $ to~\eqref{XNLW2}
 such that 
  $(\wt u, \dt \wt u)  \in C([0, \wt T_{\max}); \vec H^s(\T^2))$, $0 < s < \frac 14$.
Then, we have
\begin{align}
%\wt  \PP (  \wt T_{\max} < \infty  ) = 1
\wt T_{\max} < \infty
\qquad \text{and}
\qquad
\lim_{t \to \wt T_{\max}-}
\|(\wt u(t), \dt \wt u(t))\|_{\vec H^s} = \infty, %\quad \text{}
\label{F1x}
\end{align}

\noi
$\wt \PP$-almost surely.
See also Remarks \ref{REM:uniq1}
and \ref{REM:uniq2}.

From the almost sure convergence 
\eqref{data2zz}
of the extended enhanced data sets, 
the decompositions \eqref{exp5}  
and 
\eqref{exp6}, 
and  
the continuity of the solution
map for \eqref{NLW5x}
stated in 
Proposition~\ref{PROP:LWP1}, 
we obtain 
almost sure local-in-time convergence 
of 
the solution $(\wt u_N, \dt \wt u_N)$ to~\eqref{XNLW1}
 to the solution $(\wt u, \dt \wt u)$ of 
the  stochastic NLW~\eqref{XNLW2}.
By iterating the local-in-time approximation argument
as in \cite[Section 5]{LLOT}, 
we can obtain 
the following almost sure convergence result on every interval $[0,T]$ with
$
0<T<\wt T_{\max}$, 
where $\wt T_{\max}$ denotes the maximal existence time
 for the solution $\wt u$ to \eqref{XNLW2}.

\begin{theorem}
\label{THM:uNu} 

Let $\al \ge \frac 12$.
There exist
a new probability space
$(\wt \O, \wt \F, \wt \PP)$
and
a sequence $\big\{ \big(\wt u_0^{(N)}, \wt u_1^{(N)}\big)\big\}_{N \in \N}$
of $\D'(\T^2) \times \D'(\T^2)$-valued random variables on $\wt \O$
with 
\begin{align*}
\Law \big((\wt u_0^{(N)}, \wt u_1^{(N)})\big) = \mu_{1-\al}\otimes \mu_{-\al}, \quad N \in \N, 
%\label{Law0}
\end{align*}

\noi
where $\mu_s$ is as in \eqref{mu1}, 
such that, 
given a random time $T$ with $0 < T < \wt T_{\max}$, $\wt \PP$-almost surely, 
the solution $(\wt u_N, \dt \wt  u_N)$ to~\eqref{XNLW1}
with the renormalized initial data
$\big(C_{\al, N} \pi_N \wt u_0^{(N)}, C_{\al, N} \pi_N \wt  u_1^{(N)}\big)$
converges $\wt \PP$-almost surely to the solution $(\wt  u, \dt \wt  u)$ of the 
stochastic NLW \eqref{XNLW2}
in $C([0, T]; \vec H^{-\frac 12 - \eps}(\T^2))$
 as $N\to \infty$.
Here,  $\wt T_{\max}$ denotes the maximal existence time
 for the solution $(\wt u, \dt \wt u)$ to \eqref{XNLW2}.

\end{theorem}

We present  the proof of 
Theorem \ref{THM:1}
in Section \ref{SEC:THM}.

\subsection{Organization of the paper}

In Section \ref{SEC:2}, 
we introduce basic notations
and recall\,/\,establish preliminary lemmas.
In particular, in Subsection \ref{SUBSEC:count}, 
we establish the counting lemma
and the angular averaging lemma, 
which in particular play an important role
in the proof of Theorem \ref{THM:ZNZ} presented in Section \ref{SEC:Psecond}.
In Sections \ref{SEC:rand2} and \ref{SEC:para}, 
we establish the basic regularity and convergence
properties
of the stochastic terms and the paracontrolled operators, 
following the argument in \cite{GKO2, OO, OOT1}.
In Section \ref{SEC:conv1}, 
we present the proof of Proposition \ref{PROP:LWP1}
on local well-posedness of the paracontrolled system \eqref{NLW5x}, 
while
we present the proof of Proposition
\ref{PROP:B1}
on local well-posedness and  almost sure finite time blowup
for the limiting equation \eqref{NLW1b}
in Section \ref{SEC:B}.
In Section \ref{SEC:Sk}, 
we justify the application of the Skorokhod representation theorem 
for \eqref{data2zz}-\eqref{Law8}.
In Section \ref{SEC:THM}, 
we first establish almost sure convergence on the new probability space
(Theorem \ref{THM:uNu})
in Subsection \ref{SUBSEC:main1}
and then prove convergence in law on the original probability space
(Theorem \ref{THM:1}) 
in Subsection \ref{SUBSEC:main2}.
In Appendix \ref{SEC:A}, 
 we present the proof of Corollary \ref{COR:2}
on probabilistic ill-posedness.

\section{Preliminaries}
\label{SEC:2}

\subsection{Notations}
\label{SUBSEC:notas}
We write $A\les B$ if $A\le CB$ for some constant $C>0$,
$A\sim B$ if $A\les B\les A$,
and $A\ll B$ if $A\le cB$ for some small constant $c>0$.
We use $\les_\b$  to indicate dependence on $\b$.

For nonzero vectors $a_1,a_2\in\R^2$,
let $\angle(a_1,a_2)\in[0,\pi]$ denote their usual angle,
and let $\Arg(a_1,a_2)\in[0,2\pi)$ denote
 the counterclockwise oriented angle from $a_1$ to $a_2$.

In this paper, we restrict our attention to real-valued
functions.
We write $\ft f$ to denote the spatial or temporal  Fourier transform of $f$
but its meaning is clear from the context.
For $N\in\N$, the Dirichlet projector $\pi_N$ is defined by
\begin{align}
\ft{\pi_Nf}(n)=\ind_{|n|\le N}  \cdot \ft f(n).
\label{Diri1}
\end{align}

Unless otherwise stated, function spaces consist of
real-valued functions or distributions.
We set $\R_+=[0,\infty)$ and $\mathring{\R}_+=(0,\infty)$.
We equip $\D(\T^2)=C^\infty(\T^2)$
and $\D(\mathring{\R}_+\times\T^2)
=C_c^\infty(\mathring{\R}_+\times\T^2)$
with their usual test-function topologies,
and denote their continuous duals by $\D'(\T^2)$
and $\D'(\mathring{\R}_+\times\T^2)$, respectively.

Given a Banach space $B$
of functions\,/\,distributions on $\T^2$, 
 we write
$C_TB=C([0,T];B)$ and $C_IB=C(I;B)$
for a time interval $I\subset\R_+$.
As in Theorem \ref{THM:ZNZ}, 
all path spaces on $\R_+$
are endowed with the compact-open topology in time.

Given Banach spaces $X$ and $Y$, 
we denote by $\L(X;Y)$ the space of bounded linear operators from $X$ to $Y$.

For $0<\be<1$ and $T > 0$, 
let $\C^\be([0, T]; B)$ be the Lipschitz function space, 
defined by the norm:
\begin{align}
\|f\|_{\C_T^\be B}
=\|f\|_{C_TB}+\|f\|_{C_T^\be B}, 
\label{HO2}
\end{align}

\noi
where $\|f\|_{C_T^\be B}$ denotes the  H\"older semi-norm:
\begin{align*}
\|f\|_{C_T^\be B}
=\sup_{0\le t_1 <t_2\le T}
\frac{\|f(t_1)-f(t_2)\|_B}{|t_1-t_2|^\be}.
%\label{HO1}
\end{align*}

We also recall the integral representation
of the Bessel function $J_0$ of order $0$:
\begin{align}
J_0(r)=\frac1{2\pi}\int_0^{2\pi}
\cos(r\sin\ta)d\ta,
\qquad r\in\R,
\label{Bessel1}
\end{align}

\noi
and its asymptotic behavior:
\begin{align}
J_0(r)=\sqrt{\frac2{\pi|r|}}
\cos\Big(|r|-\frac\pi4\Big)
+O(|r|^{-\frac32})
\label{Bessel2}
\end{align}

\noi
as $|r|\to\infty$;
see \cite[Appendices~B.2 and B.8]{Gra}.
In particular, we have $|J_0(r)|\les\jb{r}^{-1/2}$.

\subsection{Sobolev and Besov spaces}
\label{SUBSEC:space}
For $s\in\R$, we define $H^s(\T^2)$ by the norm:
\begin{align*}
\|f\|_{H^s}
=\|\jb{n}^s\ft f(n)\|_{\ell_n^2}.
\end{align*}

\noi
For $1\le p\le\infty$, we define $W^{s,p}(\T^2)$
as the completion of $\D(\T^2)$ under the norm:
\begin{align*}
\|f\|_{W^{s,p}}
=\|\jb{\nb}^sf\|_{L^p}
=\big\|\F^{-1}(\jb{n}^s\ft f(n))\big\|_{L^p}.
\end{align*}

\noi
With this convention, $W^{s,p}(\T^2)$ is separable
even when $p=\infty$.
We have $H^s(\T^2)=W^{s,2}(\T^2)$.

For $j\in\N_0=\N\cup\{0\}$,
 we define the Littlewood-Paley projector $\P_j$ as the
Fourier multiplier operator with a symbol $\varphi_j$.
See \cite[Subsection 3.1]{OOT1}
for the precise definition. 
For $s\in\R$ and $1\le p,q\le\infty$,
we define the Besov space
$B^s_{p,q}(\T^2)$
as the completion of $\D(\T^2)$ under the norm:
\begin{align*}
\|f\|_{B^s_{p,q}}
=\big\|2^{sj}\|\P_jf\|_{L^p}\big\|_{\ell_j^q(\N_0)}.
\end{align*}

\noi
In particular, $H^s=B^s_{2,2}$ with equivalent norms.
Recall that 
\begin{align}
\| u \|_{B^s_{p, \infty}}
\les \| u \|_{W^{s, p}}
   \les \| u \|_{B^{s+\eps}_{p, \infty}}
 \label{besov1}
\end{align}

\noi
for any $s \in \R$, $1 \le p \le \infty$, and $\eps > 0$.

Recall the following paraproduct decomposition due to Bony~\cite{Bony}:
\begin{align}
\begin{split}
fg&=f\pl g+f\pe g+f\pg g\\
&:=\sum_{j<k-2}\P_jf\,\P_kg
+\sum_{|j-k|\le2}\P_jf\,\P_kg
+\sum_{k<j-2}\P_jf\,\P_kg.
\end{split}
\label{para1}
\end{align}

\noi
The first and third terms are paraproducts,
and the second is the resonant product.
We also write $f\pge g=f\pg g+f\pe g$.
We recall the following basic estimates;
see \cite{BCD,MW1}
for details of the proofs in the non-periodic
case (which can be easily extended to the current periodic setting).

\begin{lemma}\label{LEM:para}
 Let $s_1,s_2\in\R$ and
$1\le p,p_1,p_2,q\le\infty$ satisfy
$\frac1p=\frac1{p_1}+\frac1{p_2}$.
Then, we have
\begin{align*}
\|f\pl g\|_{B^{s_2}_{p,q}}
&\les\|f\|_{L^{p_1}}\|g\|_{B^{s_2}_{p_2,q}}.
\end{align*}

\noi
If $s_1<0$, then we have
\begin{align*}
\|f\pl g\|_{B^{s_1+s_2}_{p,q}}
\les\|f\|_{B^{s_1}_{p_1,q}}
\|g\|_{B^{s_2}_{p_2,q}}.
\end{align*}

\noi
If $s_1+s_2>0$, then we have
\begin{align*}
\|f\pe g\|_{B^{s_1+s_2}_{p,q}}
\les\|f\|_{B^{s_1}_{p_1,q}}
\|g\|_{B^{s_2}_{p_2,q}}.
\end{align*}

\end{lemma}

\subsection{Counting estimates
and angular averaging by Abel summation}
\label{SUBSEC:count}

We first collect the lattice point counting estimates used 
in the proof of Theorem \ref{THM:ZNZ} in Section \ref{SEC:Psecond}.

\begin{lemma}\label{LEM:count}
\textup{(i)} 
For dyadic  $M \ge 1$, a unit vector $\nu\in \mathbb{S}^1$,  and $0 \le \ta \le 2\pi$, 
 we have 
\begin{align}
\#\{n\in\Z^2: |n|\sim M,\ \angle(n,\nu)\le\ta\}
\les M+M^2\ta.
\label{Rvol}
\end{align}

%For $r\ge1$, let
%\begin{align}
%Q(r)=\{n\in\Z^2:|n|\le r\}.
%\label{Gauss1}
%\end{align}
%
%\noi
%Then, we have 
%\begin{align}
%\#Q(r)=\pi r^2+O(r^{\frac23}), \qquad \text{as $r \to \infty$.}
%\label{Gauss2}
%\end{align}

\smallskip

\noi
\textup{(ii)} For $r\ge1$ and $a\in\R^2$ with $|a|\ge1$,
we have
\begin{align}
\#\{n\in\Z^2:|n|\le r<|n-a|\}\les r|a|.
\label{Gauss2b}
\end{align}

\smallskip

\noi
\textup{(iii)} For $r\ge 1$, $a\in\R^2\setminus\{0\}$,  and
$0\le\ta_1<\ta_2\le2\pi$, let
\begin{align*}
S(r,a,\ta_1,\ta_2)
=\{n\in\Z^2:0<|n|\le r,
\ \ta_1\le\Arg(a,n)<\ta_2\}.
%\label{Gauss3}
\end{align*}

\noi
Then, we have
\begin{align}
\#S(r,a,\ta_1,\ta_2)
=\frac{\ta_2-\ta_1}{2}r^2+O(r)
\label{Gauss4}
\end{align}

\noi
for $r\ge1$, 
uniformly in $a,\ta_1,\ta_2$.
\end{lemma}

\begin{proof}
As for \eqref{Rvol}, 
 the angular sector has area $O(M^2\ta)$
and boundary length $O(M)$.
Then, counting unit squares that meet this angular sector
yields
\eqref{Rvol}.
For \eqref{Gauss2b}, note that the set is contained in
\begin{align*}
\big\{n\in\Z^2:  \max(r-|a|, 0)<|n|\le r\big\}\cup\{0\}.
\end{align*}

\noi
When $|a|\le r$, 
counting unit squares that meet this annulus gives
$O(r|a|)$.
When $|a|>r$, the bound follows from
$\#\{n\in\Z^2:|n|\le r\}\les r^2\le r|a|$.
The same counting argument proves~\eqref{Gauss4} since
only squares meeting the boundary contribute to the error,
and the boundary has length at most $(2\pi+2)r$.
\end{proof}

In order to handle  the logarithmically divergent sum at 
the threshold case $\al=\frac12$, 
we exploit 
angular averaging 
by  Abel summation.

\begin{lemma}[angular averaging]
\label{LEM:angavg}
Let $a\in\R^2\setminus\{0\}$ and
$\Phi\in C^1([0,2\pi])$.
Then,
\begin{align}
\lim_{N\to\infty}\frac1{\log\jb{N}}
\sum_{1\le|m|\le N}\frac{\Phi(\Arg(a,m))}{|m|^2}
=\int_0^{2\pi}\Phi(\ta)d\ta.
\label{angavg1}
\end{align}
\end{lemma}

\begin{proof}
For $r\ge1$, $a\in\R^2\setminus\{0\}$,  and $0\le\theta\le2\pi$, set $B_a(r,0) = 0$ and 
\begin{align*}
B_a(r,\theta)
=\#\{m\in\Z^2:1 \le|m|\le r,\ 0 \le \Arg(a,m)<\theta\} \quad \text{for $0 < \ta \le 2\pi$}.
\end{align*}

\noi
Then, from  \eqref{Gauss4}, we have 
\begin{align}
B_a(r,\theta)=\frac12\theta r^2+O(r),
\label{ave1}
\end{align}

\noi
uniformly in $\theta$ and $a$.
Let $M=B_a(r,2\pi)$, and enumerate the lattice points in
$\{m\in\mathbb Z^2:1\le |m|\le r\}$ as $m_1,\dots,m_M$.
Then, set $\alpha_j=\Arg(a,m_j)$
for $j=1,\dots,M$.
Note that  we have
\begin{align}
B_a(r,\theta)
=
\sum_{j=1}^M \mathbf 1_{\{\al_j<\theta\}}.
\label{ave2}
\end{align}

\noi
Then, using  \eqref{ave2}, we have
\begin{align}
\begin{split}
A_\Phi(r)
&:=\sum_{1\le|m|\le r}\Phi(\Arg(a,m))
= \sum_{j = 1}^M \Phi(\al_j)\\
& =
\sum_{j=1}^M
\bigg[
\Phi(2\pi)
-\int_0^{2\pi}
\Phi'(\theta)\mathbf 1_{\{\alpha_j<\theta\}}d\theta
\bigg]\\
&=\Phi(2\pi)B_a(r,2\pi)
-\int_0^{2\pi}\Phi'(\theta)B_a(r,\theta)d\theta\\
&=\frac{r^2}{2}\int_0^{2\pi}\Phi(\theta)d\theta
+O(\|\Phi\|_{C^1}\cdot r), 
\end{split}
\label{A60}
\end{align}

\noi
where the last step follows from 
 \eqref{ave1} and integration by parts in $\ta$.

Hence, by 
 applying 
the radial Abel summation (= summation by parts with respect to the radial variable $|m|$) 
with  \eqref{A60},  
 we obtain
\begin{align*}
\begin{split}
\sum_{1\le|m|\le N}\frac{\Phi(\Arg(a,m))}{|m|^2}
&=\frac{A_\Phi(N)}{N^2} %- A_\Phi(1)
+2\int_1^N\frac{A_\Phi(r)}{r^3}dr\\
&=\log N\int_0^{2\pi}\Phi(\theta)d\theta
+O(\|\Phi\|_{C^1}), 
\end{split}
%\label{A61}
\end{align*}

\noi
from which  \eqref{angavg1} follows.
\end{proof}

\subsection{Probabilistic tools}
\label{SUBSEC:probtool}

We use the probabilistic tools collected in
\cite[Subsection~2.3]{LLOT}.
In Subsection \ref{SUBSEC:uniq}, 
we use 
 \cite[Lemma~2.4]{LLOT}
on convergence of distribution-valued random variables,
and the fourth moment theorem \cite[Lemma~2.7]{LLOT}
to prove uniqueness of the limit
for Theorem~\ref{THM:ZNZ}.
We also use 
the Wiener chaos estimate \cite[Lemma~2.5]{LLOT}
and 
\cite[Lemma~2.6]{LLOT}
for establishing 
 convergence of various stochastic terms.

Let $W$ be the isonormal Gaussian process over
$L^2(\{0,1\}\times\T^2)$ given by
\begin{align}
W(f)=\sum_{j=0}^1\sum_{n\in\Z^2}
\ft f_j(n)g_{j,n},
\label{WW1}
\end{align}

\noi
where $f_j=f(j,\cdot)$,
$g_{0,n}=g_n$, and $g_{1,n}=h_n$
with $g_n$ and $h_n$ as in \eqref{data1}.
Here, the underlying Hilbert space is real,
and the sum converges in $L^2(\O)$.
We denote by $\H_k$ the $k$th homogeneous Wiener chaos
associated with $W$ and set
$\H_{\le k}=\bigoplus_{j=0}^k\H_j$.
Let $I_k$ denote the $k$th multiple Wiener integral,
normalized by
\begin{align*}
I_k(f^{\otimes k})=H_k(W(f)),
\qquad \|f\|_{L^2(\{0,1\}\times\T^2)}=1,
%\label{WW2}
\end{align*}

\noi
where $H_k$ is the $k$th Hermite polynomial.
For a general kernel $f$, we set $I_k(f)=I_k(f_{\rm sym})$,
where $f_{\rm sym}$ denotes its symmetrization.
See, for example, \cite[(B.1)]{OWZ}
for the definition of symmetrization.
We say that $f$ is symmetric 
if $f = f_{\rm sym}$.

For symmetric functions $f_j\in L^2
((\{0,1\}\times\T^2)^{k_j})$, $j=1,2$,
and $1\le r\le\min(k_1,k_2)$,
we write
\begin{align}
\begin{split}
&(f_1\otimes_r f_2)(x_1,\dots,x_{k_1+k_2-2r})\\
&=\sum_{i_1,\dots,i_r\in\{0,1\}}
\int_{(\T^2)^r}
f_1(x_1,\dots,x_{k_1-r},
(i_1,y_1),\dots,(i_r,y_r))\\
&\hphantom{=\sum\int}\times
f_2(x_{k_1-r+1},\dots,x_{k_1+k_2-2r},
(i_1,y_1),\dots,(i_r,y_r))
dy_1\cdots dy_r.
\end{split}
\label{N3a}
\end{align}

\noi
The Wiener chaos estimate and the fourth moment theorem
in \cite[Lemmas~2.5 and 2.7]{LLOT}
apply to this isonormal Gaussian process as well;
see also \cite[Theorems~5.2.7 and 6.2.3]{NP}.
The fourth moment theorem 
 was first introduced in \cite{NP05}
 by Nualart and Peccati
 and was  further extended in \cite{NO08, NP09}, 
which plays an important role in proving
uniqueness of the limit of $Z_N$.
See Lemma \ref{LEM:uniq2}.

We also recall the following version
of the Skorokhod representation theorem;
see \cite[Theorem~6.7]{PB71}.
Note that Lemma \ref{LEM:Sk}
is slightly different from the version stated in \cite[Lemma 2.10]{LLOT}.

\begin{lemma}[Skorokhod representation theorem]
\label{LEM:Sk}
Let $(\M,d)$ be a metric space.
Suppose that Borel probability measures
$\{\mu_n\}_{n\in\N}$ on $\M$
converge weakly to a Borel probability measure $\mu$
as $n\to\infty$.
Assume that there exists a separable Borel subset
$S\subset\M$ such that $\mu(S)=1$.
Then, there exist a probability space
$(\wt\O,\wt\F,\wt\PP)$
and $\M$-valued Borel random variables $X_n$ and $X$
such that
\begin{align*}
\Law(X_n)=\mu_n
\qquad\text{and}\qquad
\Law(X)=\mu,
\end{align*}

\noi
and
\begin{align*}
d(X_n,X)\longrightarrow0
\quad\wt\PP\text{-almost surely}.
\end{align*}

\end{lemma}

\section{Convergence of the second Picard iterates}
\label{SEC:Psecond}

In this section, we prove  Theorem~\ref{THM:ZNZ}.
The normalization in \eqref{CN} makes the limiting covariance
in \eqref{CZNd} equal to $J_0((t_1-t_2)|n|)$
as in \eqref{nxib}.
We note that the Gaussian processes $\ze_n$ in \eqref{nxib} are well defined.
Indeed, from \eqref{Bessel1},  we have 
\begin{align*}
\sum_{j,k=1}^m c_j\cj{c_k}\, J_0((t_j-t_k)|n|)
=\frac1{2\pi}\int_0^{2\pi}
\bigg|\sum_{j=1}^m c_je^{it_j|n|\cos\theta}\bigg|^2
 d\theta\ge0
%\label{A63}
\end{align*}

\noi
for $c_j\in\mathbb C$ and $t_j\in\R$,
showing that  the kernel is positive definite.
Then, we construct $\{\ze_n\}_{n \in \Z^2}$ 
by taking independent complex-valued Gaussian processes
(with one representative of each pair $\{n,-n\}$, $n\ne0$)
with independent real and imaginary parts,
each with covariance $\frac12J_0((t-r)|n|)$,
and impose that $\ze_{-n} = \cj{\ze_n}$, $n \in \Z^2\setminus\{0\}$.
The zero mode $\ze_0$ is a real-valued standard Gaussian, constant in time.
Moreover, for each fixed $n$,
it follows from \eqref{Bessel1} that 
\begin{align*}
\E\big[|\ze_n(t)-\ze_n(r)|^2\big]
=2\big(1-J_0((t-r)|n|)\big)
\les |n|^2|t-r|^2,
\end{align*}

\noi
and thus  these processes admit continuous versions.
The fourth moment argument in Subsection~\ref{SUBSEC:uniq} 
below shows that the Duhamel integral of these Gaussian processes 
(namely, $Z$ in~\eqref{Z1}) is the limit in law of $Z_N$.

In Subsection \ref{SUBSEC:Riea}, we establish two preliminary
Riemann sum approximations
(Lemmas~\ref{LEM:COVa} and~\ref{LEM:COVb})
for the nearly resonant contributions.
In Subsection \ref{SUBSEC:uniq},
using Lemmas~\ref{LEM:COVa} and~\ref{LEM:COVb}
together with the fourth moment theorem,
we establish uniqueness of the limit.
In Subsection~\ref{SUBSEC:tight},
by exploiting multilinear dispersive smoothing for $Z_N$,
we establish tightness of $\{Z_N\}_{N\in\N}$.

\subsection{Riemann sum approximations}
\label{SUBSEC:Riea}

We first record the normalization in \eqref{CN}:
\begin{align}
C_{\al,N}^4
\sum_{|m|\le N}\jb{m}^{4\al-4}=1.
\label{K1}
\end{align}

\noi
In particular, by Cauchy-Schwarz's inequality, we have
\begin{align}
C_{\al,N}^4
\sum_{\substack{n=n_1+n_2\\|n_1|,|n_2|\le N}}
\frac{1}{\jb{n_1}^{2-2\al}\jb{n_2}^{2-2\al}}
\le 1, 
\label{K2}
\end{align}

\noi
uniformly in $n\in\Z^2$ and $N\in\N$.

In the next two lemmas, we establish
the key identities
for the limiting covariance.
In Lemma \ref{LEM:COVa}, we treat the case $\al > \frac 12$, 
while we treat the threshold case $\al = \frac 12$
in Lemma \ref{LEM:COVb}.

\begin{lemma}
\label{LEM:COVa}
Let $\al>\frac12$.
Then, for any $n\in\Z^2$ and $t\in\R$, we have
\begin{align}
\lim_{N\to\infty}C_{\al,N}^4
\sum_{\substack{n=n_1+n_2\\|n_1|,|n_2|\le N}}
\frac{\cos(t(\jb{n_1}-\jb{n_2}))}
{\jb{n_1}^{2-2\al}\jb{n_2}^{2-2\al}}
=J_0(t|n|),
\label{COVa1}
\end{align}

\noi
where $J_0$ is the Bessel function of order $0$ in \eqref{Bessel1}.
\end{lemma}

\begin{proof}
Denote the expression inside the limit in \eqref{COVa1} 
by $\A_{N,1}(n)$.
When $n=0$, it follows from \eqref{K1} and \eqref{Bessel1} that 
\begin{align*}
\A_{N,1}(0)=1=J_0(0)
%\label{COVa11}
\end{align*}

\noi
for every $N\in\N$.
Thus, we fix $n\ne0$ in the following.

Let  $|m|\ge 2|n|+2$.
From  the mean value theorem, we have 
\begin{align}
\begin{split}
\bigg|\jb{n-m}-\jb{m}+\frac{n\cdot m}{|m|}\bigg|
&\les_n |m|^{-1},\\
\bigg|
\frac{1}{\jb{n-m}^{2-2\al}\jb{m}^{2-2\al}}
-\frac{1}{|m|^{4-4\al}}\bigg|
&\les_{\al,n}|m|^{4\al-5}.
\end{split}
\label{K3}
\end{align}

\noi
See also \cite[(4.13)]{GKO2}.
From Lemma \ref{LEM:count}\,(ii), we have 
\begin{align}
\#\{m\in\Z^2:|m|\le N<|n-m|\}\les_n N.
\label{K3a}
\end{align}

\noi
Define $\A_{N,2}(n)$ by
\begin{align}
\A_{N,2}(n)
=C_{\al,N}^4\sum_{1\le |m|\le N}
\frac{\cos(tn\cdot m/|m|)}{|m|^{4-4\al}}.
\label{COVa8}
\end{align}

\noi
We now study the difference  $\A_{N,1}(n) - \A_{N,2}(n)$.
We first note that 
the contribution from  $|m|<2|n|+2$ is $O_{\al, n}(1)$ before multiplication by
$C_{\al,N}^4$.
For $|m|\ge 2|n|+2$, it follows from~\eqref{K3} and
$|\cos x-\cos y|\le |x-y|$ that the error in each summand
is bounded by
\begin{align}
\bigg|
\frac{\cos(t(\jb{n-m}-\jb{m}))}
{\jb{n-m}^{2-2\al}\jb{m}^{2-2\al}}
-\frac{\cos(tn\cdot m/|m|)}
{|m|^{4-4\al}}
\bigg|
\les_{\al,t, n}|m|^{4\al-5}.
\label{K3b}
\end{align}

\noi
With $n_2 = m$ (and thus $n_1 = n - m$), 
the summation for $\A_{N, 1}$ is over the set
$E_{N, 1} = \{|m|\le N,\ |n-m|\le N\}$, 
while 
the summation for $\A_{N, 2}$ is over the set
$E_{N, 2} = \{1\le |m|\le N\}$.
Using~\eqref{K3a}, we have  
\begin{align}
\sum_{m \in E_{N, 2} \setminus E_{N, 1}}
\frac{|\cos(tn\cdot m/|m|)|}{|m|^{4-4\al}}
\les_{\al,  n} N^{4\al-3}
\label{K3c}
\end{align}

\noi
for sufficiently large $N \gg 1$, 
where we used the fact that $|m|\sim N$
on $E_{N, 2} \setminus E_{N, 1}$
for $N \gg1 $ (recall that $n$ is fixed).
Hence, by \eqref{K3b} and \eqref{K3c} with \eqref{CN}, we obtain
\begin{align}
\begin{split}
|\A_{N,1}(n)-\A_{N,2}(n)|
&\les_{\al, t, n} C_{\al,N}^4
\bigg(1+ \sum_{1\le |m|\le N}|m|^{4\al-5}
+N^{4\al-3} \bigg)\\
&\les_{\al}
\begin{cases}
N^{2-4\al},&\frac12<\al<\frac34,\\
N^{-1}\log\jb{N},&\al=\frac34,\\
N^{-1},&\al>\frac34, 
\end{cases}
\end{split}
\label{K4}
\end{align}

\noi
where the right-hand side tends to $0$ as $N \to \infty$.

Next, we  determine the limit of $\A_{N,2}(n)$
by viewing it as a Riemann sum approximation of 
the integral of the function
\begin{align}
F(x)
=
\frac{\cos(tn\cdot x/|x|)}{|x|^{4-4\alpha}},
\quad 0<|x|\le1, 
\label{K5a}
\end{align}

\noi
with $x = \frac m N$.
Let $0<\eta<1$. The contribution 
to \eqref{COVa8}
(with $C_{\al,N}^4$ replaced by $N^{2-4\al}$)
from $1\le |m|\le\eta N$
is bounded by
\begin{align}
N^{2-4\al}
\sum_{1\le |m|\le\eta N}|m|^{4\al-4}
\les_\al \eta^{4\al-2}, 
\label{K5}
\end{align}

\noi
uniformly in $N \in \N$.
%which tends to $0$ as $\eta \to 0$.
On the region $\eta\le|x|\le1$, we apply  the usual
Riemann sum approximation:
\begin{align}
\lim_{N \to \infty} 
N^{2-4\al}
\sum_{\eta N \le |m|\le N}
\frac{\cos(tn\cdot m/|m|)}{|m|^{4-4\al}}
= 
\int_{\eta < |x|\le1}
\frac{\cos(tn\cdot x/|x|)}{|x|^{4-4\al}}dx.
\label{COVa9x}
\end{align}

\noi
Hence, by noting that the integral of $F$ in \eqref{K5a}
on $\{0 < |x| \le \eta\}$ is $O(\eta^{4\al - 2})$
(which tends to $0$ as $\eta \to 0$)
and by first letting $N\to\infty$ with $\eta>0$ fixed and then letting $\eta\to0$, 
it follows from \eqref{K5}, \eqref{COVa9x}, 
and 
the polar coordinates
that 
\begin{align}
\begin{split}
\lim_{N \to \infty}
 N^{2-4\al}
\sum_{1\le |m|\le N}
\frac{\cos(tn\cdot m/|m|)}{|m|^{4-4\al}}
&
= 
\int_{0 < |x|\le1}
\frac{\cos(tn\cdot x/|x|)}{|x|^{4-4\al}}dx\\
&=\frac{1}{4\al-2}
\int_0^{2\pi}\cos(t|n|\cos\theta)d\theta\\
&=\frac{\pi}{2\al-1}J_0(t|n|), 
\end{split}
\label{COVa9}
\end{align}

\noi
where we use \eqref{Bessel1} in the last step.
By repeating  the same argument with
\begin{align}
\big|\jb{m}^{4\alpha-4}
- 
|m|^{4\alpha-4}\big|\les_\al |m|^{4\alpha-6}
\label{K6x}
\end{align}

\noi
for $|m|\ge 1$
and 
 the cosine replaced by $1$,
we obtain
\begin{align}
\lim_{N \to \infty}N^{2-4\al}C_{\al,N}^{-4}
=
\lim_{N \to \infty} N^{2-4\al}\sum_{|m|\le N}\jb{m}^{4\al-4}
=  \frac{\pi}{2\al-1}.
\label{K6}
\end{align}

\noi
Hence, 
from \eqref{COVa8}, \eqref{COVa9}, and \eqref{K6}, 
we obtain 
\begin{align}
\lim_{N \to \infty} \A_{N,2}(n) =  J_0(t|n|).
\label{K6a}
\end{align}

\noi
Therefore, from 
 \eqref{K4} and \eqref{K6a}, 
 we conclude that 
$\A_{N,1}(n)$ converges to $ J_0(t|n|)$ 
as $N \to \infty$, proving \eqref{COVa1}.
\end{proof}

Next, we treat the threshold case $\al = \frac 12$.

\begin{lemma}
\label{LEM:COVb}
For any $n\in\Z^2$ and $t\in\R$, we have
\begin{align}
\lim_{N\to\infty}C_{\frac12,N}^4
\sum_{\substack{n=n_1+n_2\\|n_1|,|n_2|\le N}}
\frac{\cos(t(\jb{n_1}-\jb{n_2}))}
{\jb{n_1}\jb{n_2}}
=J_0(t|n|), 
\label{COVb1}
\end{align}

\noi
where $J_0$ is the Bessel function of order $0$ in \eqref{Bessel1}.
\end{lemma}

\begin{proof}
The case $n=0$ follows from \eqref{K1}.
Fix $n\ne0$ and 
denote the expression inside the limit in \eqref{COVb1} 
by $\B_{N,1}(n)$.
Define $\B_{N,2}(n)$ by
\begin{align}
\B_{N,2}(n)
=C_{\frac 12 ,N}^4\sum_{1\le |m|\le N}
\frac{\cos(tn\cdot m/|m|)}{|m|^{2}}.
\label{COVb2}
\end{align}

\noi
Then, 
by \eqref{K3}, \eqref{K3a}, and the same argument as in the proof of \eqref{K4}, we have 
\begin{align}
\begin{split}
|\B_{N,1}(n) - \B_{N,2}(n)|
&\les_{t, n}C_{\frac12,N}^4
\bigg(1+\sum_{m\ne0}|m|^{-3}+N^{-1}\bigg)\\
&\les (\log\jb{N})^{-1}.
\end{split}
\label{COVb3}
\end{align}

From  Lemma~\ref{LEM:angavg} with $a=n$ and
$\Phi(\theta)=\cos(t|n|\cos\theta)$, 
and  \eqref{Bessel1}, we have 
\begin{align}
\lim_{N \to \infty}
\frac{1}{\log\jb{N}}\sum_{1\le|m|\le N}
\frac{\cos(tn\cdot m/|m|)}{|m|^2}
=  2\pi J_0(t|n|).
\label{K7}
\end{align}

\noi
By 
applying 
\eqref{K6x} with $\al = \frac 12$ and then 
Lemma~\ref{LEM:angavg} 
with $\Phi\equiv1$, we have 
\begin{align}
\lim_{N \to \infty}\frac{C_{\frac12,N}^{-4}}{\log\jb{N}}
=\lim_{N \to \infty}\frac{1}{\log\jb{N}}\sum_{|m|\le N}\jb{m}^{-2}
=  2\pi.
\label{K8}
\end{align}

\noi
Therefore,
dividing \eqref{K7} by \eqref{K8} and using \eqref{COVb3} 
with \eqref{COVb2},
we obtain
 \eqref{COVb1}.
\end{proof}

\subsection{Uniqueness of the limit}
\label{SUBSEC:uniq}

We first record the regularity of the limiting process
$Z$ defined in \eqref{Z1}.
For $h\in\R$, we use $\dl_h$ to denote 
the following difference operator:
\begin{align}
\dl_h f (t)=f(t+h)-f(t).
\label{dl1}
\end{align}

\begin{lemma}\label{LEM:randZ}
Let $\al \ge \frac 12$ and $Z$ be as in \eqref{Z1}.
Then, given any $s < \frac 14$, we have
\begin{align*}
Z\in C(\R_+;W^{s,\infty}(\T^2))
\cap C^1(\R_+;W^{s-1,\infty}(\T^2)), 
\end{align*}

\noi
almost surely.
\end{lemma}

\begin{proof}
We only discuss the regularity of $Z$
since the regularity of $\dt Z$
follows from a similar consideration.

For $n\ne0$, it follows
from  \eqref{Z1}, \eqref{nxib}, and \eqref{Bessel2}
that 
\begin{align*}
\begin{split}
\E\big[|\ft Z(t,n)|^2\big]
&=\int_0^t\int_0^t
\frac{\sin((t-t_1)\jb{n})\sin((t-t_2)\jb{n})}
{\jb{n}^2}J_0((t_1-t_2)|n|)dt_2dt_1\\
&\les\jb{n}^{-\frac52}
\int_0^t\int_0^t|t_1-t_2|^{-\frac12}dt_2dt_1\\
&\les t^{\frac32}\jb{n}^{-\frac52}.
\end{split}
%\label{randZ1a}
\end{align*}

\noi
See \eqref{Z3}.  On the other hand, for $n=0$, 
from \eqref{nxib} and 
\eqref{Bessel1}, we have 
$\ft Z(t,0)=(1-\cos t)\ze_0(0)$.
Moreover, 
it follows from  the mean value theorem
(see, for example, the end of the proof of \cite[Proposition 1.6]{GKO2})
that,  
for small $\ta > 0$, we have 
\begin{align*}
\begin{split}
\E\big[|\dl_h \ft Z(t,n)|^2\big]
&\les_T
\jb{n}^{-\frac52+\ta}|h|^{\ta}
\end{split}
%\label{dtZb}
\end{align*}

\noi
for any $0 \le t \le T$ and $h \in \R$
such that $t + h \in [0, T]$.
Then, given $s <  \frac 14$, 
by 
taking small $ \ta> 0$ sufficiently small, 
the claim follows from 
\cite[Lemma~2.6]{LLOT};
see also \cite[Proposition~2.7]{OOT}.
\end{proof}

In the following proposition, 
we establish 
 convergence of the finite-dimensional
marginals of $Z_N$ in time, 
which will allow us to 
 identify the limit of $Z_N$.

\begin{proposition}
\label{PROP:uniq}
Let $\al\ge\frac12$, and let $Z_N$ and $Z$
be as in \eqref{ZN} and \eqref{Z1},
respectively.
Given $m \in \N$, let
 $0\le t_1\le  \cdots \le t_m$.
Then,
$\big\{Z_N(t_j)\big\}_{j = 1}^m$
converges in law
to
$\big\{Z(t_j)\big\}_{j = 1}^m$
in $(\D'(\T^2))^m$
as $N \to \infty$.
\end{proposition}

We first 
study the space-time covariance of 
$Z_N$.

\begin{lemma}
\label{LEM:COV}

Let $\al\ge\frac12$, and let $Z_N$ and $Z$
be as in \eqref{ZN} and \eqref{Z1}, respectively.
Then,
given any test function $\psi \in \D(\T^2)$, we have
\begin{align}
\E\big[\jb{Z_N (t), \psi}\big]  =
\E\big[\jb{Z (t), \psi}\big] = 0
\label{con1}
\end{align}

\noi
for any $t \ge 0$ and $N \in \N$,
and
\begin{align}
\lim_{N\to \infty}
\E\big[\jb{Z_N(r), \psi_1} \jb{ Z_N (t), \psi_2}\big]
=
\E\big[\jb{ Z(r), \psi_1}\jb{Z (t), \psi_2}\big]
\label{con2}
\end{align}

\noi
for any  test functions $\psi_1, \psi_2 \in \D(\T^2)$
and $t \ge r \ge 0$,
where
$\jb{f, g}$
denotes the $\D'(\T^2)$-$\D(\T^2)$ duality pairing.

\end{lemma}

\begin{proof}

Given $N\in\N$, $n \in \Z^2$, and $t, t' \ge 0$, define
\begin{align}
\begin{split}
\s_{N,n}(t,t')
&=\E\big[\ft z_N(t,n)\cj{\ft z_N(t',n)} \big]
=\ind_{|n|\le N}\cdot C_{\al,N}^2
\frac{\cos((t-t')\jb{n})}{\jb{n}^{2-2\al}}.
\end{split}
\label{sigman}
\end{align}

\noi
The renormalization constant $\s_N$ in \eqref{sigN} satisfies
\begin{align}
\s_N=\sum_{n\in\Z^2}\s_{N,n} (t, t).
\label{sigmane}
\end{align}

\noi
Then, 
from \eqref{ZN} and \eqref{sigmane}, we have
\begin{align}
\begin{split}
\ft Z_N(t,n)
&=\int_0^t S_n(t-t')\,\ft{\wick{z_N^2}}\,(t',n)dt',\\
\ft{\wick{z_N^2}}\, (t',n)
&=\sum_{n_1+n_2=n}
\big(\ft z_N(t',n_1)\ft z_N(t',n_2)
-\ind_{n=0}\cdot \s_{N,n_1}(t', t')\big), 
\end{split}
\label{fZN}
\end{align}

\noi
where $S_n(t)$ is given by 
\begin{align}
S_n(t)=\frac{\sin(t\jb{n})}{\jb{n}}.
\label{SN1}
\end{align}

In particular, we have 
\begin{align}
\jb{Z_N(t),\psi}
=\sum_{n\in \Z^2} \cj{\ft\psi(n)}
\int_0^t S_n(t-t') \,\ft{\wick{z_N^2}}\,(t',n)dt'.
\label{C2}
\end{align}

\noi
By noting that 
 $\E\big[\,\ft{\wick{z_N^2}}\,(t',n)\big]=0$, 
 we obtain 
 \eqref{con1} for $Z_N$.
The corresponding assertion for $Z$ in \eqref{con1}
follows from \eqref{Z1} and the fact that $\E[\ze_n(t')] = 0$.

By Wick's theorem (\cite[Proposition I.2]{Simon}) and \eqref{sigman},
we have 
\begin{align}
\E\Big[ \,\ft{\wick{z_N^2}}\,(t_1,n)\cj{\,\ft{\wick{z_N^2}}\,(t_2,m)}\Big]
&=\ind_{n=m}\cdot K_N(t_1-t_2,n),
\label{K9}
\end{align}

\noi
where $K_N(\tau,n)$ is given by 
\begin{align}
K_N(\tau,n)
&=2C_{\al,N}^4
\sum_{\substack{n_1+n_2=n\\|n_1|,|n_2|\le N}}
\frac{\cos(\tau\jb{n_1})\cos(\tau\jb{n_2})}
{\jb{n_1}^{2-2\al}\jb{n_2}^{2-2\al}}.
\label{K9a}
\end{align}

\noi
See \cite[(3.11)]{GKO2}.
Thus, for real-valued $\psi_1,\psi_2\in\D(\T^2)$,
we have 
\begin{align}
\begin{split}
&\E\big[\jb{Z_N(r),\psi_1}\jb{Z_N(t),\psi_2}\big]\\
&\quad =\sum_{n \in \Z^2} \cj{\ft\psi_1(n)}\ft\psi_2(n)
\int_0^r\int_0^t S_n(r-t_1)S_n(t-t_2)
 K_N(t_1-t_2,n)dt_2dt_1.
\end{split}
\label{CZN}
\end{align}

\noi
Compare this with \cite[(3.7)]{LLOT}.
In the case of BBM studied in \cite{LLOT}, 
we exploited the fast decay of the dispersive symbol 
(but no oscillation was used).
For the current NLW problem, 
we exploit dispersion for the highly oscillatory part, 
while we apply 
Lemmas \ref{LEM:COVa} and \ref{LEM:COVb} 
(based on the angular averaging via Lemma \ref{LEM:angavg}) to treat 
the nearly resonant contributions.

Using the trigonometric identity, 
we write 
\begin{align}
 K_N(\tau,n) &= K_N^+(\tau,n)
 + K_N^-(\tau,n),
\label{CZNb}
\end{align}

\noi
where $K_N^\pm(\tau,n)$ is given by 
\begin{align} 
K_N^\pm(\tau,n)
&=C_{\al,N}^4
\sum_{\substack{n_1+n_2=n\\|n_1|,|n_2|\le N}}
\frac{\cos(\tau(\jb{n_1}\pm\jb{n_2}))}
{\jb{n_1}^{2-2\al}\jb{n_2}^{2-2\al}}.
\label{CZNb2}
\end{align}

We first study the contribution from 
$K_N^+$.
Fix $n\in\Z^2$ and $0\le t_1,t\le T$.
For $a>0$, integration by parts with \eqref{SN1} yields
\begin{align}
\bigg|\int_0^t S_n(t-t_2)
\cos((t_1-t_2)a)dt_2\bigg|
\le\frac{|S_n(t)|+\int_0^t|S_n'(t-t_2)|dt_2}{a}
\le\frac{2T}{a}.
\label{SN2}
\end{align}

\noi
Since $n \in \Z^2$ is fixed, 
it follows from 
 \eqref{CN} that the contribution 
 from $\min(|n_1|,|n_2|)\les\jb{n}$
 to~\eqref{CZNb2} vanishes as $N \to \infty$.
 Thus, we assume 
$\min(|n_1|,|n_2|)\gg\jb{n}$
in the following.  In particular, we have $\jb{n_1} \sim \jb{n_2}$.
Now, 
applying \eqref{SN2} with 
$a = \jb{n_1}+ \jb{n_2}$ and  \eqref{CN},
we obtain
\begin{align}
\begin{split}
\bigg|\int_0^t S_n(t-t_2) K_N^+(t_1-t_2,n)dt_2\bigg|
&\les_{T, n}C_{\al,N}^4
\bigg(1+\sum_{|m|\le N}\jb{m}^{4\al-5}\bigg)\\
&\too 0, 
\end{split}
\label{CZNc}
\end{align}

\noi
as $N \to \infty$.

On the other hand, 
from Lemmas \ref{LEM:COVa} and \ref{LEM:COVb}, we have 
\begin{align}
\lim_{N \to \infty}  K_N^-(\tau,n) = J_0(\tau|n|).
\label{CZNd}
\end{align}

\noi
Moreover, from  \eqref{K2}, we have $|K_N^\pm(\tau,n)|\le1$.
In particular, the summand in \eqref{CZN}
is bounded by
\begin{align*}
2T^2\jb{n}^{-2}|\ft\psi_1(n)||\ft\psi_2(n)|,
\end{align*}

\noi
which is summable in $n \in \Z^2$.
Therefore, by the Lebesgue dominated convergence 
theorem
with \eqref{CZNc},  \eqref{CZNd}, 
\eqref{Z1},  and \eqref{nxib},
we obtain 
\begin{align*}
\begin{split}
&\lim_{N\to\infty}
\E\big[\jb{Z_N(r),\psi_1}\jb{Z_N(t),\psi_2}\big]\\
&\quad =\sum_n\cj{\ft\psi_1(n)}\ft\psi_2(n)
\int_0^r\int_0^t S_n(r-t_1)S_n(t-t_2)
J_0((t_1-t_2)|n|)dt_2dt_1\\
& \quad =
\E\big[\jb{ Z(r), \psi_1}\jb{Z (t), \psi_2}\big].
\end{split}
%\label{CZNe}
\end{align*}

\noi
This proves \eqref{con2}.
\end{proof}

In view of  \cite[Lemma~2.4]{LLOT}
on the characterization of convergence in law
of the finite-dimensional marginal of a distribution-valued
process, 
Proposition~\ref{PROP:uniq}
follows once we establish  the next lemma.

\begin{lemma}
\label{LEM:uniq2}
Let $\al\ge\frac12$, and let $Z_N$ and $Z$
be as in \eqref{ZN} and \eqref{Z1}, respectively.
Given $m\in\N$, let
 $0\le t_1\le  \cdots \le t_m$.
Then, given any test functions
$\{\psi_j\}_{j = 1}^m \subset \D(\T^2)$,
$\big\{\jb{Z_N(t_j), \psi_j}\big\}_{j = 1}^m$
converges in law
to the Gaussian vector
$\big\{\jb{Z(t_j), \psi_j}\big\}_{j = 1}^m$
as $N \to \infty$.
\end{lemma}

\begin{proof}
Fix $t\ge0$ and $\psi\in\D(\T^2)\setminus\{0\}$.
From 
\eqref{C2} with \eqref{fZN}, 
we see that  $\jb{Z_N(t),\psi}$
is an element in  the second homogeneous Wiener chaos.
With the convention \eqref{WW1}, we have
\begin{align}
\jb{Z_N(t), \psi}
= I_2(\fn^t),
\label{B1}
\end{align}

\noi
where $\fn^t = \fn^t((i_1, x_1), (i_2, x_2)) \in
L^2((\{0,1\}\times\T^2)^2)$ is
a symmetric function
whose Fourier transform is given by
\begin{align}
\begin{split}
\F_x ({\fn^t})((i_1,n_1), (i_2,n_2))
&
 =
 \ind_{ |n_1|, |n_2|\leq N}
 \cdot
 C_{\al, N}^2
 \int^t_0
 \frac{\sin((t-t')\jb{n_1+ n_2})}{\jb{n_1+ n_2}} \\
& \quad
\times
\frac{\cj{\ft \psi(n_1+n_2)}}{\jb{n_1}^{1-\al} \jb{n_2}^{1-\al}}
J_{n_1,n_2}^{(i_1,i_2)}(t') dt'
\end{split}
\label{B2}
\end{align}

\noi
for $n_1, n_2\in \Z^2$.
Here,  $J_{n_1,n_2}^{(i_1,i_2)}$ is defined by
\begin{align}
\begin{split}
J_{n_1,n_2}^{(0, 0)} (t)
& = \cos(t\jb{n_1})\cos(t\jb{n_2}), \qquad
 J_{n_1,n_2}^{(0, 1)}(t)
 =
 \cos(t\jb{n_1}) \sin(t\jb{n_2}), \\
 J_{n_1,n_2}^{(1, 0)}(t)
& = \sin(t\jb{n_1}) \cos(t\jb{n_2}),
\qquad
 J_{n_1,n_2}^{(1, 1)}(t)
 =
\sin(t\jb{n_1}) \sin(t\jb{n_2}).
\end{split}
\label{B3a}
\end{align}

As in the proof of \cite[Lemma~3.4]{LLOT}, 
our goal is to show
\begin{align}
\lim_{N \to \infty} \|\fn^t\otimes_1\fn^t\|_{L^2((\{0, 1\}\times \T^2)^2)} = 0, 
\label{SN3}
\end{align}

\noi
where the contraction $\fn^t\otimes_1\fn^t$ is as in \eqref{N3a}, 
and apply 
the fourth moment theorem
\cite[Lemma~2.7]{LLOT}.
Let $H = L^2(\{0,1\}\times\T^2)$
and regard
\[
\fn^t\in L^2((\{0,1\}\times\T^2)^2) \simeq  H\otimes H\]

\noi
as the kernel of the following operator $\TT_N^t:H\to H$:
\[(\TT_N^t g)({\bf x})
= \int_{\{0,1\}\times\T^2} \fn^t({\bf x}, {\bf y})g({\bf y}) d{\bf y},\]

\noi
where the integral here includes the sum over the \(\{0,1\}\) index.
Since $\fn^t$ is symmetric, we have 
\[(\fn^t\otimes_1 \fn^t)({\bf x},{\bf y})
=
\int_{\{0,1\}\times\T^2} \fn^t({\bf x},{\bf z})\fn^t({\bf y},{\bf z}) d{\bf z},  \]

\noi
which 
 is precisely the kernel of $\TT_N^t(\TT_N^t)^* = (\TT_N^t)^2$
 (here we used the fact that $\fn^t$ is real-valued
 and symmetric, which implies that $\TT_N^t$ is self-adjoint).
Consequently, we have 
\begin{align*}
\|\fn^t\otimes_1\fn^t\|_{L^2}
=
\|(\TT_N^t)^2\|_{\HS}.
\end{align*}

\noi
Together 
with the following standard operator inequality:
\[
\|AB\|_{\HS}
\le
\|A\|_{\mathrm{op}}\|B\|_{\HS}, 
\]

\noi
we have
\begin{align}
\|\fn^t\otimes_1\fn^t\|_{L^2}
\le
\|\TT_N^t\|_{\op}
\|\TT_N^t\|_{\HS}
=
\|\TT_N^t\|_{\op}
\|\fn^t\|_{L^2}.
\label{SN4}
\end{align}

From  the Wiener isometry
(see, for example, \cite[Definition 2.9]{CLO}), \eqref{B1},  \eqref{CZN}
(with \eqref{SN1}, \eqref{CZNb}, and \eqref{CZNb2}), 
 and \eqref{K2}, we have 
\begin{align}
2\|\fn^t\|_{L^2}^2
=\E \big[  |\jb{Z_N(t),\psi}|^2  \big]
\les t^2\sum_{n \in \Z^2} \frac{|\ft\psi(n)|^2}{\jb{n}^2}.
\label{SN5}
\end{align}

\noi
By \eqref{B2}, \eqref{B3a},  and Peetre's inequality,
\begin{align*}
\begin{split}
|\F_x(\fn^t)((i,n),(j,m))|
&\le 
\ind_{|n|,|m|\le N}
\cdot
t\, C_{\al,N}^2
\jb{n}^{\al-1}\jb{m}^{\al-1}
\frac{|\ft\psi(n+m)|}{\jb{n+m}}\\
&\les_\al 
\ind_{|n|,|m|\le N}\cdot     t\, C_{\al,N}^2 
\jb{n}^{2(\al-1)}
b(n+m), 
\end{split}
%\label{K10}
\end{align*}

\noi
where $b(n)$ is given by
\begin{align*}
b(n)=\jb{n}^{|\al-1|-1}|\ft\psi(n)|.
\end{align*}

\noi
The same estimate holds with $n$ and $m$ interchanged.
Since $\psi\in\D(\T^2)$, we have $b\in\l^1(\Z^2)$.
Hence, it follows from Schur's test with \eqref{CN} that 
\begin{align}
\|\TT_N^t\|_{\op}
\les_{\al, t}
C_{\al,N}^2\|b\|_{\l^1}\max_{|n|\le N}\jb{n}^{2\al-2}
\les
\begin{cases}
(\log\jb{N})^{-\frac12},&\al=\frac12,\\
N^{1-2\al},&\frac12<\al\le1,\\
N^{-1},&\al>1, 
\end{cases}
\label{K11}
\end{align}

\noi
which tends to $0$ as $N \to \infty$.
Therefore, \eqref{SN3} follows from \eqref{SN4}, 
\eqref{SN5},  and \eqref{K11}.

Finally, from 
the vector-valued fourth moment theorem
\cite[Lemma~2.7]{LLOT} with \eqref{SN3}
and Lemma~\ref{LEM:COV},
we conclude that 
$\big\{ \jb{Z_N(t_j), \psi_j}\big\}_{j=1}^m
= \big\{I_2(f_{N,\psi_j}^{t_j})\big\}_{j=1}^m$
converges in law to 
$\big\{\jb{Z(t_j),\psi_j}\big\}_{j=1}^m$
as $N \to \infty$.
This concludes the proof of  Lemma~\ref{LEM:uniq2}.
\end{proof}

\subsection{Tightness via multilinear dispersive smoothing}
\label{SUBSEC:tight}

\begin{proposition}\label{PROP:rand1}
Let $\al\ge\frac12$ and $Z_N$ be as in \eqref{ZN}.
Then, for any $s<\frac14$, the sequence $\{Z_N\}_{N\in\N}$
is tight in
$C(\R_+;W^{s,\infty}(\T^2))
\cap C^1(\R_+;W^{s-1,\infty}(\T^2))$.

\end{proposition}

Theorem \ref{THM:ZNZ}
follows from 
Proposition~\ref{PROP:rand1}
and
 the Prokhorov theorem \cite[Lemma~2.9]{LLOT}
together with 
Proposition~\ref{PROP:uniq}
on  uniqueness of the limit.
We omit details.  See the proof of \cite[Theorem 1.6]{LLOT}.

\begin{proof}[Proof of Proposition~\ref{PROP:rand1}]
We follow the proof of \cite[Proposition~1.6]{GKO2}
and the two-dimensional computation in
\cite[Subsection~4.2]{OO}, 
where multilinear smoothing for the second Picard iterates
was established.

Fix $T>0$.
Our first goal is to prove the following bound:
\begin{align}
 \E\big[|\dt^j\ft Z_N(t,n)|^2\big]
\les_T \jb{n}^{2j -\frac 52}
\label{J0}
\end{align}

\noi
for any $0 \le t \le T$ and $j=0,1$, uniformly in $N \in\N$.
Expanding the sine and cosine factors in~\eqref{CZN}
(with $\psi_1 = \psi_2 = e_n$
and the complex conjugate on $\jb{Z_N(t),\psi_2}$)
and integrating in one time variable as in
\cite[(4.17)-(4.20)]{OO}, we obtain\footnote{The diagonal contribution from $n_1 = n_2 = \frac n2$
needs to be treated separately, but 
one can check that it satisfies \eqref{J1}.
A similar comment applies
to the terms analogous to $I^{(k_1,k_2)}$, $(k_1, k_2) \ne (1, 1)$ appearing in \cite[(4.14)]{OO}. }
\begin{align}
\begin{split}
& \E\big[|\dt^j\ft Z_N(t,n)|^2  \big]\\
& \quad \les_T C_{\al,N}^4 \jb{n}^{2j-2}
\sum_{ \eps_0, \eps_1, \eps_2\in\{-1,1\}}
\sum_{\substack{n=n_1+n_2\\|n_1|,|n_2|\le N}}
\frac{1}{\jb{n_1}^{2-2\al}\jb{n_2}^{2-2\al}
(1+|\kk_1(\bar n)|)}
\end{split}
\label{J1}
\end{align}

\noi
for any $0 \le t \le T$, $n \in \Z^2$,  and $j=0,1$, where
\begin{align*}
\kk_1(\bar n)=\eps_0\jb{n}
-\eps_1\jb{n_1}-\eps_2\jb{n_2}.
%\label{SN5b}
\end{align*}

\noi
Since $\ft Z_N(t,n)=0$ for $|n|>2N$,
we may assume $|n|\le2N$ and, by symmetry,
we assume $|n_1|\ge|n_2|$ in the following.
In the following, we suppress the dependence of implicit constants on $T$.

When  $n=0$, 
it follows from \eqref{CN} that the right-hand side of \eqref{J1} is 
$O(1)$.
Next, we consider the case $0<|n|\le2N$. 
Given dyadic $M \ge 1$, 
decompose $n_2$ into dyadic blocks $\jb{n_2}\sim M$.
Since $n_1=n-n_2$ and $|n_1|\ge|n_2|$, we have
\begin{align}
\jb{n_1}\sim L:=\max(\jb{n},M)\qquad
\text{and}\qquad 
a:=\frac{\jb{n} M}{L}=\min(\jb{n},M)\ge1.
\label{SN5c}
\end{align}

\noi
In the following, we first prove the following bound:
\begin{align}
\sum_{\substack{n = n_1+n_2\\
|n_1|\ge|n_2|}}
\frac{\ind_{\jb{n_2}\sim M}}{1+|\kk_1(\bar n)|}
&\les
M+\frac{\max(\jb{n},M)^{\frac12}M^{\frac32}}
{\jb{n}^{\frac12}}.
\label{SN5a}
\end{align}

\noi
When  $M=1$, 
\eqref{SN5a} holds trivially and thus 
we assume $M \ge 2$ in the following.
We consider the following four cases
for 
$(\eps_0, \eps_1, \eps_2)$: 
\begin{align}
(\pm 1, \mp1, \mp1), \quad 
 (\pm 1, \mp1, \pm1), 
\quad 
 (\pm 1, \pm1, \mp1), \quad
 (\mp 1, \mp1, \mp1).
\label{SN6}
\end{align}

\noi
In the first two cases, we have $|\kk_1(\bar n)|\ges \jb{n}$
since $\jb{n_1} \ge \jb{n_2}$.
Thus, we have 
\begin{align*}
\text{LHS of } \eqref{SN5a}
\les\frac{M^2}{\jb{n}}, 
%\les M^2a^{-\frac12}.
\end{align*}

\noi
satisfying \eqref{SN5a} in these two cases.
We now consider the last two cases in \eqref{SN6}.
Set
\begin{align*}
\ta_-=\angle(n_2,-n)
\qquad \text{and}\qquad
\ta_+=\angle(n_2,n).
\end{align*}

\noi
Then, proceeding as in \cite[(4.23)]{OO} with \eqref{SN5c}, 
we have 
\begin{align}
\jb{n}-\jb{n_1}+\jb{n_2}
&=
\frac{(\jb{n}+\jb{n_2})^2-  \jb{n_1}^2}{\jb{n}+\jb{n_1}+ \jb{n_2}}
\gtrsim\frac{\jb{n}M}{L}(1-\cos\ta_-)
\gtrsim a\ta_-^2.
\label{J1a}
\end{align}

\noi
A similar computation with \cite[(4.22)]{OO} 
and
$\jb{n}\jb{n_2} \ge   1+|n||n_2|$,  which follows from Cauchy-Schwarz's inequality, 
yields
\begin{align}
- \jb{n} + \jb{n_1} + \jb{n_2}
&=\frac{2(\jb{n}\jb{n_2}-n\cdot n_2)-1}{\jb{n_1}+\jb{n}-\jb{n_2}}
\gtrsim\frac{\jb{n}M}{L}(1-\cos\ta_+)
\gtrsim a\ta_+^2.
\label{J1ax}
\end{align}

\noi
Here,  we used $1-\cos\ta\sim\ta^2$ for $0\le\ta\le\pi$.

For the last two cases in \eqref{SN6}, 
we split 
the region of the summation over $n_2$
in \eqref{SN5a}
into the angular regions:
\[\Ta_k = \{ \ta_\pm
\sim 2^ka^{-\frac12}\}, \quad k \in \N \cup\{0\} \text{\ \ such that } 2^k a^{-\frac 12} \les 1\]
with the understanding that
$\Ta_0 = \{ \ta_\pm \les a^{-\frac12}\}$ when $k = 0$.
Then, from \eqref{J1a} and \eqref{J1ax}, 
we have 
\begin{align}
1+ |\kk_1(\bar n)| \ges 2^{2k} .
\label{J1ay}
\end{align}

\noi
Then, from \eqref{J1ay}
and 
Lemma \ref{LEM:count}\,(i) with \eqref{SN5c}, we have 
\begin{align}
\begin{split}
\sum_{\substack{n = n_1+n_2\\
|n_1|\ge|n_2|}}
\frac{\ind_{\jb{n_2}\sim M}}{1+|\kk_1(\bar n)|}
&\les\sum_{k= 0}^\infty 2^{-2k}
\big(M+M^2\min(1,2^ka^{-\frac12})\big)\\
&\les M+M^2a^{-\frac12}\\
&=M+\frac{\max(\jb{n},M)^{\frac12}M^{\frac32}}
{\jb{n}^{\frac12}}.
\end{split}
\label{Zb4}
\end{align}

\noi
Hence, 
 \eqref{SN5a} holds for every choice of signs in \eqref{SN6}.

From \eqref{J1} and \eqref{SN5a} with $\jb{n_1}\sim\max(\jb{n}, M)$
and \eqref{CN}, 
we have 
\begin{align*}
\begin{split}
\E\big[|\ft Z_N(t,n)|^2\big]
&\les C_{\al,N}^4
\sum_{\substack{1\le M \les N\\\mathrm{dyadic}}}
\frac{M^{2\al-1}}{\jb{n}^2\max(\jb{n},M)^{2-2\al}}\\
&\quad+C_{\al,N}^4
\sum_{\substack{1\le M \les N\\\mathrm{dyadic}}}
\frac{M^{2\al-\frac12}}
{\jb{n}^{\frac52}\max(\jb{n},M)^{\frac32-2\al}}\\
&\les \jb{n}^{-3}+\jb{n}^{-\frac52}
\les \jb{n}^{-\frac52}, 
\end{split}
%\label{Zb}
\end{align*}

\noi
where 
 the penultimate step follows from carrying out separate analysis 
 for $M \le \jb{n}$ and $M > \jb{n}$.
This proves \eqref{J0} for $j=0$.
A similar computation yields
\eqref{J0} with $j = 1$.

From 
\eqref{K9} and \eqref{K9a} with \eqref{K2}, we have
\[\E\Big[ |\, \ft{\wick{z_N^2}}\,(t,n)|^2\Big] \les 1.\]

\noi
Then, 
proceeding with an interpolation argument
as in the last part of the proof of \cite[Proposition 1.6]{GKO2}
(see \cite[(3.24)-(3.26)]{GKO2}), we obtain
the following estimates on the time increments:
\begin{align}
\E\big[|\dl_h (\dt^j\ft Z_N)(t,n)|^2\big]
&\les |h|^{\ta}\jb{n}^{2j -\frac52+\ta}
\label{K17}
\end{align}

\noi
for small $\ta > 0$ and any $0 \le t, t+h\le T$ and $n \in \Z^2$, 
uniformly in $N \in \N$.
Now, let 
 $s<s_1<\frac14$ and small $\be > 0$.
Then, it follows from  
\cite[Lemma~2.6]{LLOT} with 
 \eqref{J0} and \eqref{K17}
 (with sufficiently small $\ta > 0$)
that 
\begin{align*}
\sup_{N\in \N}\Big\|
\|Z_N\|_{\C^\be_TW_x^{s_1,\infty}}
+\|\dt Z_N\|_{\C^\be_TW_x^{s_1-1,\infty}}
\Big\|_{L^p(\O)}\le C_p < \infty, 
\end{align*}

\noi
where $\C^\be([0, T]; W^{s_1,\infty}(\T^2))$ is as in \eqref{HO2}.
By the compact embedding
$W^{s_1,\infty}(\T^2) \Subset W^{s,\infty}(\T^2)$
and the Arzel\`a-Ascoli theorem, these bounds give tightness in
$C([0, T]; W^{s,\infty}(\T^2))\cap C^1([0, T]; W^{s-1,\infty}(\T^2))$.
Finally, a diagonal argument over $T\in\N$ proves the claim;
see 
the proofs of 
\cite[Proposition 5.4]{OTh1}
and 
\cite[Proposition~3.5]{LLOT}.
\end{proof}

\section{Convergence of the other stochastic distributions}

\label{SEC:rand2}

In this section, we study regularity and convergence properties
of $z_N$ and $R_N$ defined in~\eqref{zN}
and \eqref{Zrz}, respectively.
Since these results follow
from a straightforward modification of the arguments
in \cite{GKO2, OO, LLOT}, 
we keep our presentation brief in the following.

\begin{lemma}
\label{LEM:randln}
Let $\al \ge \frac12$.
Given any $s < - \frac 12$,
the random linear solution $z_N$ defined in
\eqref{zN} converges to
$0$ in $C(\R_+;W^{s,\infty}(\T^2))
\cap C^1(\R_+;W^{s-1,\infty}(\T^2))$, 
almost surely,
as $N \to \infty$.
\end{lemma}

\begin{proof}
We only discuss convergence of $z_N$, 
since the claimed convergence
of $\dt z_N$ follows from a straightforward modification.
We first consider the case
 $\al>\frac12$.
A direct computation with~\eqref{zN}, \eqref{CN}, 
and the mean value theorem (as in the proof of \cite[Lemma~3.1\,(i)]{GKO2}) 
yields 
\begin{align*}
\E\big[|\ft z_N(t,n)|^2\big]
&\les N^{-\g}\jb{n}^{-1+\g},\\
\E\big[|\dl_h\ft z_N(t,n)|^2\big]
&\les  N^{-\g}\jb{n}^{-1+\g+\ta}|h|^\ta
\end{align*}

\noi
for any 
small  $\g,\ta>0$ with $\g<2\al-1$ 
and 
for any $t, h \in \R$, $n \in \Z^2$, and $N  \in \N$.
Here, $\dl_h$ is as in \eqref{dl1}.
Then, given $s < - \frac 12$, 
by 
taking small $\g,   \ta> 0$ sufficiently small
such that $- 1 + \g  +  \ta < - 2 - 2s$, 
the claim follows from 
\cite[Lemma~2.6]{LLOT};
see also \cite[Proposition~2.7]{OOT}.

Next, we consider the case  $\al=\frac12$.
Let  $y_N=C_{\frac12,N}^{-1}z_N$.
Then, proceeding as above, we have
\begin{align*}
\E\big[|\ft y_N(t,n)|^2\big]
& \les
 \jb{n}^{-1},\\
\E\big[|\ft y_M(t,n) - \ft y_N(t,n)|^2\big]
& 
\les N^{-\g} \jb{n}^{-1+\g}, \\
\E\big[|\dl_h\ft y_N(t,n)|^2\big]
& \les  \jb{n}^{-1+\ta}|h|^{\ta},\\
\E\big[|\dl_h(\ft y_M- \ft y_N)(t,n)|^2\big]
& \les N^{-\g}  \jb{n}^{-1+\g + \ta}|h|^{\ta}
\end{align*}

\noi
for any small $\g, \ta > 0$ and 
for any $t, h \in \R$, $n \in \Z^2$, and $M \ge N \ge 1$.
Then, by noting that $C_{\frac 12, N} \to 0$ as $N \to \infty$, 
the claim in this case also follows from 
\cite[Lemma~2.6]{LLOT}.
\end{proof}

\begin{lemma}\label{LEM:rand2}
Let $\al\ge\frac12$.
Given any $s<-\frac14$, 
the resonant product 
$R_N$ 
defined  in \eqref{Zrz}
converges to $0$
in $C(\R_+;W^{s,\infty}(\T^2))$
in probability as $N\to\infty$.
When $\al>\frac12$, the convergence also holds almost surely.
\end{lemma}

\begin{proof}
We use the decomposition in \cite[(4.1) and (4.3)]{GKO2},
or equivalently \cite[(4.45) and (4.47)]{OO},
with the stochastic convolution and its covariance
replaced by $z_N$ and $\s_{N,n}$ in \eqref{zN} and~\eqref{sigman}, respectively.
Write
\begin{align*}
R_N=\RR_{N,2}+\sum_{j=1}^4\RR_{N,1j}.
%\label{zZ3}
\end{align*}

Fix $T>0$.
Our first goal is to prove the following bound (for sufficiently small $\g > 0$):
\begin{align}
\begin{split}
\E\big[|\ft R_N(t,n)|^2\big]
&\les_T b_N \jb{n}^{-\frac32+\g}
\end{split}
\label{zZ3a}
\end{align}

\noi
for any $0 \le t \le T$ and  $n \in \Z^2$, 
where $b_N$ is given by 
\begin{align}
b_N&=
\begin{cases}
N^{-\g},&\al>\frac12,\\
(\log\jb{N})^{-\frac12},&\al=\frac12.
\end{cases}
\label{zZ3b}
\end{align}

\noi
In the following, 
we suppress the dependence of implicit constants on $T$. 
Proceeding 
with Wick's theorem
as in the proof of \cite[Proposition 1.8]{GKO2}
(see the discussion between \cite[(4.1) and (4.4)]{GKO2};
see also \cite[(4.48)]{OO}
and the discussion preceding it), we have
\begin{align}
\begin{split}
\E\big[|  \ft\RR_{N,2}(t,n) |^2\big]
&\les  \ind_{|n|\les1}\cdot  C_{\al,N}^2,\\
\E\big[| \ft\RR_{N,14}(t,n) |^2\big]
&\les  \ind_{|n|\les N} \cdot  
C_{\al,N}^6    \jb{n}^{6\al-8}, \\
\E\big[|\ft\RR_{N,12}(t,n)|^2\big]
&\les  \ind_{|n|\les N} \cdot  
C_{\al,N}^6     \jb{n}^{2\al-2}  
\sum_{    \jb{n}  \les \jb{m}\les N}\jb{m}^{4\al-6},
\end{split}
\label{zZ6}
\end{align}

\noi
and 
\begin{align}
\begin{split}
\E\big[|\ft\RR_{N,11}(t,n)|^2\big]
&\les 
\ind_{|n|\les N}
\sum_{\substack{n = m+n_3\\|m|\sim |n_3|}} \E\big[|\ft Z_N(t,m)|^2\big]
  \E\big[|\ft z_N(t,n_3)|^2\big]\\
&\les 
\ind_{|n|\les N}\cdot 
C_{\al,N}^2
\sum_{\substack{n = m+n_3\\ |m|, |n_3|\les N\\|m|\sim |n_3|}}
 \frac{1}{\jb{m}^\frac 52 \jb{n_3}^{2-2\al}}\\
& \les
\ind_{|n|\les N}\cdot  C_{\al,N}^2
\sum_{ \jb{n} \lesssim\jb{m}\lesssim N}
\jb{m}^{2\al-\frac92}, 
\end{split}
\label{zZ6a}
\end{align}

\noi
all satisfying  \eqref{zZ3a}.
In \eqref{zZ6a},  we used \eqref{J0},  \eqref{zN}, 
and the fact that $\jb{m}\ges \jb{n}$ under $|m| \sim |n - m|$
(and we also used this  fact in the third bound in \eqref{zZ6}).

It remains to treat the following term (see \cite[(4.3)]{GKO2}):
\begin{align}
\ft\RR_{N,13}(t,n)
=2\int_0^t
\ft z_N(t',n)
\sum_{\substack{|n_2|\le N\\|n_2|\sim |n+n_2|\ne0}}
%\chi(n+m,-m)
\frac{\sin((t-t')\jb{n+n_2})}{\jb{n+n_2}}
\s_{N,n_2}(t,t')dt',
\label{zZ7}
\end{align}

\noi
where
$\s_{N,n_2}$ is as in \eqref{sigman}.
Strictly speaking, 
in \eqref{zZ7}, we should have $\sum_{|j-k|\le2}\varphi_j(n_2)\varphi_k(n+n_2)$
inside the $n_2$-summation (where $\varphi_j$
is the multiplier for the Littlewood-Paley projector $\P_j$), 
 but for simplicity of presentation, 
 we simply use $|n_2|\sim |n+n_2|$ without referring to $\varphi_j$'s.
As in \cite[(4.6)]{GKO2}, we have 
\begin{align}
\begin{split}
&  \E\big[|\ft\RR_{N,13}(t,n)|^2  \big]\\
& \quad \les
\ind_{|n|\les N}\cdot 
\frac{ C_{\al,N}^6}{ \jb{n}^{2- 2\al}}
\sum_{ \eps_1, \eps_2, \eps_3\in\{-1,1\}}
\sum_{\substack{|n_2|\le N\\|n+n_2|\sim |n_2|}}
\frac{1}{\jb{n+n_2} \jb{n_2}^{2-2\al}
(1+|\kk_2(\bar n)|)}\\
& \quad \quad \times 
\sum_{\substack{|n_2'|\le N\\|n+n_2'|\sim |n_2'|}}
\frac{\ind_{|n_2|\ges |n_2'|}}{\jb{n+ n_2'} \jb{n_2'}^{2-2\al}}, 
\end{split}
\label{zZ7a}
\end{align}

\noi
where $\kk_2(\bar n)$ is given by 
\begin{align}
\kk_2(\bar n)=\eps_1\jb{n+n_2}
+\eps_2\jb{n_2}-\eps_3\jb{n}.
\label{zZ7b}
\end{align}

Under $|n+n_2|\sim |n_2|$, we have $|n_2|\ges |n|$.
When $|n_2| \sim |n|$, 
it follows from \eqref{SN5a} that 
\begin{align}
\sum_{\substack{|n_2| \le N \\   |n_2|\sim |n|}}
\frac{ 1}{1+|\kk_2(\bar n)|}
&\les \jb{n}^\frac 32.
\label{zZ7c}
\end{align}

\noi
Using \eqref{zZ7c}, the contribution from
$|n_2|\sim |n|$ in \eqref{zZ7a} is bounded by
\begin{align*}
&
\les \ind_{|n|\les N}\cdot  C_{\al,N}^6
\jb{n}^{-\frac 72 + 4\al}
\sum_{\substack{|n_2'|\le N\\|n+n_2'|\sim|n_2'|}}
\frac{\ind_{|n|\ges |n_2'|}}
{\jb{n+n_2'}\jb{n_2'}^{2-2\al}}\\
&\les
\ind_{|n|\les N}\cdot  C_{\al,N}^6\jb{n}^{6 \al-\frac92}
\les  b_N\jb{n}^{-\frac 32  +\g}, 
\end{align*}

\noi
yielding
 \eqref{zZ3a}.

Next, we consider the case 
$|n_2| \gg |n|$.
For $\al>\frac12$ and $0\le\ta<2\al-1$, we have
\begin{align}
C_{\al,N}^2
\sum_{|m|\le N}
\jb{m}^{2\al-3-\ta}
\les_{\al,\ta} N^{-\ta}.
\label{D6}
\end{align}

\noi
When $\al=\frac12$ and $\ta>0$, we instead have
\begin{align}
C_{\frac12,N}^2
\sum_{|m|\le N}
\jb{m}^{-2-\ta}
\les_\ta C_{\frac12,N}^2.
\label{D6a}
\end{align}

\noi
Then, 
for the non-resonant case, namely,
$|\kk_2(\bar n)| \ges \jb{n_2}$,
we first note that
\begin{align}
\begin{split}
\sum_{\substack{|n_2'| \le N\\
|n+n_2'| \sim |n_2'|}}
\frac{ \ind_{|n_2| \ges |n_2'|} }
{ \jb{n+n_2'} \jb{n_2'}^{2-2\al} }
& \les_\al
\sum_{\substack{\jb{n} \les M \les \jb{n_2}\\
\mathrm{dyadic}}}
M^{2\al-1}\\
& \les_\al
\begin{cases}
\jb{n_2}^{2\al-1},&   \al>\frac12,\\
\log \jb{n_2},
&\al=\frac12.
\end{cases}
\end{split}
\label{D6b}
\end{align}

\noi
Hence, using \eqref{D6b}, \eqref{D6}, and \eqref{D6a},
the corresponding contribution to \eqref{zZ7a}
is bounded by
\begin{align*}
& \les C_{\al,N}^6 \jb{n}^{2\al-2}
\sum_{\substack{|n_2| \le N\\
|n+n_2| \sim |n_2|\\
|n_2| \gg |n|}}
\jb{n_2}^{4\al-5}
\Big(1+\ind_{\al=\frac12}
\log\jb{n_2}\Big)\\
& \les 
\begin{cases}
N^{-\g}\jb{n}^{-2+\g}, & \al > \frac 12,\\
(\log\jb{N})^{-\frac 32} \jb{n}^{-2+\g}, 
& \al = \frac 12, 
\end{cases}
\end{align*}

\noi
yielding  \eqref{zZ3a}.

We now consider the nearly resonant case
$|\kk_2(\bar n)| \ll \jb{n_2}$.
In this case, from \eqref{zZ7b}
and $|n_2| \sim |n+n_2| \gg |n|$, we have
$\eps_1 = - \eps_2$.
This corresponds to Case 3 in the proof of
\cite[Proposition 1.8]{GKO2}.
Fix  small $\g_0> 0$ (to be chosen later), 
and 
we first consider the case
\begin{align*}
\jb{n}\ges\jb{n_2}^{\g_0}.
\end{align*}

\noi
Given small $\dl \ge 0$
 (with $\dl=0$ when $\al>\frac12$
and $\dl>0$ when $\al=\frac12$), we have
\begin{align*}
\jb{n}^{-2\dl}
\les
\jb{n_2}^{-\g_0\dl}
\jb{n_2'}^{-\g_0\dl},
\end{align*}

\noi
since $|n_2'| \les |n_2|$.
Hence,
the  contribution
from this case to \eqref{zZ7a}
is bounded by 
\begin{align}
\begin{split}
& \les
C_{\al,N}^6
\jb{n}^{2\al -2+2\dl}
\sum_{\substack{|n_2| \le N\\
|n+n_2| \sim |n_2|}}
\frac{ 1 }
{ \jb{n+n_2}
\jb{n_2}^{2-2\al+\g_0\dl}
(1+|\kk_2(\bar n)|) }\\
& \qquad\qquad \times
\sum_{\substack{|n_2'| \le N\\
|n+n_2'| \sim |n_2'|}}
\frac{ \ind_{|n_2| \ges |n_2'|} }
{ \jb{n+n_2'}
\jb{n_2'}^{2-2\al+\g_0\dl} }.
\end{split}
\label{LK1}
\end{align}

\noi
As in \eqref{D6b}, the inner sum over $n_2'$ satisfies
\begin{align}
\sum_{\substack{|n_2'| \le N\\
|n+n_2'| \sim |n_2'|}}
\frac{ \ind_{|n_2| \ges |n_2'|} }
{ \jb{n+n_2'}
\jb{n_2'}^{2-2\al+\g_0\dl} }
\les_\al
\begin{cases}
\jb{n_2}^{2\al-1-\g_0\dl},&\al>\frac12,\\
\jb{n}^{-\g_0\dl},&\al=\frac12.
\end{cases}
\label{LK2}
\end{align}

\noi
For each dyadic
$M \sim \jb{n_2} \gg \jb{n}$, 
it follows from 
\eqref{SN5a} that 
\begin{align}
\sum_{\substack{|n_2| \le N\\
\jb{n_2} \sim M\\
|n+n_2| \sim |n_2|}}
\frac{ 1 }
{ 1+|\kk_2(\bar n)| }
\les
M+\frac{ M^2 }{ \jb{n}^{\frac12} }.
\label{LK3}
\end{align}

\noi
Hence, from  \eqref{LK2} and \eqref{LK3}, 
we have
\begin{align*}
\begin{split}
\eqref{LK1}
& \les
C_{\al,N}^6
\jb{n}^{2\al-\frac52+2\dl}
\sum_{\substack{\jb{n} \ll M \le N\\
\mathrm{dyadic}}}
M^{4\al-2-2\g_0\dl}\\
& \les
N^{1-2\al-2\g_0\dl}
\jb{n}^{2\al-\frac52+2\dl}\\
& \les
N^{-\g}
\jb{n}^{-\frac32+\g}
\end{split}
%\label{zZ13c}
\end{align*}

\noi
for $\al > \frac 12$, 
which 
 yields \eqref{zZ3a}
for $0<\g<2\al-1$.
When  $\al = \frac12$, we have
\begin{align*}
\begin{split}
\eqref{LK1}
 \les (\log\jb{N})^{-\frac32}
\jb{n}^{-\frac32+2\dl}, 
\end{split}
%\label{zZ13cx}
\end{align*}

\noi
 yielding \eqref{zZ3a}.

Finally, we consider the case:
$\jb{n}\ll\jb{n_2}^{\g_0}$.
In this case, 
there is no good lower bound on $\kk_2(\bar n)$ in \eqref{zZ7a}.
Instead of using \eqref{zZ7a}, 
 we go back to the expression corresponding to
the right-hand side of \cite[(4.4)]{GKO2}.
Recalling that $\eps_1 = - \eps_2$
in \eqref{zZ7b} (with a slightly different notation than \cite{GKO2}),  
the relevant contribution to
$\E \big[ |\ft \RR_{N,13}(t,n) |^2 \big]$
is bounded by 
\begin{align}
\begin{split}
& \les \frac{ C_{\al,N}^6 }{ \jb{n}^{2-2\al} }
\sum_{\substack{|n_2'| \le N\\
|n+n_2'| \sim |n_2'|}}
\frac{1}
{\jb{n+n_2'} \jb{n_2'}^{2-2\al} }\\
&\quad\times
\bigg|
\sum_{\substack{|n_2| \le N\\
|n+n_2| \sim |n_2| \\\jb{n} \ll \jb{n_2}^{\g_0}}}
\frac{ \ind_{|n_2| \ges |n_2'|} }
{ \jb{n+n_2} \jb{n_2}^{2-2\al} }
\int_0^t
e^{ i t_1 \eps_3 \jb{n} }
\sin\big(
(t-t_1)(\jb{n+n_2}-\jb{n_2})
\big)
\,dt_1
\bigg|.
\end{split}
\label{zZ15}
\end{align}

\noi
We now consider the contribution corresponding to
\cite[(4.10)]{GKO2}:
\begin{align}
\int_0^t
e^{ i t_1\eps_3\jb{n}}
\sum_{\substack{|n_2| \le N\\
|n+n_2| \sim |n_2|\\
\jb{n} \ll \jb{n_2}^{\g_0}}}
\ind_{|n_2| \ges |n_2'|} \cdot 
\frac{
\sin\big((t-t_1)(\jb{n+n_2}-\jb{n_2})\big)
}
{ 
\jb{n+n_2}
\jb{n_2}^{2-2\al} }
dt_1.
\label{zZ16}
\end{align}

\noi
In estimating \eqref{zZ16}, 
we use the sine cancellation
in \cite[(4.11)-(4.14)]{GKO2};
see also \cite[(7.22)-(7.24)]{OOT1}.
As in \cite{GKO2}, 
the condition
$\jb{n}\ll\jb{n_2}^{\g_0}$
allows us to remove
the resonant condition $|n+n_2| \sim |n_2|$ in \eqref{zZ16}, 
which in turn allows us to use the symmetry 
$n_2\leftrightarrow-n_2$
together with 
\begin{align*}
\jb{n\pm n_2}-\jb{n_2}
&=\pm\frac{n\cdot n_2}{\jb{n_2}}
+O\bigg(\frac{\jb{n}^2}{\jb{n_2}}\bigg),
\end{align*}

\noi
which implies
\begin{align}
\bigg|
\frac{\sin(\tau(\jb{n+n_2}-\jb{n_2}))}{\jb{n+n_2}}
+\frac{\sin(\tau(\jb{n-n_2}-\jb{n_2}))}{\jb{n-n_2}}
\bigg|
&\les \frac{\jb{n}^2}{\jb{n_2}^{2}}
\label{zZ18}
\end{align}

\noi
for $0\le\tau\le T$.
Thus, proceeding as in \cite[(4.14)]{GKO2}
with  \eqref{zZ18}, we have
\begin{align}
\begin{split}
|\eqref{zZ16}|
& \les
\jb{n}^2
\sum_{\substack{|n_2| \le N\\
\jb{n} \ll \jb{n_2}^{\g_0}}}
\ind_{|n_2| \ges |n_2'|} \cdot 
\jb{n_2}^{2\al-4}.
\end{split}
\label{zZ18a}
\end{align}

\noi
Hence, using \eqref{zZ18a}
and $\jb{n}\ll\jb{n_2}^{\g_0}$
and choosing $\dl > 0$ sufficiently small, 
we have
\begin{align*}
\begin{split}
\eqref{zZ15}
& \les
\ind_{|n| \les N}\cdot C_{\al,N}^6
\jb{n}^{2\al}
\sum_{\substack{|n_2| \le N\\
\jb{n} \ll \jb{n_2}^{\g_0}}}
\jb{n_2}^{2\al-4}
\sum_{|n_2'| \les |n_2|}
\jb{n_2'}^{2\al-3}\\
& \les
\begin{cases}
N^{-2\dl}
\jb{n}^{1+2\dl-\frac1{\g_0}},
& \al > \frac12,\\
(\log\jb{N})^{-\frac32}
\jb{n}^{1-\frac1{2\g_0}}, 
& \al = \frac12, 
\end{cases}
\end{split}
%\label{zZ19}
\end{align*}

\noi
yielding  \eqref{zZ3a}
by taking $\g_0 > 0$ sufficiently small.

By following the argument at the end of 
\cite[Section~4]{GKO2}, we also obtain the following
increment bound:
\begin{align}
\E[|\dl_h\ft R_N(t,n)|^2]
&\les_T b_N
\jb{n}^{-\frac32+\g+\ta} |h|^{\ta}
\label{zZ1}
\end{align}

\noi
for small $\ta > 0$ and any $0 \le t, t +h \le T$ and $n \in \Z^2$, 
where 
$\dl_h$ and 
$b_N$ are as in \eqref{dl1} and \eqref{zZ3b}, respectively.
Then, for $\al > \frac12$, 
the claimed almost sure convergence
follows from 
\eqref{zZ3a}, \eqref{zZ1}, and
 \cite[Lemma~2.6]{LLOT}.
Similarly, when $\al=\frac12$,
it follows from \eqref{zZ3a}, \eqref{zZ1},
and \cite[Lemma~2.6]{LLOT} with \eqref{zZ3b}
that, for any $s < -  \frac  14$ and  any finite $p\ge1$,
we have
\begin{align*}
\big\|
\|R_N\|_{C_TW_x^{s,\infty}}
\big\|_{L^p(\O)}
\les_{T,p}
(\log\jb{N})^{-\frac14}
\too0.
\end{align*}

\noi
Hence, $R_N$ converges to $0$
in probability in $C([0, T]; W^{s, \infty}(\T^2))$.
\end{proof}

\section{Convergence of the paracontrolled operators}
\label{SEC:para}

In this section, we study the paracontrolled operators
$ \If_{\pl}^{N, (1)}$, $\If^{\Xi_1}_{\pl}$, 
and $\If_{\pl, \pe}^N  $
introduced in Subsection \ref{SUBSEC:IPa}.
Our argument follows \cite[Section~5]{GKO2} and \cite[Section~7]{OOT1},
keeping the frequency sums on $\Z^2$
and the factors $C_{\al, N}$.
We first consider the deterministic operator
$\If_{\pl}^{\Xi_1, (1)}$ in Remark~\ref{REM:para1}.
Given small $\ta > 0$, we have\footnote{As in the proof of Lemma~\ref{LEM:rand2},
we suppress the Littlewood-Paley multipliers
and indicate only the corresponding frequency restrictions.
See \cite[(5.3)]{GKO2} for a rigorous definition
of $\If_{\pl}^{\Xi_1, (1)}$, using the Littlewood-Paley multipliers $\varphi_j$.
A similar comment applies
to $\If_{\pl, \pe}^N$ studied in Lemma \ref{LEM:rand4}
and $\If_{\pl, \pe}^N[f, g]$ studied in Lemma \ref{LEM:A1}.} 
\begin{align}
\begin{split}
\If_{\pl}^{\Xi_1, (1)}(w)(t)
& = \sum_{n\in \Z^2}  e_n
\sum_{\substack{n = n_1+n_2\\
|n_2|^\ta \lesssim |n_1| \ll |n_2|}}
\int_0^t
\frac{\sin((t-t')\jb {n})}{\jb {n}}
\ft w(t', n_1) \ft \Xi_1(t', n_2) dt'.
\end{split}
\label{XX2}
\end{align}

\begin{lemma}\label{LEM:rand3}
Let $s_1, T>0$ and let $\ta>0$ be small.
Given 
$\Xi_1 \in L^3([0, T]; W^{-\frac12-\eps, \infty}(\T^2))$, 
let $\If_{\pl}^{\Xi_1, (1)}$ be as in \eqref{XX2}.
Then, for sufficiently small $\eps_0 = \eps_0 (s_1, \ta)>0$,
independent of $T > 0$, 
we have
\begin{align}
\|\If_{\pl}^{\Xi_1, (1)}(w)\|_{C_T H_x^{\frac12+2\eps_0}}
\les T^{\frac13} 
\|\Xi_1\|_{L_T^3 W_x^{-\frac12-\eps_0, \infty}}
\|w\|_{L_T^3 H_x^{s_1}}
\label{PO1}
\end{align}

\noi
for any $w \in L^3([0, T]; H^{s_1}(\T^2))$.

\end{lemma}

\begin{remark}\label{REM:Z1}
\rm 
By applying 
 Lemma~\ref{LEM:rand3} with $\Xi_1 = z_N$
and  Lemma~\ref{LEM:randln}, we obtain
almost sure convergence of $\If_{\pl}^{N, (1)}$ in \eqref{X2a} to $0$ in
the class:
\begin{align*}
 \L \big(
L^3([0, T]; H^{s_1}(\T^2));
C([0, T]; H^{\frac12+2\eps_0}(\T^2))
\big).
%\label{L1}
\end{align*}

\end{remark}

\begin{proof}[Proof of Lemma \ref{LEM:rand3}]
We follow the proof of \cite[Lemma~5.1]{GKO2};
see also \cite[Lemma~7.1 and (7.3)]{OOT1}.
Under the frequency restrictions in \eqref{XX2}, we have
\begin{align}
\jb {n}^{-\frac12+2\eps_0}
\les \jb {n_1}^{\frac{4\eps_0}{\ta}}
\jb {n_2}^{-\frac12-2\eps_0}
\les \jb {n_1}^{s_1-\eps_0}
\jb {n_2}^{-\frac12-2\eps_0}
\label{PO2}
\end{align}

\noi
for sufficiently small $\eps_0 = \eps_0(s_1, \ta)>0$.
Then, the bound \eqref{PO1} follows from \eqref{PO2} and H\"older's inequality.
\end{proof}

Next, we  consider the paracontrolled operator $\If_{\pl, \pe}^N$
defined in \eqref{X6}.
Its kernel representation is given by
\begin{align}
\If_{\pl, \pe}^N(w)(t)
=\sum_{n\in \Z^2} e_n \int_0^t \sum_{n_1\in \Z^2}
\ft w(t', n_1) \A_{n, n_1}^N(t, t') dt',
\label{POa}
\end{align}

\noi
where
\begin{align}
\A_{n, n_1}^N(t, t')
=\ind_{[0, t]}(t')
\sum_{\substack{n-n_1 = n_2+n_3\\
|n_1| \ll |n_2|^\ta \\|n_1+n_2| \sim |n_3|}}
\frac{\sin((t-t')\jb {n_1+n_2})}{\jb {n_1+n_2}}
\ft z_N(t', n_2) \ft z_N(t, n_3).
\label{POb}
\end{align}

\noi
See \cite[(5.6) and (5.7)]{GKO2}.
Note that $|n_2|, |n_3| \le N$ in \eqref{POb}.

\begin{lemma}\label{LEM:rand4}
Let $\al \ge \frac12$, $0<s_2<1$, and $T>0$.
Then, for sufficiently small $\ta=\ta (s_2)>0$
and any finite $q>1$, the operator $\If_{\pl, \pe}^N$
converges to $0$ almost surely
in $\L(q, T, s_2)$ defined in~\eqref{L2}.
\end{lemma}

\begin{proof}
We follow the proofs of \cite[Proposition~1.11]{GKO2}
and \cite[Proposition~2.6]{OOT1}.
The covariance in the present setting is given by
$\s_{N,n}(t,t')$ in 
\eqref{sigman}.
As in \cite[(7.6)]{OOT1}, write
\begin{align}
\A_{n, n_1}^N
=\A_{n, n_1}^{N, (1)}
+\ind_{n=n_1}\cdot \A_n^{N, (2)},
\label{PO3}
\end{align}

\noi
where $\A^{N, (1)}_{n, n_1}=\A^N_{n, n_1}-\E[\A^N_{n, n_1}]
=\A^N_{n, n_1}-\ind_{n = n_1}\cdot \E[\A^N_{n, n}]$
and $\A_n^{N, (2)}=\E[\A_{n, n}^N]$.

Fix $q>1$. Let $\A(q', T)$ be the Banach space
of kernels with norm:
\begin{align*}
\|A\|_{\A(q', T)}
:=\|\jb {n}^{s_2-1}A_{n, n_1}(t')
\|_{L^{q'}_{t'}([0, T];\ell^2_{n, n_1})}.
%\label{POc}
\end{align*}

\noi
Here, $q'$ denotes the H\"older conjugate of $q$.
Let $\If_{\pl, \pe}^{N, (1)}$ denote the operator
with kernel $\A^{N, (1)}$. Then,
proceeding as in  \cite[(7.9)]{OOT1}, we have
\begin{align}
\|\If_{\pl, \pe}^{N, (1)}\|_{\L(q, T, s_2)}
\les \|\A^{N, (1)}\|_{C([0, T];\A(q', T))}.
\label{PO3a}
\end{align}

\noi
Proceeding as in 
\cite[(7.10)]{OOT1}
with \eqref{PO3}, \eqref{POb}, and \eqref{zN}
and summing over $n_2 \in \Z^2$ with $|n_2|\sim M$, we have 
\begin{align}
\begin{split}
\E \Big[\|\jb {n}^{s_2-1}
\A_{n, n_1}^{N, (1)}(t, t')\|_{\ell^2_{n, n_1}}^2\Big]
&\les C_{\al, N}^4
\sum_{\substack{1\le M\le2N\\\mathrm{dyadic}}}
M^{4\al-4}
\sum_{\substack{|n|\lesssim M\\|n_1|\lesssim M^\ta}}
\jb {n}^{2s_2-2}\\
&\les C_{\al, N}^4
\sum_{\substack{1\le M\le2N\\\mathrm{dyadic}}}
M^{4\al-4+2\ta+2s_2}, 
\end{split}
\label{PO4}
\end{align}

\noi
where we used the fact that $s_2 > 0$ in the second step.
We first consider the case
 $\al>\frac12$.
It follows from \eqref{PO4}, 
 \eqref{CN}, the Wiener chaos estimate
\cite[Lemma~2.5]{LLOT}, and Minkowski's integral inequality
with the fact that $s_2 < 1$ that 
there exists small  
$\g >0$ such that 
\begin{align*}
\big\|\|\A^{N, (1)}(t, \cdot)\|_{\A(q', T)}
\big\|_{L^p(\O)}&\les_{T, p, q}N^{-\g}, 
%\label{PO5a}
\end{align*}

\noi
A similar computation with the mean-value theorem yields
the following bound on 
the time-increment: 
\begin{align*}
\big\|\|\dl_h\A^{N, (1)}(t, \cdot)\|_{\A(q', T)}
\big\|_{L^p(\O)}&\les_{T, p, q}N^{-\g}|h|^\ta
%\label{PO5}
\end{align*}

\noi
for sufficiently small $\g, \ta > 0$, 
where $\dl_h$ is as in \eqref{dl1}.
Then, 
the Kolmogorov continuity criterion 
(\cite[Exercise 8.2]{Bass})
and the Borel-Cantelli lemma (with large $p\gg1$)
yield almost sure convergence of 
$\A^{N, (1)}$  to $0$
in $C([0, T];\A(q', T))$.
Hence, in view of \eqref{PO3a}, 
we conclude that $\If_{\pl, \pe}^{N, (1)}$ converges
almost surely to $0$ in $\L(q, T, s_2)$.

Next, we consider the case 
 $\al=\frac12$.
 Let 
$\B^N=C_{\frac12, N}^{-2}\, \A^{N, (1)}$.
For $M\ge N\ge 1$, 
we have $\max(|n_2|, |n_3|)  > N$
(in \eqref{POb})
for 
the difference $\B^M-\B^N$.
Then, proceeding as above, we have
\begin{align*}
\begin{split}
\big\|\|(\B^M-\B^N)(t, \cdot)\|_{\A(q', T)}
\big\|_{L^p(\O)}&\les_{T, p, q}N^{-\g}, \\
\big\|\|\dl_h(\B^M-\B^N)(t, \cdot)\|_{\A(q', T)}
\big\|_{L^p(\O)}
&\les_{T, p, q}N^{-\g}|h|^\ta.
\end{split}
%\label{PO5c}
\end{align*}

\noi
Thus,  from 
the Kolmogorov continuity criterion 
and the Borel-Cantelli lemma (with large $p\gg1$),
we see  that 
$\B^N$ converges almost surely to some limit
in $C([0, T];\A(q', T))$.
Hence, from \eqref{CN1} and 
 \eqref{PO3a}, 
we conclude almost sure convergence of $\If_{\pl, \pe}^{N, (1)}$
to $0$ in 
 $\L(q, T, s_2)$.

For the deterministic part, the frequency restriction is
$|n| \ll |n_2|^\ta$, which is invariant under $n_2 \mapsto -n_2$.
The suppressed multiplier also has this symmetry.
As in \cite[(7.7)]{OOT1}, decompose the diagonal kernel
$\A_n^{N, (2)}$ as follows:
\begin{align*}
\A_n^{N, (2)}(t, t')&=
\A_n^{N, +}(t, t')
+ \A_n^{N, -}(t, t'), 
\end{align*}

\noi
where 
\begin{align*}
\A_n^{N, \pm}(t, t')
&=
\ind_{[0, t]}(t')\cdot 
C_{\al, N}^2\sum_{\substack{|n_2|\le N\\|n|\ll|n_2|^\ta}}
\frac{\sin((t-t')(\jb {n+n_2}\pm\jb {n_2}))}
{2\jb {n+n_2}\jb {n_2}^{2-2\al}}.
%\label{PO6}
\end{align*}

\noi
Let $\If_{\pl, \pe}^{N, \pm}$ denote the operators
with kernels $\A^{N, \pm}$.
Interpolating  \eqref{zZ18}
and the trivial bound by $\jb{n_2}^{-1}$, we have
\begin{align}
\bigg|
\frac{\sin(\tau(\jb{n+n_2}-\jb{n_2}))}{\jb{n+n_2}}
+\frac{\sin(\tau(\jb{n-n_2}-\jb{n_2}))}{\jb{n-n_2}}
\bigg|
&\les \frac{\jb{n}^{2\g}}{\jb{n_2}^{1+\g}}
\label{PO6a}
\end{align}

\noi
for any $0 < \g < 1$.
Then, by 
 the symmetrization
$n_2 \leftrightarrow-n_2$
and \eqref{PO6a}
(see \cite[(7.22)--(7.24)]{OOT1}), 
we have
\begin{align}
|\A_n^{N, -}(t, t')|
\les_T C_{\al, N}^2\jb {n}^{2\g}
\sum_{|n_2|\le N}\jb {n_2}^{2\al-3-\g}.
\label{PO8}
\end{align}

\noi
By taking $\g>0$ sufficiently small
and applying  \eqref{D6} and \eqref{D6a} to \eqref{PO8}, we obtain
\begin{align*}
\sup_{t, t'\in[0, T]}\sup_n
\jb {n}^{s_2-1}|\A_n^{N, -}(t, t')|
\les_T
\begin{cases}
N^{-\g}, &\al>\frac12, \\
C_{\frac12, N}^2, &\al=\frac12, 
\end{cases}
%\label{PO9}
\end{align*}

\noi
which tends to $0$ as $N \to \infty$.
Hence, from \eqref{POa} with $n = n_1$, 
we conclude 
that 
$\If_{\pl, \pe}^{N, -}$ converges to $0$
in $\L(q, T, s_2)$.

We now consider the non-resonant part
$\A_n^{N, +}$.
Let 
$\Kc_n^N(t)=\ind_{[0, T]}(t)\cdot \A_n^{N, +}(t, 0)$.
We use the Fourier-transform argument in
\cite[(7.15)-(7.20)]{OOT1}.
For $\phi_{n, \tau}(n_2)=|\jb {n+n_2} +\jb {n_2}-\tau|$,
the corresponding two-dimensional counts are
\begin{align*}
\#\{n_2 \in \Z^2 :\phi_{n, \tau}(n_2)\sim M, \ |n|\ll|n_2|^\ta \}
\les
\begin{cases}
M\tau, &M\ll \tau, \\
\tau^2, &M\sim \tau, \\
M^2, &M\gg \tau,
\end{cases}
%\label{PO15}
\end{align*}

\noi
for $M\ge1$ and $\tau\ge1$.
The region $\phi_{n, \tau}\le1$ contains $O(\tau)$ points.
Then, by carrying out  the dyadic summation 
as in \cite[(7.18) and (7.19)]{OOT1}
with these counts, we obtain
\begin{align}
|\ft \Kc_n^N(\tau)|\les_T
\begin{cases}
N^{-\g}\displaystyle
\frac{\log(2+|\tau|)}{\jb {\tau}^{1-\g}},
&\al>\frac12,  \rule[-5mm]{0pt}{0pt}\\
C_{\frac12, N}^2\displaystyle
\frac{\log(2+|\tau|)}{\jb {\tau}},
&\al=\frac12,
\end{cases}
\label{PO16}
\end{align}

\noi
uniformly in $n\in \Z^2$, with sufficiently small $\g>0$
in the first case.
Here, when $\al > \frac 12$, we used
$C_{\al, N}^2\jb {n_2}^{2\al-2}
\les N^{-\g}\jb {n_2}^{-1+\g}$ for $|n_2|\le N$.

Let  $1<q\le2$.
By applying 
Hausdorff-Young's inequality
and Minkowski's inequality, as in
\cite[(7.20) and the math display immediately after (7.20)]{OOT1}
and \eqref{PO16} 
 with sufficiently small  $\g>0$, 
 we have
\begin{align*}
\|\If_{\pl, \pe}^{N, +}(w)\|_{C_TL_x^2}
\les \|w\|_{L_T^qL_x^2}
\sup_{n\in \Z^2} \|\ft \Kc_n^N\|_{L_\tau^q}
=o(1)\|w\|_{L_T^qL_x^2}
\end{align*}

\noi
as $N \to \infty$, 
showing convergence of 
$\If_{\pl, \pe}^{N, +}$ to  $0$
in $\L(q, T, s_2)$.
For $q>2$, we use the embedding
$L^q([0, T]; L^2(\T^2)) \subset L^2([0, T]; L^2(\T^2))$.

Lastly, 
for each $N \in \N$, the kernel is continuous on $0\le t'\le t$
and has finite spatial Fourier support.
Thus, the Lebesgue dominated convergence theorem 
(applied to $t_k \to t$) yields continuity
of $\If_{\pl, \pe}^{N}(w)$ in time
for  $w \in L^q([0, T]; L^2(\T^2))$.

This concludes the proof of Lemma \ref{LEM:rand4}.
\end{proof}

\section{Local well-posedness
of the paracontrolled system}
\label{SEC:conv1}

In this section, we 
establish
local well-posedness 
of the deterministic 
paracontrolled 
system~\eqref{NLW5x}
(Proposition~\ref{PROP:LWP1}).

\begin{proof}[Proof of Proposition~\ref{PROP:LWP1}]
Fix  $\frac14<s_1<\frac12<s_2<\frac34$
with $s_2\le2s_1$, 
and let  $\ta = \ta (s_2)>0$ and
$\eps = \eps (s_1, s_2)>0$ be sufficiently small.
Given $ \Xi = (\Xi_1, \Xi_2, \Xi_3, \Xi_4)
 \in \Xc^{\eps, s_2}$, 
define
$\Phi^\Xi = (\Phi_1^\Xi, \Phi_2^\Xi)$ 
by setting
\begin{align*}
\begin{split}
\Phi_1^\Xi (\Xs, \Ys) 
&= \I\big (\text{RHS of the $\Xs$-equation in \eqref{NLW5x}}\big), \\
\Phi_2^\Xi (\Xs, \Ys) 
&= \I\big (\text{RHS of the $\Ys$-equation in \eqref{NLW5x}}\big).\\
\end{split}
%\label{A7}
\end{align*}

Let $0 < T \le 1$.
Suppose that 
$\|(\Xs, \Ys)\|_{\vec X^{s_1, s_2}(T)}\le R$
for some $R \ge 1$.
From 
 the energy estimate:
\begin{align}
\|\I (F)\|_{X^s(T)}
\les \|F\|_{L_T^1 H_x^{s-1}}, 
\label{A4}
\end{align}

\noi
where $X^s(T)$ is as in \eqref{Xsp}, 
and Lemma~\ref{LEM:para} with \eqref{besov1}, we have
\begin{align}
\|\Phi_1^\Xi (\Xs, \Ys) \|_{X^{s_1}(T)}
\les T \|\Xi\|_{\Xc^{\eps, s_2}}
\Big(\|\Xi\|_{\Xc^{\eps, s_2}} +R\Big).
\label{AA1}
\end{align}

\noi
For simplicity of notation, set $w = \Xs + \Ys  + \Xi_2$.
Then, 
from Lemma~\ref{LEM:para} with \eqref{besov1}, we have
\begin{align}
\begin{split}
&   \|w \pg \Xi_1\|_{C_T H_x^{s_2-1}}
 + 
\|\Ys \pe \Xi_1\|_{C_T H_x^{s_2-1}}
 + \|\Xi_3\|_{C_T H_x^{s_2-1}}\\
& \quad \les  \|\Xi\|_{\Xc^{\eps, s_2}}
\Big(\|\Xi\|_{\Xc^{\eps, s_2}} +R\Big).
\end{split}
\label{A8}
\end{align}

\noi
By 
Lemmas~\ref{LEM:para}
and~\ref{LEM:rand3} with \eqref{besov1}, we have 
\begin{align}
\begin{split}
\|\If_{\pl}^{\Xi_1, (1)}(w) \pe \Xi_1 \|_{C_T H_x^{s_2 - 1}}
&\les %T^{\frac23} 
\|\Xi_1\|_{C_T W_x^{-\frac12-\eps, \infty}}^2
\|w\|_{L_T^3 H_x^{\frac 14 - \eps}}\\
& \le 
\|\Xi\|_{\Xc^{\eps, s_2}}^2
\Big(\|\Xi\|_{\Xc^{\eps, s_2}} +R\Big).
\end{split}
\label{A9}
\end{align}

\noi
From \eqref{L2}, we have 
\begin{align}
\begin{split}
\|\Xi_4(w)\|_{C_T H_x^{s_2 - 1}}
&\les  \|\Xi_4\|_{\L(q, T, s_2)}
\|w\|_{L_T^q L_x^2}
\le %T^\frac 23 
\|\Xi\|_{\Xc^{\eps, s_2}}
\Big(\|\Xi\|_{\Xc^{\eps, s_2}} +R\Big).
\end{split}
\label{A11}
\end{align}

\noi
From Sobolev's inequality  and H\"older's inequality, we have
\begin{align}
\|fg\|_{H^{s_2-1}}
\les \|fg\|_{L^{\frac2{2-s_2}}}
\le \|f\|_{L^{\frac4{2-s_2}}}
\|g\|_{L^{\frac4{2-s_2}}}
\les \|f\|_{H^{\frac{s_2}2}}
\|g\|_{H^{\frac{s_2}2}}.
\label{A5}
\end{align}

\noi
Thus, for $s_2\le2s_1$, 
we have
\begin{align}
\begin{split}
\|w^2\|_{C_T H_x^{s_2 - 1}}
&\les \Big(
\|\Xs+\Ys\|_{C_T H_x^{\frac{s_2}2}}^2
+ \|\Xi_2\|_{C_T L_x^\infty}
\|\Xs+\Ys\|_{C_T L_x^2}
+ \|\Xi_2\|_{C_T L_x^\infty}^2\Big)\\
&\le \Big(\|\Xi\|_{\Xc^{\eps, s_2}} +R\Big)^2.
\end{split}
\label{A12}
\end{align}

\noi
Putting \eqref{A8}, \eqref{A9}, \eqref{A11}, and \eqref{A12}
together with \eqref{A4}, we obtain
\begin{align}
\|\Phi_2^\Xi (\Xs, \Ys) \|_{X^{s_2}(T)}
\les T \Big(1 + \|\Xi\|_{\Xc^{\eps, s_2}}\Big)
\Big(\|\Xi\|_{\Xc^{\eps, s_2}} +R\Big)^2,
\label{AA2}
\end{align}

Hence, from \eqref{AA1} and \eqref{AA2}
and taking $T = T\big( \|\Xi\|_{\Xc^{\eps, s_2}}, R\big) > 0$
sufficiently small, we have
\begin{align*}
\|\Phi^\Xi (\Xs, \Ys)\|_{\vec X^{s_1, s_2}(T)}
\le R.
\end{align*}

\noi
A similar computation yields
the following difference estimate:
\begin{align*}
\|\Phi^\Xi (\Xs_1, \Ys_1) - \Phi^\Xi (\Xs_2, \Ys_2)\|_{\vec X^{s_1, s_2}(T)}
\le \frac 12 
\| (\Xs_1, \Ys_1) -  (\Xs_2, \Ys_2)\|_{\vec X^{s_1, s_2}(T)}.
%\label{AA3}
\end{align*}

\noi
Therefore, local well-posedness of the paracontrolled  system \eqref{NLW5x}
follows from the Banach fixed point theorem
(with $R = 1$).
A standard modification of the argument above yields
local Lipschitz continuity of the map:
$\Xi \in \Xc^{\eps, s_2}
\mapsto 
(\Xs, \Ys) \in   \vec X^{s_1, s_2}(T )$.
Note that the contraction argument 
(with $R = 1$)
yields uniqueness
of the solution in the ball of radius $1$ in 
$\vec X^{s_1, s_2}(T)$, 
but by  a standard continuity argument, 
we can extend the uniqueness to the entire class
$ \vec X^{s_1, s_2}(T )$.
\end{proof}

\section{Singularity formation for the limiting equation}
%{Finite time blowup for the limiting equation}
\label{SEC:B}

In this section, 
we present the proof of Proposition \ref{PROP:B1}
on local well-posedness
and almost sure finite time blowup
for the limiting SNLW \eqref{NLW1b}.

\begin{proof} %[Proof of Proposition \ref{PROP:B1}]
(i)
From 
 Lemma~\ref{LEM:randZ}, the linear solution $Z$
defined in \eqref{Z1} belongs to
$C(\R_+; H^s(\T^2))
\cap C^1(\R_+; H^{s-1}(\T^2))$, almost surely.
Fix a realization of $Z$ with this regularity 
and consider  the first order expansion 
\begin{align}
u =  Z+v, 
\label{F2a}
\end{align}

\noi
where the remainder term $v$ satisfies~\eqref{NLW3b}:
\begin{align}
\begin{cases}
\dt^2 v + (1 - \Dl) v = (v + Z)^2\\
(v, \dt v)|_{t = 0} = (0, 0).
\end{cases}
\label{F3}
\end{align}

Let $t_0 \in \R_+$
and consider 
the time  interval $I = [t_0, t_0 + \tau]$ with $0 < \tau \le 1$.
The Duhamel formulation of \eqref{F3}, 
starting at time $t_0$ reads as
\begin{align}
v(t) = W(t- t_0) (v(t_0), \dt v(t_0) )
+ \int_{t_0}^t \frac{\sin((t-t')\jb{\nb})}{\jb{\nb}} \big((v+Z)^2(t')\big)dt'.
\label{F4}
\end{align}

\noi
Denote the right-hand side of \eqref{F4} by $\G(v)(t)$.
By applying \eqref{A4} and  \eqref{A5}
(with $s_2 = s$), we have 
\begin{align}
\begin{split}
\|  (\G(v), \dt \G(v))  \|_{C_I \vec H_x^s}
\les \|  (v(t_0), \dt v(t_0)) \|_{\vec H^s}
+ \tau \Big(   \|v\|_{C_I H_x^s} + \|Z\|_{C_I H_x^s}  \Big)^2.
\end{split}
\label{F5}
\end{align}

\noi
A similar computation yields
the difference estimate:
\begin{align}
\begin{split}
& \|(\G(v), \dt \G(v)) - (\G(w), \dt \G(w))\|_{C_I \vec H_x^s}\\
& \quad \les 
 \tau \Big(\|v\|_{C_I H_x^s}
 + \|w\|_{C_I H_x^s}
+ \|Z\|_{C_I H_x^s}\Big)
\|v - w\|_{C_I H_x^s}.
\end{split}
\label{F5a}
\end{align}

\noi
Then, it follows from a contraction argument with \eqref{F5} and \eqref{F5a}
that \eqref{F3} is locally well-posed
on $[t_0, t_0 + \tau]$, 
provided that $0 < \tau \le 1$ is sufficiently small, 
depending only on 
$\|(v(t_0), \dt v(t_0))\|_{\vec H^s}$
and $\|Z\|_{C([0, t_0+1]; H^s_x)}$.

\vspace{1cm}

\medskip

\noi
(ii)
We first verify  \eqref{F2} under \eqref{F1}.
The local well-posedness argument in Part (i)
 yields the  blowup alternative \eqref{F5b}:
\begin{align}
 T_{\max} = \infty \qquad \text{or}\qquad
\lim_{t \to T_{\max}-}
\|(v(t), \dt v(t))\|_{\vec H^s} = \infty, 
\label{F5bx}
\end{align}

\noi
almost surely.
Hence, \eqref{F2} follows from \eqref{F5bx}, \eqref{F2a}, 
and Lemma \ref{LEM:randZ}.

\medskip

The remaining part of the proof is devoted to proving  \eqref{F1}.
We proceed
by  a contradiction argument, using  the test function method
of Zhang \cite{Zha1, Zha2};
see also \cite{Oh4, OOP}.
Namely, we suppose that 
\begin{align}
 \PP (  T_{\max} = \infty  ) >0.
\label{F1a}
\end{align}

\noi
Then, by writing \eqref{NLW1b}
in the weak formulation \eqref{F10}
and choosing a suitable family of test functions, 
we derive a contradiction to \eqref{F1a}
by making use of the non-negativity of the nonlinearity $u^2$
and the stationarity (in time) of the noise $\ze$.
See \eqref{F13b} below.

In the following, we restrict our attention to  the set of full probability, 
where \eqref{NLW1b} is locally well-posed.
Suppose that \eqref{F1} does not hold.
Namely, suppose that~\eqref{F1a} holds.
In the following, we work on the set $\Si: =\{ T_{\max} = \infty\}$.
It is easy to see that the Duhamel solution 
$u \in C(\R_+; H^s(\T^2))$
satisfies the following weak formulation:
\begin{align}
\jb{ u,  \L\eta }
= \jb{ u^2 ,  \eta}
+ \jb{\ze,  \eta}
\label{F10}
\end{align}

\noi
for any $\eta \in C^\infty_c(\R_+\times \T^2)$, 
where $\L = \dt^2 + 1 - \Dl$
and
\[ \jb{f, g}  = \int_0^\infty  \int_{\T^2} f g\, dx dt.\]

\noi
See \cite[Definition 1.1]{Oh4}
and \cite[Definition 8.1]{OOP}.
In the following, we introduce a suitable family of  non-negative test functions,
and derive a  contradiction to the assumption \eqref{F1a}.

Fix a non-zero, non-negative cutoff function
$\chi \in C_c^\infty((0, 1))$.
Then, set
\begin{align}
\begin{split}
& \psi = \chi^4,
\qquad \phi(x) = 3 + e^{i x_1} + e^{-i x_1}, \qquad \text{and}\\
%\qquad 
& \eta_k(t, x) = \psi(t-k) \phi(x),
\quad k \in \N\cup\{0\}, 
\end{split}
\label{F11}
\end{align}

\noi
where $x = (x_1, x_2) \in \T^2$.
Then, from 
 \eqref{F10} with $\eta = \eta_k$, we have 
\begin{align}
\begin{split}
\jb{\ze, \eta_k} %:\! & = \int_0^T  \int_{\T^2} \ze \eta_k\,  dx dt
& = \iint_{\{\eta_k > 0\}}
\Big(u \L\eta_k - u^2 \eta_k\Big) dx dt\\
&\le \frac{1}{4} 
\iint_{\{\eta_k > 0\}}
\frac{  |\L\eta_k |^2}{\eta_k} dx dt\\
&= \frac{1}{4} 
\int_{\{\chi > 0\}} \int_{\T^2}
\frac{|\psi''(t) \phi(x)
+ \psi(t) (1 - \Dl) \phi(x)|^2}
{\psi(t) \phi(x)} dx dt
=: K, 
\end{split}
\label{F12}
\end{align}

\noi
where, in the second step, we used
\[
u  \L\eta_k-u^2\eta_k
=
-\eta_k   \bigg(u-\frac{\L\eta_k}{2\eta_k} \bigg)^2
+\frac{|\L\eta_k|^2}{4\eta_k}
\le \frac{|\L\eta_k|^2}{4\eta_k}.
\]

\noi
We note that $K$ in \eqref{F12} is independent of $u$ and $k \in \N\cup\{0\}$.
From \eqref{F11}, we have $\phi \sim 1$
and $|(1-\Dl) \phi| \les 1$.
Then, noting that 
$\psi'' = 12 \chi^2 (\chi')^2 + 4 \chi^3 \chi''$, we have
\begin{align}
\begin{split}
K
&\les \int_{\{\chi > 0\}}
\bigg(\frac{|\psi''(t)|^2}{\psi(t)} + \psi(t)\bigg) dt\\
&\les \int_0^1
\Big(|\chi'(t)|^4
+ \chi(t)^2 |\chi''(t)|^2 + \chi(t)^4\Big) dt
< \infty.
\end{split}
\label{F13}
\end{align}

\noi
Hence, from \eqref{F12} and \eqref{F13}, we obtain
\begin{align}
\jb{\ze, \eta_k} \les 1
\label{F13a}
\end{align}

\noi
on $\Si = \{T_{\max} = \infty\}$, 
uniformly in $k \in \N \cup\{0\}$.

In the following, we will prove that 
\begin{align}
\sup_{k \in \N \cup\{0\}} \jb{\ze, \eta_k}  = \infty,
\quad \text{almost surely}, 
\label{F13b}
\end{align}

\noi
which is a contradiction to \eqref{F13a}
and thus to \eqref{F1a}, yielding 
 \eqref{F1}.

Hence, it remains to prove \eqref{F13b}.
Recall from the discussion at the beginning 
of Section~\ref{SEC:Psecond}
that 
$\ze_0(t)$ is constant in time, 
given by
a real-valued standard Gaussian, 
call it $\ze_0$.
It follows from 
 \eqref{F11} that 
\begin{align}
\langle \ze, \eta_k \rangle
&= 3\ze_0 \int_0^1 \psi(t) dt +2 X_k,
\label{F15}
\end{align}

\noi
where $X_k$ is given by 
\begin{align}
X_k &= \int_0^1 \psi(t)
\Re \ze_{(1, 0)}(k+t) dt.
\label{F15a}
\end{align}

\noi
From  \eqref{nxib},  $\ze_{-n} = \cj{\ze_n}$, 
and evenness of $J_0$, we see that 
$\{X_k\}_{k \in \N \cup\{0\}}$ is 
a centered stationary real-valued Gaussian sequence
 with  covariance
\begin{align}
\E[X_0 X_k]
= \frac{1}{2} \int_0^1 \int_0^1
\psi(t) \psi(t') J_0(k+t'-t) dt dt'.
\label{F16}
\end{align}

\noi
Here, the stationarity means that 
$\Law \big((X_j, X_{j+1}, \dots)\big)
= \Law \big((X_0, X_{1}, \dots)\big)$ for any $j \in \N$
and in particular, we have 
\begin{align*}
\E[X_j X_{k+j}] = \E[X_0 X_k]  
%\label{F16a}
\end{align*}

\noi
for any $k, j \in \N \cup\{0\}$.
Then, from \eqref{F16} with the fact that 
$J_0(r) \ge \cos 1 > 0$ for $|r| \le 1$
(which follows from  \eqref{Bessel1}), 
we have 
\begin{align}
\E[X_0^2]
\ge \frac{\cos 1}{2}
\bigg(\int_0^1 \psi(t) dt\bigg)^2 > 0.
\label{F17}
\end{align}

\noi
On the other hand, by \eqref{Bessel2} and \eqref{F16},
we have
\begin{align}
\E[X_0 X_k] \too 0, 
\label{F18}
\end{align}

\noi
as $ k \to \infty$.
We now claim that
\begin{align}
\sup_{k \in \N \cup\{0\}} X_k = \infty, 
\quad \text{almost surely}.
\label{F19}
\end{align}

\noi
Once we prove \eqref{F19}, 
the desired identity \eqref{F13b} follows from 
\eqref{F15} and \eqref{F19} 
with the almost sure finiteness of 
the Gaussian random variable $\ze_0$.

Given $R > 0$, let 
$Y_k = \ind_{\{X_k > R\}}$ and
$p_R = \E[Y_k] = \PP(X_0 > R)$.
Since $X_0$ is a mean-zero Gaussian random variable
with positive variance (see \eqref{F17}), we have $p_R > 0$.
Noting from~\eqref{F15a}
that  $(X_0, X_k)$ is jointly Gaussian, 
it follows from \eqref{F18} 
that, as $k \to \infty$,  $(X_0, X_k)$ converges in law to a pair of
independent centered Gaussians with variance $\E[X_0^2]$
(which can be seen from the convergence
of the characteristic functions; see also \cite[Lemma~2.7]{LLOT}, where the condition (ii) is vacuously satisfied.)
In particular, we have 
\begin{align}
\E[Y_0 Y_k]
= \PP(X_0 > R, \, X_k > R)
\longrightarrow p_R^2, 
\label{F20}
\end{align}

\noi
as $k \to \infty$.
By expanding the square with $\E[Y_k]  = \E[Y_k^2] = p_R$
and applying the stationarity of the family $\{Y_k\}_{k \in \N \cup\{0\}}$ 
(which follows from the stationarity of $\{X_k\}_{k \in \N \cup\{0\}}$) 
and \eqref{F20}, we obtain
\begin{align}
\begin{split}
\E \Bigg[
\bigg|\frac{1}{M} \sum_{k=1}^M Y_k - p_R\bigg|^2
\Bigg]
&\le \frac{p_R(1-p_R)}{M}
+ \frac{2}{M^2} \sum_{1\le k < \l \le M}
 \big|\E[Y_0 Y_{\l - k} ] - p_R^2\big|\\
& = \frac{p_R(1-p_R)}{M}
+ \frac{2}{M^2} \sum_{j=1}^{M-1}
(M-j) \big|\E[Y_0 Y_j] - p_R^2\big|\\
&\too 0, 
\end{split}
\label{F21}
\end{align}

\noi
as $M \to \infty$.
On the event $A_R : = \{\sup_{k \in \N\cup\{0\}} X_k \le R\}$,
we have $Y_k = 0$ for any  $k \in \N\cup\{0\}$.
Hence, from~\eqref{F21}, we have
\begin{align*}
\PP(A_R)
& =  \frac{1}{p_R^2}
\E \Bigg[\ind_{A_R}
\bigg|\frac{1}{M} \sum_{k=1}^M Y_k - p_R\bigg|^2
\Bigg]
\le \frac{1}{p_R^2}
\E \Bigg[
\bigg|\frac{1}{M} \sum_{k=1}^M Y_k - p_R\bigg|^2
\Bigg] \\
& \too 0, 
\end{align*}

\noi
as $M \to \infty$.
Namely, we have 
\begin{align}
\PP(A_R) = 0
\label{F22}
\end{align}

\noi
 for any fixed $R > 0$.
Finally, by writing
\begin{align}
\bigg\{\sup_{k \in \N \cup\{0\}} X_k < \infty \bigg\}
= \bigcup_{R\in \N} A_R
\label{F23}, 
\end{align}

\noi
the claim \eqref{F19} follows from 
\eqref{F23} and \eqref{F22}.
\end{proof}

\begin{remark}\label{REM:uniq2}
\rm

In this remark, we prove \eqref{FF2a}
claimed in 
Remark \ref{REM:uniq1}.
%
%In  the proof of 
%Proposition \ref{PROP:B1}\,(i), we proved local well-posedness
%of~\eqref{F3} (namely, \eqref{NLW3b})
%for the remainder term $v$ in \eqref{F2a} (namely, \eqref{dec2}), 
%%of the limiting equation \eqref{NLW1b}.
%where the local-in-time solution $v$
%belongs to $X^s(\tau)$, 
%$0 < s < \frac 14$, for an almost surely positive 
%local existence time $\tau$.
%Here,  $X^s(\tau)$ is as in \eqref{Xsp}.
%By using the blowup alternative \eqref{F5b}
%and iterating the local well-posedness argument, 
%we can extend $v$ to the maximal time interval $[0, T_{\max})$, 
%where $T_{\max} < \infty$ almost surely 
%and hence
%%\begin{align}
%%\limsup_{t \to T_{\max}-}
%%\|(v(t), \dt v(t))\|_{\vec H^s} = \infty.
%%\label{FF1}
%%\end{align}
%
%
%\vspace{1cm}
%

Fix $0 < T < T_{\max}$.
In the following, all the estimates hold almost surely
and we drop ``almost surely'' for simplicity of presentation.
From \eqref{A4},  \eqref{A5}
(with $s_2$ replaced by $2s$), and Proposition \ref{PROP:B1}, 
we have
\begin{align}
\|v\|_{X^{2s}(T)}
= \|\I(u^2)\|_{X^{2s}(T)}
\les T   \|u\|_{C_T H_x^s}^2
< \infty.
\label{FF3}
\end{align}

\noi
By applying  \eqref{A4} and \eqref{A5}
(with $s_2$ replaced by $4s < 1$) with  
\eqref{FF3} and the fact that $Z \in C([0, T]; L^\infty(\T^2))$, 
we have 
\begin{align*}
\|v\|_{X^{4s}(T)}
&\les T   \Big(
\|v\|_{C_T H_x^{2s}}^2
+ \|Z\|_{C_T L_x^\infty}
  \|v\|_{C_T L_x^2}
+ \|Z\|_{C_T L_x^\infty}^2
\Big)
< \infty.
\end{align*}

\noi
By iterating this argument, we have 
\begin{align}
\|v\|_{X^{2^k s}(T)}
&\les T   \Big(
\|v\|_{C_T H_x^{2^{k-1}s}}^2
+ \|Z\|_{C_T L_x^\infty}
  \|v\|_{C_T L_x^2}
+ \|Z\|_{C_T L_x^\infty}^2
\Big)
< \infty, 
\label{FF4}
\end{align}

\noi
where $\frac{s_2}{2} \le 2^k s < s_2 $.
Note that we used the fact that 
$ 2^k s < 1$ in obtaining \eqref{FF4}.
Then, by applying  \eqref{A4} and \eqref{A5} once again, 
we obtain
\begin{align*}
\|v\|_{X^{s_2}(T)}
&\les T   \Big(
\|v\|_{C_T H_x^{\frac{s_2}{2}}}^2
+ \|Z\|_{C_T L_x^\infty}
  \|v\|_{C_T L_x^2}
+ \|Z\|_{C_T L_x^\infty}^2
\Big)
< \infty.
\end{align*}

\noi
Since the choice of 
 $T < T_{\max}$ was arbitrary,
 we conclude \eqref{FF2a}.
Lastly, we note that the argument above shows that 
the solution $v$ to the limiting equation \eqref{F3} (namely, \eqref{NLW3b})
is unique in $X^s(T)$ for 
any  $0 < T < T_{\max}$
and
any small $s>0$.

\end{remark}

\section{On the application of
the Skorokhod representation theorem}
\label{SEC:Sk}

In this section, we go over the application of  the Skorokhod representation 
theorem (Lemma~\ref{LEM:Sk})
discussed in Subsection~\ref{SUBSEC:full}, 
yielding 
almost sure convergence of the extended enhanced data set 
in \eqref{data2z} on a new probability space.

\subsection{Measurability of the extended enhanced data}
\label{SUBSEC:Sk1}

In this subsection, we introduce the metric space $(\M, d_{\M})$
and
show that the extended enhanced data \eqref{data2x} is a Borel measurable $\M$-valued random variable.
In Subsection~\ref{SUBSEC:Sk0}, 
we will prove the weak convergence of the extended enhanced data 
set in \eqref{data2x} 
on this metric space $(\M, d_{\M})$.

Fix $1<q<\infty$ and $s, s_2\in \R$. 
Let $\B$ be the class of 
 linear operators
from $L^q_{\rm loc}(\R_+;L^2(\T^2))$ into  $C(\R_+;H^{s_2-1}(\T^2))$.
We say that
 $L$ is  adapted
if it satisfies \eqref{adapt1} for any $t \in \R_+$.
Set
\begin{align}
\begin{aligned}
\L_{a}(q,\infty, s_2) =\big\{
& L \in \B: 
\text{$L$ is adapted and}\\
%&\hphantom{=\big\{A:{}}
& L|_{L^q_TL^2_x}\in\L(q,T,s_2)\text{ for every }T>0\big\}, 
\end{aligned}
\label{JW1}
\end{align}

\noi
where $\L(q,T,s_2)$ 
is  as in \eqref{L2}.
In view of adaptedness, 
we simply use $L$ for $L|_{L^q_TL^2_x}$ in the following.
We equip 
$\L_{a}(q,\infty,s_2)$ with the following Fr\'echet metric:
\begin{align}
d_\infty(L,\wt L)&=\sum_{k=1}^\infty2^{-k}
\Big(1\wedge\|L-  \wt L\|_{{\L(q,k,s_2)}}\Big).
\label{JW2}
\end{align}

\noi
Let $\mathcal M$ be the product space given by 
\begin{align}
\begin{aligned}
\mathcal M 
& =  \vec H^s(\T^2)\times 
\Ec_z\times \Ec_Z\times \Ec_R
\times \L_a(q,\infty, s_2), 
\end{aligned}
\label{A34}
\end{align}

\noi
where
\begin{align}
\begin{split}
\Ec_z
&=C(\R_+;W^{-\frac12-\eps,\infty}(\T^2))
\cap C^1(\R_+;W^{-\frac32-\eps,\infty}(\T^2)),\\
\Ec_Z
& = 
C(\R_+;W^{\frac14-\eps,\infty}(\T^2))
\cap C^1(\R_+;W^{-\frac34-\eps,\infty}(\T^2)), \\
\Ec_R
&=C(\R_+;W^{-\frac14-\eps,\infty}(\T^2)), 
\end{split}
\label{A24}
\end{align}

\noi
endowed with the compact open topology in time.
Given $k \in \N$, define
$\Ec_z(k)$, $\Ec_Z(k)$, and $\Ec_R(k)$
 as in \eqref{A24}
by replacing 
$\R_+$ in \eqref{A24} with $[0, k]$.
Then, 
$\Ec_z$, $\Ec_Z$, and $\Ec_R$ in \eqref{A24}
are endowed with 
the  following Fr\'echet metric:
\begin{align*}
d_{\Ec_A}
(f,\wt f)&=\sum_{k=1}^\infty2^{-k}
\Big(1\wedge\|f-  \wt f\|_{\Ec_A(k)}\Big), \qquad 
\text{$A = z$, $Z$, or $R$.}
\end{align*}

\noi
We note that $\Ec_z$, $\Ec_Z$, and $\Ec_R$
are separable; see \cite[Remark 5.6]{LOZ}.
We endow $\M$ with the following metric:
\begin{align}
d_{\M}(\Xi, \wt \Xi)
= \|\pmb u - \wt{\pmb u}\|_{\vec H^s}
+ d_{\Ec_z}(z, \wt z)
+ d_{\Ec_Z}(Z, \wt Z)
+ d_{\Ec_R}(R, \wt R)
+ d_\infty(L, \wt L)
\label{JW4}
\end{align}

\noi
 for $\Xi =(\pmb u, z, Z, R, L)$
and $\wt\Xi =(\wt{\pmb u}, \wt z, \wt Z, \wt R, \wt L)$.

We first show that $\If_{\pl,\pe}^N$ in 
\eqref{X6} and \eqref{POa}
is measurable with respect to the Borel
$\s$-algebra on $(\L_{a}(q,\infty, s_2), d_\infty)$.
Given 
 $N\in\N$, set 
\begin{align*}
E_{N}=C(\R_+;\pi_NL^2(\T^2)).
\end{align*}

\noi
For $f,g\in E_{N}$,
let $\If_{\pl,\pe}^{N}[f,g]$
be the operator obtained by replacing
the two entries of $z_N$ in~\eqref{POb}
with $f$ and $g$.
More precisely, set
\begin{align}
\If_{\pl,\pe}^{N}[f,g](w)(t)
&=\sum_{n \in \Z^2} e_n\int_0^t\sum_{n_1\in \Z^2}
{\ft w(t',n_1)}\A^{N}_{n,n_1}[f, g](t,t')dt',
\label{A27}
\end{align}

\noi
where the kernel $\A^{N}_{n,n_1}[f, g]$ is given by 
\begin{align}
\A^{N}_{n,n_1}[f, g](t,t')
&=\ind_{[0, t]}(t')
\sum_{\substack{n-n_1 = n_2+n_3\\
|n_1| \ll |n_2|^\ta \\|n_1+n_2| \sim |n_3|}}
\frac{\sin((t-t')\jb{n_1+n_2})}{\jb{n_1+n_2}}
{\ft f(t',n_2)}{\ft g(t,n_3)}.
\label{JW3}
\end{align}

\noi
By comparing \eqref{POa} and  \eqref{POb}
with \eqref{A27} and \eqref{JW3}, we have 
\[\If_{\pl,\pe}^N = \If_{\pl,\pe}^N[z_N,z_N]\qquad \text{and}\qquad
\A^{N}_{n,n_1} = \A^{N}_{n,n_1}[z_N, z_N].\]

\noi
Note that $\If_{\pl,\pe}^N[f,g](w)$
is continuous in time
and is adapted.
We first establish the following preliminary lemma.

\begin{lemma}
\label{LEM:A1}
Given 
$0<s_2<1$, small $\ta>0$, and  $N\in \N$, 
let 
$\If_{\pl, \pe}^N[f,g]$ be  as in \eqref{A27}.
Then, given any 
$1<q<\infty$ and  $T>0$, we have 
\begin{align}
\|\If_{\pl,\pe}^N[f,g]\|_{{\L(q,T,s_2)}}
\les T^{\frac1{q'}}N^{s_2+\ta}
\|f\|_{C_TL_x^2}\|g\|_{C_TL_x^2}
\label{A28}
\end{align}

\noi
for any $f,g\in E_{N}$.
Consequently, 
the map $\Phi_{N,T}:
C([0, T];\pi_NL^2(\T^2))
\to {\L(q,T,s_2)}$
given by\footnote{Here, 
the operator in \eqref{A29} is restricted
on $L^q([0, T];L^2(\T^2))$, 
thus depending only on 
the information of $f$  on the time interval $[0, T]$.
}
\begin{align}
\Phi_{N,T}(f)=\If_{\pl,\pe}^N[f,f]
\label{A29}
\end{align}

\noi
is Lipschitz continuous on bounded sets in the operator norm.
\end{lemma}

\begin{proof}
Let $0\le t'\le t \le T$.
From  \eqref{JW3}
 and Young's  inequality, we have 
\begin{align*}
|\A^{N}_{n,n_1}[f, g](t,t')|
&\les \|f(t')\|_{L_x^2}\|g(t)\|_{L_x^2},
\end{align*}

\noi
and thus we have
\begin{align}
\begin{split}
\big\|\jb{n}^{s_2-1}
\A^{N}_{n,n_1}[f, g](t,t')\big\|_{\l^2_{n,  n_1}}^2
&\les\sum_{|n|\les N}\jb{n}^{2s_2-2}
\sum_{|n_1|\les N^\ta}
\|f(t')\|_{L_x^2}^2\|g(t)\|_{L_x^2}^2\\
&\les N^{2s_2+2\ta}
\|f(t')\|_{L_x^2}^2\|g(t)\|_{L_x^2}^2
\end{split}
\label{A30}
\end{align}

\noi
for $f, g \in E_N$, 
where we used $s_2>0$
in the second inequality.
Hence, by 
\eqref{A27}, 
Cauchy-Schwarz's inequality in $n_1$, 
Minkowski's integral inequality,
\eqref{A30}, and 
H\"older's inequality in~$t'$,
we have 
\begin{align}
\begin{split}
\|\If_{\pl,\pe}^N[f,g](w)\|_{C_T H^{s_2-1}_x}
&\les N^{s_2+\ta}\|g\|_{C_TL_x^2}
\int_0^T\|w(t')\|_{L_x^2}\|f(t')\|_{L_x^2}dt'\\
&\le T^{\frac1{q'}}N^{s_2+\ta}
\|f\|_{C_TL_x^2}\|g\|_{C_TL_x^2}
\|w\|_{L_T^qL_x^2}.
\end{split}
\label{A39b}
\end{align}

\noi
This proves  \eqref{A28}.
The claimed local Lipschitz continuity of $\Phi_{N,T}$ in \eqref{A29}
then follows from 
the bilinearity of $\If_{\pl,\pe}^N[f,g]$
and 
\eqref{A28} with
\eqref{L2}.
\end{proof}

Next we prove  measurability of
the extended enhanced data  set 
$(u_0, u_1, z_N,  Z_N, R_N, 
\If_{\pl, \pe}^N) $
in~\eqref{data2x}.
For $\pmb a=(a_0,a_1)\in\vec H^s(\T^2)$, set\footnote{Strictly speaking, 
in the definition of 
$L_{N,k}[\pmb a]$, 
we should have $z_N[\pmb a]|_{[0, k]}$ instead of $z_N[\pmb a]$.
However, for simplicity of notation, we drop the restriction onto $[0, k]$ when there is no confusion.}
\begin{align}
\begin{split}
z_N[\pmb a](t)
&=W(t) \big(C_{\al,N}\pi_N a_0, C_{\al,N}\pi_N a_1\big),\\
Z_N[\pmb a]&=\I\big(z_N[\pmb a]^2-\s_N\big),\\
R_N[\pmb a]&=Z_N[\pmb a]\pe z_N[\pmb a],\\
L_{N,k}[\pmb a]
&=\Phi_{N,k}(z_N[\pmb a]),\quad k \in \N, 
\end{split}
\label{A38}
\end{align}

\noi
where $C_{\al, N}$, $W(t)$,  $\I$, 
$\s_N$, 
and $\Phi_{N, k}$ 
are as in 
\eqref{CN},  \eqref{Wx}, 
\eqref{Duhamel}, 
\eqref{sigN}, 
\eqref{A29}, respectively.

\begin{lemma}
\label{LEM:DataC}
Let $\al \ge \frac 12$. Let  $s  \in \R$, 
 $0<s_2<1$,  
 and 
 $1<q<\infty$.
Then, for each  $N \in \N$, 
the map $\Psi_N: \vec H^s(\T^2) \to \mathcal M$
defined by 
\begin{align}
\Psi_N (\pmb a) 
= \Big(\pmb a, z_N[\pmb a], Z_N[\pmb a], R_N[\pmb a],
\If_{\pl,\pe}^N\big[z_N[\pmb a],z_N[\pmb a]\big]\Big)
\label{HH1}
\end{align}

\noi
is continuous, where $\M$ is as in \eqref{A34}.
In particular, if $s < - \al$, 
for the Gaussian random initial data  $(u_0, u_1)$  in \eqref{data1}, 
%with $\al \ge \frac 12$, 
the extended enhanced data  set in \eqref{data2x}\textup{:}
\begin{align}
(u_0, u_1, z_N,  Z_N, R_N, 
\If_{\pl, \pe}^N) 
= \Psi_N(u_0, u_1)
\label{JWW0d}
\end{align}
is  an $\mathcal M$-valued Borel random variable
\textup{(}after modifying on a null set\textup{)}.

\end{lemma}

\begin{proof}
Since 
$\Ec_z$, $\Ec_Z$ and $\Ec_R$ in \eqref{A24}
are all equipped with the compact-open 
topology in time, we only need to show the result
on the time interval $[0,k]$ for any $k \in \N$.

Let  $\pmb a, \pmb b\in\vec H^s(\T^2)$.
Since $N \in \N$ is fixed, 
a crude estimate with \eqref{A38}, 
Bernstein's and H\"older's inequalities
yields 
\begin{align}
\begin{split}
& \big\|z_N[\pmb a]-z_N[\pmb b]\big\|_{\Ec_z(k)}
+ \big\|Z_N[\pmb a]-Z_N[\pmb b]\big\|_{\Ec_Z(k)}\\
& \quad + 
\|R_N[\pmb a]-R_N[\pmb b]\|_{\Ec_R(k)}
\les_{\al, N, k, s}
\Big(1 + \|\pmb a\|_{\vec H^s}+\|\pmb b\|_{\vec H^s}\Big)^2
\|\pmb a-\pmb b\|_{\vec H^s}.
\end{split}
\label{A39a}
\end{align}

\noi
From \eqref{A29} and \eqref{A38}, 
we have  
\[\If_{\pl,\pe}^N \big[z_N[\pmb a],z_N[\pmb a]  \big] \big|_{[0,k]}
= \Phi_{N,k}\big(z_N[\pmb a]|_{[0,k]}\big)
= L_{N,k}[\pmb a]
\]
 for every $k\in \N$.
Then, from the Lipschitz continuity of 
$\Phi_{N,k}$ (Lemma \ref{LEM:A1}; see also \eqref{A39b})
and \eqref{A39a}, we have 
\begin{align}
\begin{split}
\|L_{N,k}[\pmb a]-L_{N,k}[\pmb b]\|_{{\L(q,k, s_2)}}
&\les_{\al, N, k, s}
\Big(\|\pmb a\|_{\vec H^s}+\|\pmb b\|_{\vec H^s}\Big)
\|\pmb a-\pmb b\|_{\vec H^s}.
\end{split}
\label{A41}
\end{align}

\noi
Hence, from \eqref{A39a} and \eqref{A41}, 
we conclude that
$\Psi_N$ is 
 continuous from $\vec H^s(\T^2)$
into $\M$.
Consequently,  with  the Gaussian random initial data  $(u_0, u_1)$  in \eqref{data1}, 
the extended enhanced data 
set
$\Psi_N(u_0, u_1) = (u_0, u_1, z_N,  Z_N, R_N, 
\If_{\pl, \pe}^N) $ in \eqref{data2x}
is measurable with respect to 
the Borel $\s$-algebra on $(\mathcal M, 
d_{\mathcal M})$, 
provided that $s < -\al$.
\end{proof}

\subsection{Joint convergence in law of the extended enhanced data set}
\label{SUBSEC:Sk0}

As a preparation to 
apply
Skorokhod representation theorem (Lemma~\ref{LEM:Sk}), 
we prove 
 joint convergence in law of the extended enhanced data set 
$(u_0, u_1, z_N,  Z_N, R_N, 
\If_{\pl, \pe}^N) $
in \eqref{data2x} to 
$(u_0, u_1, 0,  Z, 0, 0)$
in~\eqref{data2y}.

\begin{lemma}
\label{LEM:SEP}
Let  $\al \ge\frac 12$, $s < - \al$, $0 < s_2 < 1$, and $1 < q < \infty$.
Let $(u_0, u_1)$
be  the Gaussian random initial data    in \eqref{data1}.
Given $N \in \N$, let $\mu_N$
denote the law of 
the extended enhanced  data set 
 $(u_0, u_1, z_N,  Z_N, R_N, 
\If_{\pl, \pe}^N) $  in \eqref{data2x}, 
viewed
as an $\mathcal M$-valued Borel 
random variable. 
Then, as $N\to\infty$, $\mu_N$ converges weakly to  
\begin{align}
\mu = \Law(u_0, u_1, 0, Z, 0 ,0)
\label{JC1a}
\end{align}
as Borel probability measures on $\M$.
Moreover,  
there exists a closed separable subset $\S
\subset \M$ such that $\mu(\S) =1$.

\end{lemma}

Note that in the statement of Lemma \ref{LEM:SEP}, 
we used Lemma \ref{LEM:DataC}
to view $\mu_N$ as a Borel probability measure on $\M$.

\begin{proof}[Proof of  Lemma \ref{LEM:SEP}]
Arguing as in   
\cite[Remark 3.6]{LLOT} with Lemma \ref{LEM:COV}, Lemma \ref{LEM:uniq2}, 
and Proposition~\ref{PROP:rand1}, 
we see that  $(u_0, u_1, Z_N)$ converges in law to 
$(u_0, u_1, Z)$ in $\vec H^s(\T^2) \times \Ec_Z$
as $N \to \infty$. 
More precisely, given 
$0\le t_1\le\cdots\le t_m$,
and $\phi_0,\phi_1,\psi_1,\dots,\psi_m\in\D(\T^2)$, 
consider
\begin{align}
Y_N = \Big(\jb{u_0,\phi_0},\jb{u_1,\phi_1}, \jb{Z_N(t_1),\psi_1},\dots,\jb{Z_N(t_m),\psi_m}\Big).
\label{JC1}
\end{align}

\noi
Thus, 
by a straightforward adaptation of the proofs of
Lemmas \ref{LEM:COV} and \ref{LEM:uniq2}, 
with the vector-valued fourth moment
theorem \cite[Lemma~2.7]{LLOT}, 
we see that $Y_N$
converges in law to the following Gaussian vector: 
\begin{align}
Y = \Big(\jb{u_0,\phi_0},\jb{u_1,\phi_1}, \jb{Z(t_1),\psi_1},\dots,\jb{Z(t_m),\psi_m}\Big).
\label{JCC1}
\end{align}

\noi
Note the first two entries in \eqref{JC1} are independent and 
are in the first order Wiener chaos~$\H_1$, while all other entries in \eqref{JC1} are in the second order (homogeneous) Wiener chaos $\H_2$, 
in particular orthogonal to the first two entries.
Then, by taking $N\to\infty$, 
we see that the covariance matrix of $Y$ has zero cross covariance between the first two entries and the remaining entries in \eqref{JCC1}.
Hence, from 
 the joint Gaussianity of $Y$, 
 we see that 
 $(u_0,u_1)$ and $Z$ are independent.

From  Lemmas~\ref{LEM:randln}, \ref{LEM:rand2} 
and \ref{LEM:rand4} with \eqref{JW1} and \eqref{JW2},
we see that 
$(z_N, R_N, \If_{\pl,\pe}^N)$ converges in probability 
to $(0, 0, 0)$ in $\Ec_z \times \Ec_R \times \L_{a}(q,\infty,s_2)$
as $N \to \infty$.
Hence, 
we 
conclude that $\mu_N$ converges weakly to $\mu$ as $N \to \infty$.
The limiting law $\mu$ in \eqref{JC1a} is concentrated on the closed separable set
\begin{align*}
\mathcal S
=\vec H^s(\T^2)\times\{0\}\times \Ec_Z
\times\{0\}\times\{0\}\subset\mathcal M,
%\label{A35}
\end{align*}

\noi
and 
we have $\mu(\mathcal S) =1$ since
$\big(\Law(u_0, u_1)\otimes \Law(Z) \big)(\vec H^s(\T^2) \times \Ec_Z) =1$.
\end{proof}

\subsection{Application of Skorokhod representation theorem}
\label{SUBSEC:Sk2}

In this subsection, we apply the Skorokhod representation theorem 
(Lemma \ref{LEM:Sk}) and show that 
\eqref{Law5}, \eqref{Law6}, \eqref{Law7},  and \eqref{Law8} hold
on a new probability space $(\wt \O, \wt \F, \wt \PP)$.

\begin{proposition}
\label{PROP:Sk}
Let $\al\ge\frac12$, $s<-\al$, 
$0<s_2<1$, and 
$1<q<\infty$, 
and let $\M$ be as in \eqref{A34}.
Then, there exist a probability space $(\wt\O,\wt\F,\wt\PP)$,
$\mathcal M$-valued Borel random variables
\begin{align}
\wt V_N=\big(\wt u_0^{(N)},\wt u_1^{(N)},
\wt z_N,\wt Z_N,\wt R_N,\wt\If_{\pl,\pe}^N\big),
\quad N\in \N
\label{JWW1}
\end{align}

\noi
and an $\vec H^s(\T^2)\times \Ec_Z$-valued random variable
$\big(\wt u_0,\wt u_1,\wt Z\big)$ such that the following hold\textup{:}

\smallskip
\noi
\textup{(i)} 
Let $(u_0,u_1,z_N,Z_N,R_N,\If_{\pl,\pe}^N)$ be as in \eqref{data2x}.
Then, 
for every $N\in\N$, we have
\begin{align}
\Law(\wt V_N)
&=\mu_N = \Law(u_0,u_1,z_N,Z_N,R_N,\If_{\pl,\pe}^N), \quad N \in \N, 
\label{JWW0ax}
\\
\Law\big((\wt u_0,\wt u_1),\wt Z\big)
&=\Law(u_0,u_1)\otimes\Law(Z), 
\label{JWW0a}
\end{align}

\noi
which also implies \eqref{Law2}\textup{;} for each $N \in \N$, we have 
\begin{align}
\begin{split}
\Law \big(\wt u_0^{(N)}, \wt u_1^{(N)}\big)
= 
\Law (\wt u_0, \wt u_1)
= \Law (u_0, u_1)
=   \mu_{1-\al}\otimes \mu_{-\al}, 
\end{split}
\label{JWW0b}
\end{align}

\noi
where $\mu_s$ is as in \eqref{mu1}.
Moreover, 
$\wt V_N$ converges to 
\begin{align*}
\wt V = (\wt u_0,\wt u_1,0,\wt Z,0,0)
\end{align*}

\noi
 in $\M$, 
$\wt \PP$-almost surely, 
as $N\to \infty$, 
and 
 the new enhanced data sets on the new probability space $\wt \O$\textup{:}
\begin{align}
\wt\Xi_N=\big(\wt z_N,\wt Z_N,\wt R_N,\wt\If_{\pl,\pe}^N\big)
\qquad \text{and}
\qquad
\wt\Xi=(0,\wt Z,0,0)
\label{JWW0c}
\end{align}

\noi
belong 
to $\Xc^{\eps,s_2}$ defined in \eqref{L3}, $\wt \PP$-almost surely,  for every $N\in\N$,
and $\wt\Xi_N$ converges to $\wt\Xi$ in $\Xc^{\eps,s_2}$, 
$\wt \PP$-almost surely, as $N \to \infty$.

\medskip
\noi
\textup{(ii)}
For every $N\in\N$,  \eqref{Law5} holds
$\wt\PP$-almost surely.

\medskip
\noi
\textup{(iii)}
Let $\wt \zeta$ be as in \eqref{Law6}.
Then, we have $\Law(\wt\ze)=\Law(\ze)$,
where $\ze$ is as in \eqref{XI1} and \eqref{nxib}.
Moreover, 
we have $(\wt Z,\dt\wt Z)|_{t=0}=(0,0)$ and $\wt Z=\I(\wt\ze)$, 
$\wt\PP$-almost surely.

\end{proposition}

\begin{proof}
\textup{(i)}
As in 
Lemma \ref{LEM:SEP},
let $\mu_N = \Law \big(
(u_0, u_1, z_N,  Z_N, R_N, 
\If_{\pl, \pe}^N) \big)$.
In view of 
Lemma \ref{LEM:SEP},
we  apply 
the Skorokhod representation theorem (Lemma~\ref{LEM:Sk}) 
to $\{\mu_N\}_{N \in \N}$ converging weakly 
to $\mu$ in \eqref{JC1a}
as Borel probability measures
on the metric space $(\M,d_{\M})$.
This yields  a new probability space $(\wt\O,\wt\F,\wt\PP)$, 
 $\mathcal M$-valued Borel random variables
$\wt V_N$ on $\wt \O$ of the form  \eqref{JWW1}, 
and 
an  $\mathcal M$-valued Borel random variable
$\wt V$ on $\wt \O$
such that
$\Law(\wt V_N)=\mu_N$ (proving \eqref{JWW0ax}), $\Law(\wt V)=\mu$, where 
 $\mu$ is as in \eqref{JC1a}, 
and 
$\wt V_N$ converges to $\wt V$, $\wt \PP$-almost surely, 
as $N \to \infty$.

From Lemma \ref{LEM:SEP}, we have 
 $\mu(\mathcal S)=1$ and thus  we have $\wt V\in\mathcal S$, 
$\wt\PP$-almost surely, 
which implies that the third, fifth, and sixth components
of $\wt V$ are $0$, $\wt\PP$-almost surely.
Hence, 
 we may write
\begin{align}
\wt V
=(\wt u_0,\wt u_1,0,\wt Z,0,0),
\label{JWW2}
\end{align}

\noi
where $(\wt u_0,\wt u_1) \in \vec H^s(\T^2)$ and $\wt Z \in \Ec_Z$ 
are Borel random variables.
From \eqref{JC1a} and \eqref{JWW2}
with $\Law( \wt V) = \mu$, 
we obtain \eqref{JWW0a} and \eqref{JWW0b}
as marginal laws.

Given  $L\in\L_{a}(q,\infty,s_2)$
defined in \eqref{JW1}, 
we have
$L\in\L_{a}(q,1,s_2)$, 
where  $\L_{a}(q,1, s_2)$ denotes  the closed subspace
of adapted operators in $\L(q,1, s_2)$
defined in~\eqref{L2}.
We also note that 
thanks to the adaptedness, 
the $\wt \PP$-almost sure convergence of $\wt\If_{\pl,\pe}^N$ to $0$
in  $\L_{a}(q,\infty,s_2)$
implies
$\wt \PP$-almost sure convergence of $\wt\If_{\pl,\pe}^N$ to $0$
in  $\L_{a}(q,1,s_2)$.
Hence, defining
$\wt\Xi_N$ and $\wt\Xi$ by \eqref{JWW0c}, 
we see that 
$\wt\Xi_N$ and $\wt\Xi$
belong to $\Xc^{\eps,s_2}$, $\wt \PP$-almost surely,  for every $N\in\N$,
and that $\wt\Xi_N$ converges to $\wt\Xi$ in $\Xc^{\eps,s_2}$, 
$\wt \PP$-almost surely, as $N \to \infty$.

\medskip
\noi
\textup{(ii)}
Fix $N\in\N$.
As seen in \eqref{JWW0d}, 
 the extended enhanced data
set in 
\eqref{data2x} is equal to $\Psi_N(u_0, u_1)$, 
where $\Psi_N$ is as in \eqref{HH1}.
Let 
$\wt \Psi_N$ be the map obtained from 
 $\Psi_N$ by projecting onto  the second, third,  fourth, and fifth components.
Then, it follows
from Lemma~\ref{LEM:DataC}
that 
$\wt \Psi_N$ is continuous.
Note that 
$\Psi_N(u_0, u_1)$  belongs almost surely to the following set:
\begin{align}
\mathcal G_N=\Psi_N(\vec H^s(\T^2))
\label{HH2}
\end{align}

\noi
which is closed, since,  by \eqref{HH1}, it is  the graph of the continuous map
 $\wt \Psi_N$.

From \eqref{zN}, \eqref{ZN}, \eqref{Zrz}, \eqref{POa}, we have
 $\mu_N(\mathcal G_N)=1$.
Thus, 
we have 
$\wt\PP(\wt V_N\in\mathcal G_N)=1$ since the law of 
$\wt V_N$ in \eqref{JWW1} under $\wt \PP$ equals $\mu_N$ by Part~(i).
Let $\mathcal G=\bigcap_{N\in \N}\{\wt V_N\in\mathcal G_N\}$.
By noting that $\wt \PP(\mathcal G) = 1$, 
it follows from  \eqref{JWW1}, 
 \eqref{HH2}, and \eqref{HH1} with \eqref{A38}
that 
\eqref{Law5} holds.

\medskip
\noi
\textup{(iii)}
The map $\Gamma_1: \Ec_Z\to \D'(\mathring{\R}_+\times\T^2)$
given by 
\begin{align*}
\Gamma_1(Z) = (\dt^2+1-\Dl)Z
\end{align*}

\noi 
and the map $\Gamma_2: \Ec_Z \to W^{\frac14-\eps,\infty}(\T^2)\times W^{-\frac34-\eps,\infty}(\T^2)$ given by 
\begin{align*}
\Gamma_2(Z) = (Z, \dt Z)|_{t=0}
\end{align*}

\noi
are both continuous, 
in particular, measurable. 
From \eqref{Z1}, we have $\G_1(Z) = \ze$
and $\G_2(Z)=(0,0)$, almost surely.
Since $\Law(\wt Z)=\Law(Z)$ by Part~(i), 
we then have 
\begin{align}
\Law(\wt\ze)
& =\Law\big(\G_1(\wt Z)\big)
=\Law\big(\G_1( Z)\big)
=\Law(\ze), \notag \\
(\wt Z,\dt\wt Z)|_{t=0} & = \G_2(\wt Z) =(0,0), \quad \text{$\wt \PP$-almost surely}.
\label{HH3}
\end{align}

\noi
Finally, 
$\wt Z=\I(\wt\ze)$
follows from \eqref{Law6} and \eqref{HH3}.
\end{proof}

We conclude this section by discussing 
the necessity of the first two components of $\wt V_N$ 
in~\eqref{JWW1} being $N$-dependent.

\begin{remark}\rm
\label{REM:IV1}
In Proposition \ref{PROP:Sk}, 
the first two components of $\wt V_N$ in \eqref{JWW1} 
need to 
be $N$-dependent.
Otherwise, $\wt Z_N$, satisfying \eqref{Law5}, 
does not converge to any limit almost surely.

Suppose,  for contradiction,  that
the first two components of $\wt V_N$ 
are given 
by  a pair 
$(v_0, v_1)\in\vec H^s(\T^2)$ for any $N \in \N$.  
Fix $0 < t < 2\pi$.
A direct computation shows 
\begin{align}
\wt \E\Big[\ft {\wt Z}_N(t,0) \ft{\wt  Z}_M(t,0)\Big]
=\frac{C_{\al,M}^2}{C_{\al,N}^2} \wt \E\Big[|\ft {\wt Z}_N(t,0)|^2\Big]
\label{A52}
\end{align}

\noi
for any $M \ge N \ge 1$.
From 
Lemma~\ref{LEM:COV} and \eqref{nxib}, we have, as $N \to \infty$, 
\begin{align}
\wt \E\Big[|\ft {\wt Z}_N(t,0)|^2\Big]\longrightarrow
\int_0^t\int_0^t
\sin(t-t_1)\sin(t-t_2)dt_1dt_2
=(1-\cos t)^2, 
\label{A53}
\end{align}

\noi
which is positive for  $0 < t< 2\pi$.
On the other hand, 
by taking 
 $M=N^2$, a direct computation with \eqref{CN} yields
\begin{align}
\frac{C_{\al,N^2}^2}{C_{\al,N}^2}
\longrightarrow\rho_\al
:=\begin{cases}
0,&\al>\frac12, \\
2^{-\frac12},&\al=\frac12, 
\end{cases}
\label{A54}
\end{align}

\noi
as $N \to \infty$.
Here, when $\al = \frac 12$, we used the following asymptotics:
\[
\sum_{|n|\le N}\frac1{\langle n\rangle^2}
=
2\pi\log N+O(1),
\]

\noi
which follows from Lemma \ref{LEM:count}\,(iii)
and summation by parts;
see also \eqref{K8} and \cite[Lemma~3.2]{HRW}.

Now, let  $K_N=\ft {\wt Z}_{N^2}(t,0)- \ft {\wt Z}_N(t,0)$.
Then, from \eqref{A52}, \eqref{A53}, and \eqref{A54}, we have 
\begin{align}
\| K_N\|_{L^2(\wt \O)}^2
\longrightarrow2(1-\rho_\al)(1-\cos t)^2>0, 
\label{A55}
\end{align}

\noi
as $N \to \infty$.
Moreover, noting that $K_N$ belongs
to the second order Wiener chaos, it follows from  the Wiener chaos estimate
and \eqref{A55} that 
\begin{align*}
\sup_{N\in \N} \|K_N\|_{L^4(\wt \O)}
\les \sup_{N\in \N} \|K_N\|_{L^2(\wt \O)}<\infty,
\end{align*}

\noi
which implies that $\{K_N^2\}_{N\in \N}$ is uniformly integrable
via Chebyshev's inequality.  
Now if $\big\{\ft {\wt Z}_N(t,0)\big\}_{N\in \N}$ were Cauchy in probability (which 
is a consequence of almost sure convergence),
then $K_N^2$ would converge to $0$ in probability.
By Vitali's convergence theorem
(see, for example, \cite[Theorem 4.6.3]{Durret}), 
we have 
$\|K_N\|_{L^2(\wt \O)}\to0$, 
which is a contradiction to  \eqref{A55}.

Note that, in the discussion above, 
we used our assumption that the first two components of 
$\wt V_N$ are given by a common pair $(v_0,v_1)$, independent of $N \in \N$, 
in deriving \eqref{A52} and hence \eqref{A55}, 
leading to a contradiction.

\end{remark}

\section{Proof of the main result}
\label{SEC:THM}

In this section, we present the proof of  Theorem \ref{THM:1}.
In Subsection \ref{SUBSEC:main1}, 
 we first establish 
 almost sure convergence on the new probability space
(Theorem \ref{THM:uNu}).
In Subsection \ref{SUBSEC:main2},
we then transfer the 
 almost sure convergence on the new probability space
 to convergence in law on the original probability space, 
 thus establishing 
Theorem \ref{THM:1}.

\subsection{Proof of Theorem \ref{THM:uNu}}
\label{SUBSEC:main1}

In this subsection, we present the proof of  Theorem~\ref{THM:uNu}.
While this result follows
from 
iterating the local-in-time approximation argument, 
we need to account for 
 the contribution
of the paracontrolled operators from the past
when we restart an iteration at a positive time;
see~\eqref{G3c}.

In the following, we work on the new probability space
$(\wt  \O, \wt   \F, \wt  \PP)$
constructed in  Proposition~\ref{PROP:Sk}.
Fix a realization in the set of full probability
where the convergence \eqref{data2zz} holds
and the limiting equation \eqref{NLW3b} 
(with $Z$ replaced by $\wt Z$) is  well-posed
on the maximal time interval $[0, \wt T_{\max})$ on the new probability space $\wt \O$.
For simplicity of notation,
we suppress the tildes 
and the dependence on the fixed $\o$ in the following
pathwise argument.

Fix a target time $0 < T < T_{\max}$.
Let  $s_1$ and $s_2$ be
as in Proposition \ref{PROP:LWP1}.
Given an interval $I \subset \R_+$,
we use $X^s(I)$ and $\vec X^{s_1, s_2}(I)$
to denote the spaces in \eqref{Xsp},
with $[0, T]$ replaced by $I$.
It follows from 
Remark \ref{REM:uniq2} 
that there exists a unique solution 
$v \in X^{s_2}(T)$ to the limiting equation \eqref{NLW3b}
(and \eqref{NLW5a} with $(X, Y) = (0, v)$)
on the time interval $[0, T]$.
Set
\begin{align}
K
= 1 + \|v\|_{X^{s_2}(T)}
+ \|Z\|_{C_T W_x^{\frac14-\eps, \infty}}
< \infty.
\label{G1}
\end{align}

\noi
From \eqref{data2zz},
\eqref{JW2}, and \eqref{A24}, we have
\begin{align}
\begin{split}
\dl_N
: \!& = 
\|z_N\|_{\Ec_z(T)}
+ \|Z_N - Z\|_{\Ec_Z(T)}\\
&\quad + \|R_N\|_{\Ec_R(T)}
%\|z_N\|_{C_T W_x^{-\frac12-\eps, \infty}}
%+ \|Z_N - Z\|_{C_T W_x^{\frac14-\eps, \infty}}\\
%&+ \|R_N\|_{C_T W_x^{-\frac14-\eps, \infty}}
 + \|\If_{\pl, \pe}^N\|_{\L(q, T, s_2)}
\too 0, 
\end{split}
\label{G2}
\end{align}

\noi
as $N \to   \infty$, 
where  
$\Ec_z(T)$, $\Ec_Z(T)$, and $\Ec_R(T)$
are defined as in \eqref{A24}
by replacing $\R_+$ with $[0, T]$.

In particular, we may assume $\dl_N \le 1$
by taking $N$ sufficiently large.
In the following, we restrict our attention to $N \gg1 $
such that $\dl_N \le 1$ holds.
Whenever 
the solution $(X_N, Y_N)$ to 
the paracontrolled system \eqref{NLW5} exists on $[0, t]$,
write
\begin{align}
E_N(t)
= \|(X_N, Y_N - v)\|_{\vec X^{s_1, s_2}(t)}.
\label{G3}
\end{align}

\noi
Note that 
we have $E_N(0) = 0$.
Moreover, from 
 Proposition~\ref{PROP:LWP1} and \eqref{G2}, we have 
\begin{align}
\lim_{N \to \infty} E_N(T_*) =  0, 
\label{G3a}
\end{align} 

\noi
where $T_* = T_* \big(\| \Xi\|_{\Xc^{\eps, s_2}}\big)$
with $\Xi = (0, Z, 0, 0)$.
Our primary goal in the following is to show that 
\begin{align}
\lim_{N \to \infty} E_N(T) =  0.
\label{G3b}
\end{align}

Suppose that $(X_N, Y_N)$ has been constructed
on $[0, t_0]$ for some  $ t_0  \in [0,  T)$
and that 
\begin{align}
E_N(t_0) \le 1.
\label{G4}
\end{align}

\noi
Set $w_N = X_N + Y_N + Z_N$.
In considering the continuation $(X_N, Y_N)$ on 
$I = [t_0, t_0 + \tau ]$ (for small $\tau > 0$), 
we need to write the paracontrolled terms in \eqref{NLW5}
as
\begin{align}
\begin{split}
& 4 \If_{\pl}^{N, (1)}(w_N ) \pe z_N (t) 
  + 4 \If_{\pl, \pe}^N (w_N)(t)\\
& \quad =  \Big( 4 \If_{\pl}^{N, (1)}(\ind_{[0, t_0]}\cdot w_N ) \pe z_N (t)
  + 4 \If_{\pl, \pe}^N (\ind_{[0, t_0]}\cdot w_N)(t)\Big)\\
& \quad \quad 
+   \Big( 4 \If_{\pl}^{N, (1)}(\ind_{(t_0, t_1]}\cdot w_N ) \pe z_N (t)
  + 4 \If_{\pl, \pe}^N (\ind_{(t_0, t_1]}\cdot w_N) (t)\Big)\\
& =:   \1^N_{t_0}(X_N, Y_N)(t) + \II^N_{t_0}(w_N)(t)
\end{split}
\label{G3c}
\end{align}

\noi
for $t_0 \le t \le t_1$, 
where 
$\ind_{[t_0, t_1]}\cdot w_N$
is yet to be constructed via a contraction argument.
%Denote by $\1^N_{t_0}(t)$
%the sum of the first two terms on the right-hand side
%of \eqref{G3c} which is well defined for $t_0 \le t \le T$.
Note that 
$\1^N_{t_0}(X_N, Y_N)$ is determined
by the solution $(X_N, Y_N)$ on $[0, t_0]$
and the enhanced data set $\Xi_N$ 
on $[0, T]$, defined in \eqref{data2}.
By Lemmas \ref{LEM:para} and \ref{LEM:rand3}
with 
\eqref{L2}, \eqref{G1}, \eqref{G2}, \eqref{G3}, and \eqref{G4},
we have
\begin{align}
\begin{split}
\|\1^N_{t_0}(X_N, Y_N)\|_{C([t_0, T]; H^{s_2-1}_x)}
&  \les_T
\|z_N\|_{C_T W_x^{-\frac12-\eps, \infty}}^2
\|w_N\|_{L^3_{t_0} H^{\frac14-\eps}_x}\\
&\quad
+ \|\If_{\pl, \pe}^N\|_{\L(q, T, s_2)}
\|w_N\|_{L^q_{t_0} L^2_x}\\
&   \les_T
(\dl_N^2 + \dl_N)
\big( K + E_N(t_0) \big)
\les K \dl_N
\end{split}
\label{G5}
\end{align}

\noi
for any $N \gg 1$, 
where the second step follows
from writing
$
w_N
= X_N + (Y_N-v) + v + (Z_N-Z) + Z
$
with the fact that $\dl_N \le 1 \le K$
(and we used the assumption that $E_N(t_0) \le 1$ in the third step).
Hence, under the assumption \eqref{G4}, we have 
\begin{align*}
\lim_{N \to \infty} \|\1^N_{t_0}\|_{C([t_0, T]; H^{s_2-1}_x)} = 0.
\end{align*}

Next, we extend the solution $(X_N, Y_N)$,
a priori defined on $[0, t_0]$, onto $I = [t_0, t_1]$
with $t_1 = (t_0 + \tau) \wedge T$.
Define a subclass of  
$\vec X^{s_1, s_2}(t_1)$
by setting 
\begin{align*}
\vec X^{s_1, s_2}
_{t_0, (X_N, Y_N)}(t_1)
= \big\{(\Xs, \Ys) \in  \vec X^{s_1, s_2}(t_1):
(\Xs, \Ys)|_{[0, t_0]} = (X_N, Y_N)\big\}
\end{align*}

\noi
and endow it with the metric:
\begin{align*}
d\big(  (\Xs_1,\Ys_1),
(\Xs_2,\Ys_2)  \big)
& =
\|  (\Xs_1, \Ys_1) - (\Xs_2, \Ys_2)\|_{\vec X^{s_1, s_2}(t_1)}\\
& = \|  (\Xs_1, \Ys_1) - (\Xs_2, \Ys_2)\|_{\vec X^{s_1, s_2}(I)}.
\end{align*}

\noi
Consider the Duhamel formulation of 
\eqref{NLW5}, starting at $t_0$,
with \eqref{G3c}.
More precisely, with $\ws_N = \Xs + \Ys + Z_N$, 
set 
\begin{align}
\begin{split}
\Phi_{1, t_0}^N(\Xs, \Ys) (t)
&  =  W(t-t_0)\big(X_N(t_0), \dt X_N(t_0)\big) + 2 \I_{t_0}(\ws_N\pl z_N)(t), \\
\Phi_{2, t_0}^N(\Xs, \Ys) (t)
&  =  W(t-t_0)\big(Y_N(t_0), \dt Y_N(t_0)\big)\\
& \quad 
+ \I_{t_0}\Big( \ws_N^2 + 2\ws_N \pg z_N 
  + 2 \Ys \pe z_N + 2 R_N \\
& \hphantom{XXXXl}  + \1^N_{t_0}(X_N, Y_N)
  + \II^N_{t_0}(\ws_N)\Big)(t)
\end{split}
\label{G7x}
\end{align}

\noi
for $t \ge t_0$, where $\I_{t_0}$ is as in $\I$
with the left limit of the integral replaced by $t_0$.
Then, define $\Phi^N_{t_0}$ by setting
\[
\Phi^N_{t_0}(t) =
\begin{cases}
 (X_N, Y_N)(t), & \text{for $0 \le t \le t_0$},\\
 (\Phi^N_{1, t_0}, \Phi^N_{2, t_0})(t), 
 & \text{for $t \ge t_0$}.
\end{cases}\]

\noi
From \eqref{G1}, \eqref{G3},  and \eqref{G4}, 
the $\vec X^{s_1, s_2}(I)$-norm of the linear solutions in~\eqref{G7x}
is bounded by 
a constant multiple of $K$.
Then, by applying the estimates
\eqref{A4}, \eqref{AA1}, \eqref{A8},
\eqref{A9}, \eqref{A11}, and \eqref{A12}
on $I$ together with \eqref{G4} and \eqref{G5}, we obtain
\begin{align}
\begin{split}
\|\Phi^N_{t_0}(\Xs, \Ys)\|_{\vec X^{s_1, s_2}(I)}
& \le C_0 K + C(T) \tau K^2,\\
\|\Phi^N_{t_0}(\Xs_1, \Ys_1)
- \Phi^N_{t_0}(\Xs_2, \Ys_2)\|_{\vec X^{s_1, s_2}(I)}
&  \le C(T) \tau K
\|(\Xs_1, \Ys_1) - (\Xs_2, \Ys_2)\|_{\vec X^{s_1, s_2}(I)}
\end{split}
\label{G7}
\end{align}

\noi
for any
$  (\Xs, \Ys)$, $(\Xs_j,\Ys_j)$, $j = 1, 2$, 
belonging to the closed subset:
\begin{align*}
B_K = \big\{
(\Xs,  \Ys )
\in \vec X^{s_1,s_2}_{t_0,(X_N,Y_N)}(t_1):
\|  (\Xs, \Ys)  \|_{\vec X^{s_1,s_2}(I)}
\le 2C_0K
\big\}.
\end{align*}

\noi
In the second estimate in \eqref{G7},
we used the fact that 
the terms $\1^N_{t_0}(X_N, Y_N) $ cancel
each other since 
$(\Xs_1,\Ys_1)$ and $(\Xs_2,\Ys_2)$ share the same past on $[0, t_0]$.
Note that 
the constants in~\eqref{G7}
are independent of $t_0$ and $N$
under the assumption that 
 $\dl_N \le 1$ and $E_N(t_0) \le 1$ as in~\eqref{G4}.
By choosing $0 < \tau = \tau(T, K) \le 1$ sufficiently small, 
we see that $\Phi^N_{t_0}$ is a contraction in $B_K$, 
allowing us to extend $(X_N, Y_N)$
onto $[t_0, t_1]$.
We point out that the local existence time
$\tau = \tau(T, K)$
 can be chosen
independently of $t_0$ and $N$, provided that $N \ge N_*(t_0) \gg 1$
(which is needed to guarantee~\eqref{G4}).

We now compare $(X_N, Y_N)$ with  $(0, v)$
on $I =[t_0, t_1]$ with $t_1 = (t_0 +\tau)\wedge T$.
Proceeding as in \eqref{G7} with \eqref{G3}, 
 we have
\begin{align}
\begin{split}
\|(X_N, Y_N-v)\|_{\vec X^{s_1, s_2} (I)}
&  \le
C E_N(t_0)
+ C(T) \tau K \dl_N
+ C \tau
\|\1^N_{t_0}(X_N, Y_N)\|_{C_IH^{s_2-1}_x}\\
&\quad
+ C(T) \tau K 
\|(X_N, Y_N-v)\|_{\vec X^{s_1, s_2}(I)}.
\end{split}
\label{G9}
\end{align}

\noi
From the choice of $\tau$ in the contraction argument above, 
we have $C(T) \tau K \le \frac 12$.
Hence, from~\eqref{G9} and \eqref{G5}, we obtain
\begin{align}
E_N\big((t_0 + \tau )\wedge T\big)
\les_{T, K} E_N(t_0) + \dl_N, 
\label{G10}
\end{align}

\noi
where the implicit constant is independent of $t_0$, 
small $\tau = \tau (T, K)> 0$, and $N$, provided that 
 $N \ge N_*(t_0) \gg 1$.

By choosing $\tau \le T_*$, where $T_*$ is as in 
\eqref{G3a}, it follows from \eqref{G10}, \eqref{G3a}, and \eqref{G2}
that there exists $N_2 \in \N$ such that 
\begin{align}
E_N\big((2\tau)\wedge T\big) \le 1
 \label{G10a}
\end{align}

\noi
for any $N \ge N_2$.
By iterating this argument, given $j \ge 3$, there exists $N_j \in \N$
such that 
\begin{align}
E_N\big((j\tau)\wedge T \big) \le 1
 \label{G10b}
\end{align}

\noi
for any $N \ge N_j$.
Under \eqref{G10a} and \eqref{G10b}, 
by iterating the argument above, 
we can construct the solution $(X_N, Y_N)$
to \eqref{NLW5} on the entire time interval $[0, T]$
with the bound:
\begin{align}
E_N(T)
\les_{T, K} E_N(0 ) + \dl_N = \dl_N, 
\label{G10c}
\end{align}

\noi
provided that $N \ge N(T, K)\gg1$.
Hence, \eqref{G3b} follows from \eqref{G10c} and \eqref{G2}.

Finally, from \eqref{exp5}, \eqref{exp6},
\eqref{G2}, and  \eqref{G3b}, 
we obtain
\begin{align*}
&\|(u_N, \dt u_N) - (u, \dt u)\|
_{C_T \vec H_x^{-\frac12-\eps}}\\
&\quad   \les
\|z_N\|_{X^{-\frac12-\eps}(T)}
+ \|Z_N-Z\|_{X^{-\frac12-\eps}(T)}
+ E_N(T)
\too 0, 
\end{align*}

\noi
as $N \to \infty$.
This concludes  the proof of Theorem \ref{THM:uNu}.

\subsection{Proof of Theorem \ref{THM:1}}
\label{SUBSEC:main2}

In this subsection, we briefly discuss  the proof of Theorem~\ref{THM:1}.

Let $T$ be a  random time such that  $0 < T < T_{\max}$
almost surely, 
where
$T_{\max}$ is 
the almost surely positive
 and finite
maximal existence time for \eqref{NLW1b}.
From 
Proposition~\ref{PROP:Sk}, we have
$\Law( \wt Z ) = \Law( Z )$.
After enlarging the new probability space if necessary, 
it follows from \cite[Theorem 8.17\,(i)]{Kal}
that 
there exists a random time $\wt T$ on the new probability space $\wt \O$
such that 
\begin{align}
\Law \big(  \wt T, \wt Z  \big)
= \Law \big(  T, Z  \big).
\label{H1}
\end{align}

\noi
Let $\wt T_{\max}$ denote the maximal existence time
of the solution $\wt u$ to \eqref{XNLW2},
as in Theorem~\ref{THM:uNu}.
Noting that 
the solution $\wt u$ to \eqref{XNLW2}
is  uniquely determined by $\wt Z$
(see \eqref{NLW3b} with $Z$ replaced by $\wt Z$), we have
\begin{align}
0 < \wt T < \wt T_{\max},  \quad \text{$\wt \PP$-almost surely}.
\label{H1a}
\end{align}

\noi
Given $N \in \N$, 
let $\wt T_{N, \max}$ denote the maximal existence time
of the solution $\wt u_N$ to \eqref{XNLW1}.
From the proof of Theorem \ref{THM:uNu}
presented  in Subsection~\ref{SUBSEC:main1}, 
we have 
\begin{align}
\liminf_{N \to \infty}\, \wt T_{N, \max}
\ge \wt T_{\max},  \quad \text{$\wt \PP$-almost surely}.
\label{H1b}
\end{align}

Fix an increasing sequence $\{\g_N\}_{N \in \N} \subset \R_+$
such that $\g_N \to 1$ as $N \to \infty$,
and  set
\begin{align}
\wt T_N
= \min \big(\wt T ,  \g_N   \wt T_{N, \max}\big).
\label{H2}
\end{align}

\noi
By \eqref{H1a}, \eqref{H1b}, and \eqref{H2}, we have, $\wt \PP$-almost surely, 
\begin{align}
\begin{split}
& 0 < \wt T_N< \wt T_{N, \max}, \  N \in \N,  
\ \ \text{and} \ \  \wt T_N = \wt T
\ \  \text{for all sufficiently large $N$}.
\end{split}
\label{H2a}
\end{align}

\noi
Then, it follows from Theorem~\ref{THM:uNu} and \eqref{H2a}
that 
\begin{align}
\big(  \wt T_N, 
\wt{\vec u}_N(\,\cdot\, \wedge \wt T_N)  \big)
\too
\big(  \wt T, \wt{\vec u}(\,\cdot\, \wedge \wt T)  \big)
\label{H3}
\end{align}

\noi
in $\R \times C(\R_+; \vec H^{-\frac12-\eps}(\T^2))$,
 $\wt \PP$-almost surely, 
 as $N \to \infty$, 
 where
 $\wt{\vec u}_N = ( \wt u_N, \dt \wt u_N )$
and $\wt{\vec u} = ( \wt u, \dt \wt u )$.

It remains to construct $T_N$ on the original probability space.
From \eqref{Law2}, we have
\begin{align}
\Law \big(  \pi_N u_0, \pi_N u_1  \big)
= \Law \big(  \pi_N \wt u_0^{(N)},
\pi_N \wt u_1^{(N)}  \big), \quad N \in \N.
\label{H3a}
\end{align}

\noi
In view of \eqref{Law1}, 
by applying 
\cite[Theorem~8.17\,(i)]{Kal} once again
(after enlarging the original probability space if necessary), 
we can choose 
a pair of random times
$( T_N, T^*_N )$
defined on the original probability space $\O$
such that
\begin{align}
\begin{split}
& \Law \big(  u_0, u_1, z_N, Z_N, R_N,
\If_{\pl, \pe}^N, T_N, T^*_N  \big)\\
& \quad =
\Law \big(  \wt u_0^{(N)}, \wt u_1^{(N)},
\wt z_N, \wt Z_N, \wt R_N,
\wt \If_{\pl, \pe}^N, \wt T_N, \wt T  \big).
\end{split}
\label{H4}
\end{align}

In proving  Theorem \ref{THM:1}, 
it suffices to consider the enhanced data set 
$\Xi_N$ 
in \eqref{data2} on the random time interval $[0, T^*_N+1]$
(measured in a variant of the 
$\Xc^{\eps, s_2}$-norm
in \eqref{L3} with  the time interval $[0, 1]$ replaced by $[0, T^*_N+1]$).
Given 
 $N \in \N$, 
 let $T^{(T^*_N +1)}_{N, \max}$ ($\le T^*_N+1$) be the maximal existence time, 
 restricted within $[0, T^*_N+1]$, 
 of the solution 
 $(X_N, Y_N)$ 
 to the paracontrolled system~\eqref{NLW5}
associated with the enhanced data set $\Xi_N$ 
in \eqref{data2}; 
 short time existence is guaranteed by 
Proposition~\ref{PROP:LWP1},
while a slight modification of the  argument
in Subsection \ref{SUBSEC:main1}
allows us to extend $(X_N, Y_N)$
onto the maximal time interval $[0, T^{(T^*_N+1)}_{N, \max}) \subset [0,  T^*_N+1]$.
Hence, 
\eqref{H2}, 
\eqref{H2a},  and \eqref{H4} yield
\begin{align}
\PP \big(  0 < T_N < T_{N, \max}^{(T^*_N+1)}  \big)
= \wt \PP \big(  0 < \wt T_N
< \wt T_{N, \max}  \big)
= 1.
\label{H4a}
\end{align}

\noi
In particular, $(X_N, Y_N)$ exists on $[0, T_N]$, almost surely.
By \eqref{exp1}, \eqref{exp4}, and \eqref{X2d},
the function $u_N$ defined by 
\begin{align}
u_N = z_N + Z_N + X_N + Y_N
\label{H4b}
\end{align}

\noi
satisfies \eqref{NLW1a} on $[0, T_N]$.
By the uniqueness of the solution to \eqref{NLW1a},
it follows from \eqref{H4} that
\begin{align}
\Law \big(  T_N, 
\vec u_N(\,\cdot\, \wedge T_N)  \big)
= \Law \big(  \wt T_N, 
\wt{\vec u}_N(\,\cdot\, \wedge \wt T_N)  \big).
\label{H5}
\end{align}

\noi
On the other hand, 
it follows from the deterministic local well-posedness argument for 
\eqref{NLW3b} (as in the proof of Proposition \ref{PROP:B1}\,(i)
presented in Section \ref{SEC:B})
with  \eqref{H1} that 
\begin{align}
\Law \big(  \wt T, \wt{\vec u}(\,\cdot\, \wedge \wt T)  \big)
= \Law \big(  T, \vec u(\,\cdot\, \wedge T)  \big).
\label{H6}
\end{align}

\noi
Therefore, putting \eqref{H5}, \eqref{H3}, and \eqref{H6} together, 
we conclude the desired convergence in law of 
$\big(  T_N, \vec u_N(\,\cdot\, \wedge T_N)  \big)$
 to
$\big(  T, \vec u(\,\cdot\, \wedge T)  \big)$.
This concludes the proof of Theorem~\ref{THM:1}.

\appendix

\section{On probabilistic ill-posedness}
\label{SEC:A}

In this section, we briefly go over the proof of 
Corollary~\ref{COR:2}.

\begin{proof}[Proof of Corollary~\ref{COR:2}]
Suppose that \eqref{disp1} holds 
for $\dl_N=C_{\al,N}$
on the time interval $[0, T_0]$
for an almost surely positive
random time $T_0$
(on 
the original probability
space $\O$).

Let $\wt u_N,\wt u$, and $0 < T < \wt T_{\max}$ be as in
Theorem~\ref{THM:uNu}, 
where  $\wt T_{\max}>0$ denotes the maximal existence time of 
the solution $\wt u$ to the limiting equation \eqref{XNLW2}.
We first identify the zero spatial mode of $\wt\ze$.
Since the map:
$F\longmapsto \ft F(\,\cdot\,,0)$
is continuous from
$\D'(\mathring{\R}_+\times\T^2)$ to
$\D'(\mathring{\R}_+)$,
Proposition~\ref{PROP:Sk}\,(iii) yields
\begin{align}
\Law\big(\ft{\wt\ze}(\,\cdot\,,0)\big)
=
\Law\big(\ft\ze(\,\cdot\,,0)\big)
= \Law(\ze_0), 
\label{A66a}
\end{align}

\noi
where $\ze_0$ is as in \eqref{nxib}.
Recall from the beginning of Section \ref{SEC:Psecond}
that $\ze_0$  is a real standard Gaussian
which is constant in time.
Thus, 
$\ft{\wt\ze}(\,\cdot\,,0)$ is also almost surely
a constant temporal distribution.
Namely, there exists a real-valued random variable $\wt g_0$
such that $\ft{\wt\ze}(\,\cdot\,,0) = \wt g_0$
as a temporal distribution.
Moreover, from \eqref{A66a}, 
we see that 
$\wt g_0$ is a real-valued standard Gaussian.

Since $\wt u$ is the solution to  \eqref{XNLW2},
writing  the Fourier zero mode of the Duhamel formulation gives
\begin{align}
\begin{split}
\ft{\wt u}(t,0)
&=(1-\cos t)\wt g_0
+\int_0^t\sin(t-t')
\|\wt u(t')\|_{L_x^2}^2dt'\\
&=\frac{t^2} 2\wt g_0  +o(t^2)
\end{split}
\label{A66}
\end{align}

\noi
as $t \to 0$,
where, in the second step,  we used the Taylor expansion for $\cos t$
and the fact that 
$\|\wt u(t)\|_{L_x^2}\to0$ as $t\to 0$.
Since $\wt g_0\ne0$ almost surely, 
it follows from \eqref{A66} that, for every $t>0$,
\begin{align}
\wt\PP\bigg(
\wt  T_{\max}>t,\,
\sup_{0\le r\le t}
|\ft{\wt u}(r,0)|>0
\bigg)
=
\wt\PP(\wt T_{\max}>t).
\label{A66b}
\end{align}

%\begin{align}
%\sup_{0\le r\le t}
%|\ft{\wt u}(r,0)|>0
%\label{A66b}
%\end{align}
%
%\noi
%almost surely on $\{T>t\}$ for every $t>0$.

Given $N \in \N$ and  deterministic $t>0$, define
\begin{align*}
F_N(t)
=
\begin{cases}
\displaystyle
\sup_{0\le r\le t}
|\ft u_N(r,0)|,
& \text{if the solution $u_N$ to \eqref{NLW1a} exists on  }[0,t],\\[6pt]
\infty,
& \text{otherwise}.
\end{cases}
\end{align*}

\noi
Define $\wt F_N(t)$ in a similar manner,  using 
the solution $\wt u_N$ to \eqref{XNLW1}.
By the deterministic uniqueness 
of solutions to \eqref{NLW1a} (with smooth initial data)
and 
 \eqref{Law2}, we have 
\begin{align}
\Law(F_N(t))=\Law(\wt F_N(t)).
\label{FX1}
\end{align}

\noi
On the one hand, 
\eqref{disp1} 
yields that 
\begin{align*}
F_N(t)\too0
\end{align*}

\noi
almost surely on $\{T_0>t\}$,
and thus the bounded convergence theorem yields
\begin{align}
\lim_{N \to \infty} \PP\big(T_0>t,\, F_N(t)>a\big) = 0.
\label{A66bb}
\end{align}

\noi
On the other hand, 
it follows from  Theorem~\ref{THM:uNu} 
that \begin{align}
\wt F_N(t)
\longrightarrow
\sup_{0\le r\le t}|\ft{\wt u}(r,0)|, 
\label{A66c}
\end{align}

\noi
almost surely on $\{T>t\}$, where $0 < T < \wt T_{\max}$
 is as in 
Theorem~\ref{THM:uNu}.
Thus, for every $a>0$, 
it follows from 
\eqref{A66c}, 
 Fatou's lemma, \eqref{FX1}, and \eqref{A66bb}
that 
\begin{align}
\begin{split}
\wt\PP\Big(
T>t,\,
\sup_{0\le r\le t}
|\ft{\wt u}(r,0)|>a
\Big)
&\le
\liminf_{N\to\infty}
\wt\PP(\wt F_N(t)>a)
=
\liminf_{N\to\infty}
\PP(F_N(t)>a)\\
& \le
\PP(T_0\le t)
+
\liminf_{N\to\infty}
\PP\big(T_0>t,\, F_N(t)>a\big)
\\
&\le
\PP(T_0\le t).
\end{split}
\label{A67}
\end{align}

\noi
Hence, it follows from 
\eqref{A66b} with  $T <  \wt T_{\max}$
and  letting $a\to 0$ in
\eqref{A67} that 
\begin{align}
\begin{split}
\wt\PP(T>t)
& = \wt\PP\Big(
T>t,\,
\sup_{0\le r\le t}
|\ft{\wt u}(r,0)|>0 
\Big)\\
& = \lim_{a\to 0}
\wt\PP\Big(
T>t,\,
\sup_{0\le r\le t}
|\ft{\wt u}(r,0)|>a
\Big)
\le\PP(T_0\le t).
\end{split}
\label{A68}
\end{align}

\noi
Finally, noting that  $T>0$, $\wt\PP$-almost surely, and
$T_0>0$, $\PP$-almost surely, 
by taking $t\to 0$
in~\eqref{A68}, we obtain  $1\le0$, 
which is a contradiction.
\end{proof}

\begin{ackno}\rm
The authors would like to thank Leonardo Tolomeo for helpful discussions. 
S.L.~would like to thank the
School of Mathematics at the University of Edinburgh for its hospitality 
during his visit in 2025, 
when this manuscript was prepared. 
G.L.~was supported by the NSFC (grant no.~12501181)
and by the Fundamental Research Funds for the Central Universities. 
S.L.~was supported by the DFG through 
the Hausdorff Center for Mathematics under Germany's Excellence Strategy - GZ 2047/1, 
Project-ID 390685813, and SFB 1720, Project-ID 539309657. 
T.O.~was supported by the European Research Council (grant no.~864138 ``SingStochDispDyn'')
and by the EPSRC Mathematical Sciences Small Grant
(grant no. EP/Y033507/1)
and also 
acknowledges support from  
the NSFC (grant no.~W2531005).
N.T.~was partially supported by the ANR
project Smooth ANR-22-CE40-0017. 

\end{ackno}

\end{document}